\documentclass[10pt,journal,compsoc]{IEEEtran}

\ifCLASSOPTIONcompsoc
  \usepackage[nocompress]{cite}
\else
  \usepackage{cite}
\fi
\ifCLASSINFOpdf
\else
\fi

\usepackage{amsfonts}
\usepackage{amsthm}
\usepackage{amsmath}

\usepackage{multirow}
\usepackage{multicol}

\usepackage{color}
\usepackage{algorithm}
\usepackage{algorithmic}

\newtheorem{assumption}{Assumption}
\newtheorem{theorem}{Theorem}
\newtheorem{lemma}{Lemma}

\newtheorem{definition}{Definition}
\newtheorem{remark}{Remark}

\usepackage{hyperref}

\usepackage{graphicx}
\usepackage[tight,footnotesize]{subfigure}
\usepackage{slashbox}
\usepackage{url}

\begin{document}
%
\title{ Nonconvex Zeroth-Order Stochastic ADMM \\ Methods with Lower Function \\ Query Complexity }
%
%
%
%

\author{Feihu~Huang,~Shangqian~Gao,~Jian~Pei, \emph{Fellow, IEEE}, and~Heng~Huang 
\IEEEcompsocitemizethanks{
\IEEEcompsocthanksitem Feihu Huang is with the Department of Electrical and Computer Engineering, University of Pittsburgh, USA, and is also with College of Computer Science and Technology, Nanjing University of Aeronautics and Astronautics,  MIIT Key Laboratory of Pattern Analysis and Machine Intelligence, Nanjing, China.
E-mail: huangfeihu2018@gmail.com
\IEEEcompsocthanksitem
Shangqian Gao is with the Department of Electrical and Computer Engineering, University of Pittsburgh, Pittsburgh, USA.
E-mail: shg84@pitt.edu
\IEEEcompsocthanksitem Jian Pei is with the Department of Computer Science, Duke University, Durham, USA.
E-mail: j.pei@duke.edu
\IEEEcompsocthanksitem Heng~Huang is with the Department of Computer Science, University of Maryland College Park, USA.
E-mail: heng@umd.edu
}}

%
%

\markboth{Journal of \LaTeX\ Class Files,~Vol.~14, No.~8, August~2015}%
{Shell \MakeLowercase{\textit{et al.}}: Bare Advanced Demo of IEEEtran.cls for IEEE Computer Society Journals}
%



\IEEEtitleabstractindextext{%
\begin{abstract}
Zeroth-order (a.k.a, derivative-free) methods are a class of effective optimization methods for solving complex machine learning problems, where gradients of the objective functions are not available or computationally prohibitive. Recently, although many zeroth-order methods have been developed, these approaches still have two main drawbacks:
1) high function query complexity; 2) not being well suitable for solving the problems with complex penalties and constraints.
To address these challenging drawbacks, in this paper, we propose a class of faster zeroth-order stochastic
alternating direction method of multipliers (ADMM) methods (ZO-SPIDER-ADMM) to solve
the nonconvex finite-sum problems with multiple nonsmooth penalties.
Moreover, we prove that the ZO-SPIDER-ADMM methods can achieve a lower function query complexity of
$O(nd+dn^{\frac{1}{2}}\epsilon^{-1})$ for finding an $\epsilon$-stationary point,
which improves the existing best nonconvex zeroth-order ADMM methods by a factor of $O(d^{\frac{1}{3}}n^{\frac{1}{6}})$,
 where $n$ and $d$ denote the sample size and data dimension, respectively.
 At the same time, we propose a class of faster zeroth-order online ADMM methods (ZOO-ADMM+) to solve the nonconvex online problems with multiple nonsmooth penalties.
We also prove that the proposed ZOO-ADMM+ methods achieve a lower function query complexity of $O(d\epsilon^{-\frac{3}{2}})$,
which improves the existing best result by a factor of $O(\epsilon^{-\frac{1}{2}})$.
 Extensive experimental results on the structure adversarial attack on black-box deep neural networks  demonstrate the efficiency of our new algorithms.
\end{abstract}

}

\maketitle

\IEEEdisplaynontitleabstractindextext

%
\IEEEpeerreviewmaketitle

\ifCLASSOPTIONcompsoc
\IEEEraisesectionheading{\section{Introduction}\label{sec:introduction}}
\else
\section{Introduction}
\label{sec:introduction}
\fi
Zeroth-order (\emph{a.k.a}, derivative-free, gradient-free) methods \cite{hu2016bandit,larson2019derivative} are a class of powerful optimization tools for many machine learning problems, which only need function values (not gradient) in the optimization.
For example, zeroth-order optimization methods have been applied to bandit feedback analysis \cite{agarwal2010optimal},
reinforcement learning \cite{choromanski2018structured}, and adversarial attacks on black-box deep neural
networks (DNNs) \cite{chen2017zoo,liu2018zeroth}.
Thus, recently the zeroth-order methods have been increasingly studied.
For example, \cite{Nesterov2017RandomGM} proposed
the zeroth-order gradient descent methods based on the Gaussian smoothing gradient estimator.
\cite{ghadimi2013stochastic} presented the zeroth-order stochastic gradient (ZO-SGD) methods
for stochastic optimization.
\cite{ghadimi2016mini} proposed a class of zeroth-order proximal stochastic gradient methods
for stochastic optimization with nonsmooth penalties.
Subsequently, \cite{liu2018admm,gao2018information} proposed the zeroth-order online (or
stochastic) alternating direction method of multipliers (ADMM) methods to also solve stochastic optimization with nonsmooth penalties. At the same time,
\cite{balasubramanian2018zeroth} developed a class of zeroth-order conditional gradient methods for constrained optimization.

\begin{table*}
  \centering
  \caption{ Convergence property comparison of the zeroth-order ADMM algorithms for finding an $\epsilon$-stationary point. C, NC, S, NS and mNS are the abbreviations of convex,
    non-convex, smooth, non-smooth and the sum of multiple non-smooth functions, respectively.
   $d$ is the dimension of data and $n$ denotes the sample size. GauGE, UniGE and CooGE are abbreviations of Gaussian distribution, Uniform distribution and Coordinate-wise smoothing gradient estimators, respectively. }
  \label{tab:1}
  \resizebox{1\textwidth}{!}{
  \begin{tabular}{c|c|c|c|c|c}
  \hline
 \textbf{Type} & \textbf{Algorithm} & \textbf{Gradient Estimator} & \textbf{Reference} &  \textbf{Problem}  & \textbf{ Function Query Complexity} \\ \hline
  \multirow{4}*{Finite-sum}
  & ZO-SVRG-ADMM & CooGE & \multirow{2}*{\cite{Huang2019zeroth}}  & \multirow{2}*{NC(S) + C(mNS)} & $O(nd + d^2n^{\frac{2}{3}}\epsilon^{-1})$  \\ \cline{2-3} \cline{6-6}
  & ZO-SAGA-ADMM & CooGE & & &  $O(nd + d^{\frac{4}{3}}n^{\frac{2}{3}}\epsilon^{-1})$ \\ \cline{2-6}
  & ZO-SPIDER-ADMM & CooGE & \multirow{2}*{Ours}  & \multirow{2}*{NC(S) + C(mNS)}
   &  \multirow{2}*{ {\color{red}{ $O( nd + dn^{\frac{1}{2}}\epsilon^{-1})$ }}} \\ \cline{2-3}
   & ZO-SPIDER-ADMM & CooGE+UniGE &   &  &  \\ \hline
  \multirow{4}*{Online} & ZOO-ADMM & GauGE & \cite{liu2018admm} &  \multirow{2}*{C(S) + C(NS)} & \multirow{2}*{$O(d\epsilon^{-2})$ } \\ \cline{2-4}
  &ZO-GADM & UniGE & \cite{gao2018information}  & &    \\ \cline{2-6}
  &ZOO-ADMM+ & CooGE & \multirow{2}*{Ours}  & \multirow{2}*{NC(S) + C(mNS)} &   \multirow{2}*{ {\color{red}{ $O(d\epsilon^{-\frac{3}{2}})$ }}} \\ \cline{2-3}
  &ZOO-ADMM+ & CooGE+UniGE &   &  &  \\ \hline
  \end{tabular}
  }
\end{table*}

\begin{figure}[!t]
  \centering
  \includegraphics[width=0.5\textwidth]{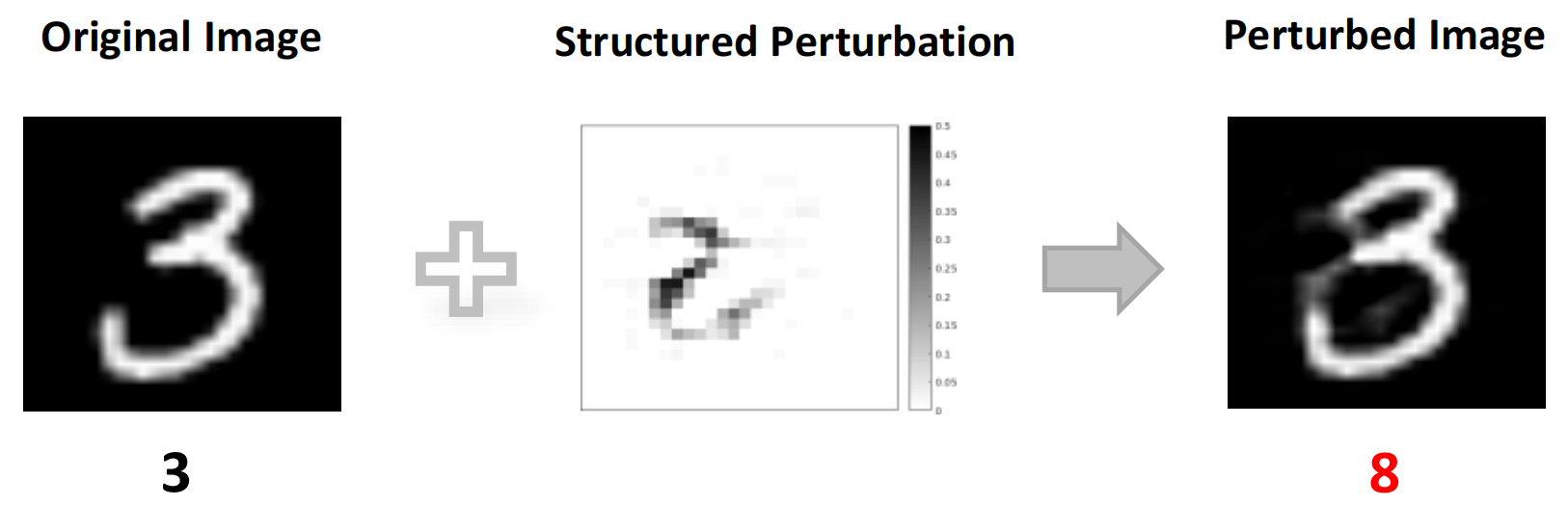}
  \caption{ Structured perturbation example not only fools the DNN, but also has some explicable.
  Black and red labels denote the initial label and the label after attack, respectively. }
  \label{fig:1}
\end{figure}

These zeroth-order methods frequently suffer from slow convergence rate and high function query complexity due to high variances generated from estimating zeroth-order (stochastic) gradients.
To alleviate this issue, some accelerated zeroth-order methods have been developed.
For example,
to accelerate the ZO-SGD,  \cite{liu2018zeroth,Liu2018StochasticZO} proposed a class of faster zeroth-order stochastic variance-reduced gradient (ZO-SVRG)
methods by using the SVRG \cite{johnson2013accelerating,allen2016variance,reddi2016stochastic}.
To reduce function query complexity, the SPIDER-SZO \cite{fang2018spider} and ZO-SPIDER-Coord \cite{Ji2019improved}
have been proposed by using stochastic path-integrated differential estimator, \emph{i.e.,} SARAH/SPIDER \cite{nguyen2017sarah,fang2018spider,wang2019spiderboost}.
For big data optimization, asynchronous parallel zeroth-order methods \cite{lian2016comprehensive,gu2018faster} and distributed zeroth-order methods \cite{hajinezhad2017zeroth} have been developed.
 More recently, \cite{huang2022accelerated} proposed a fast zeroth-order momentum method with lower query complexity.

To simultaneously deal with the non-convex loss and non-smooth regularization, \cite{ghadimi2016mini,Huang2019faster} proposed
some zeroth-order proximal stochastic gradient methods.
However, these nonconvex zeroth-order methods still are not competent for many machine learning problems with complex
nonsmooth penalties and constraints,
such as the structured adversarial attack to (black-box) DNNs in \cite{xu2018structured,Huang2019zeroth} (see Fig. \ref{fig:1}).
More recently, thus, \cite{Huang2019zeroth} proposed a class of nonconvex zeroth-order stochastic
ADMM methods (\emph{i.e.}, ZO-SVRG-ADMM and ZO-SAGA-ADMM)
to solve these complex problems.
However, these zeroth-order ADMM methods suffer from high function query complexity (please see in Table \ref{tab:1}).

In this paper, thus, we propose a class of faster zeroth-order stochastic ADMM methods with lower function query complexity
to solve the following nonconvex nonsmooth problem:
\begin{align} \label{eq:1}
\min_{x \in \mathbb{R}^d,\{y_j\in \mathbb{R}^p\}_{j=1}^m} & \quad f(x) + \sum_{j=1}^m \psi_j(y_j)  \\
\mbox{s.t.} & \quad Ax + \sum_{j=1}^mB_jy_j =c,\nonumber
\end{align}
\begin{align}
 f(x):=\left\{
\begin{aligned}
 & \frac{1}{n}\sum_{i=1}^n f_i(x)  \quad \mbox{(finite-sum)}\\
 & \mathbb{E}_{\xi} [f(x;\xi)]  \quad \mbox{(online)}
 \end{aligned} \right. \nonumber
\end{align}
where $A\in \mathbb{R}^{l \times d}$, $B_j\in \mathbb{R}^{l\times p}$ for $j\in [m]$, $c\in \mathbb{R}^l$, and $\xi$ denotes a random variable that following an unknown distribution. The  equality constraint in the problem \eqref{eq:1} encodes structure pattern of model parameters.
Here $f(x): \mathbb{R}^d\rightarrow \mathbb{R}$ is a \emph{nonconvex} and smooth
function, and $\psi_j(y_j): \mathbb{R}^p \rightarrow \mathbb{R}$ is a convex and possibly \emph{nonsmooth} function for all $j\in [m],\ m \geq 1$.
Here the explicit gradients of $f_i(x)$ and $f(x;\xi)$ are difficult or infeasible to obtain.
Under this case, we need to use the zeroth-order gradient estimators \cite{Nesterov2017RandomGM,liu2018zeroth} to estimate gradients of $f_i(x)$ and $f(x;\xi)$.
For the problem \eqref{eq:1}, its \emph{finite-sum} subproblem generally comes from the empirical loss minimization in machine learning,
and its \emph{online} subproblem generates from the expected loss minimization.
Here the non-smooth regularization functions $\{\psi_j(y_j)\}_{j=1}^m$ can encode some complex superposition structures
such as sparse and low-rank structure in robust principal component analysis (PCA) \cite{wright2009robust,candes2011robust}, and group-wise sparsity and within group sparsity structure in sparse group lasso \cite{simon2013sparse}.

In fact, the problem \eqref{eq:1} includes many machine learning problems.
For example, when $m=1$, $A=I_d$, $B_1=-I_d$ and $c=0$,
the problem \eqref{eq:1} will reduce to the standard \emph{regularized} risk minimization problem such as training deep neural networks (DNNs), defined as
\begin{align} \label{eq:2}
\min_{x\in \mathbb{R}^d}  f(x) + \psi_1(x),
\end{align}
where $f(x)$ is a nonconvex loss function, and $\psi_1(x)$ is a regularization such as $\psi_1(x)=\tau \|x\|^2$ with $\tau>0$.
When $m=1$, $A\neq I_d$ denotes a sparse graph correlation matrix,
$B_1=-I_d$ and $c=0$, i.e., the equality constraint $Ax=y_1$, and $\psi_1(y_1)=\tau\|y_1\|_1$, the problem \eqref{eq:1} will degenerate to the risk minimization problem with graph-guided fused lasso \cite{kim2009multivariate} regularization, defined as
\begin{align}
    \min_{x\in \mathbb{R}^d}  f(x) + \tau\|Ax\|_1,
\end{align}
where $f(x)$ is a loss function, and $\tau>0$. Because the term  $\|Ax\|_1$ consider the sparse relationships between different features, it can improve the performances of tasks such as SVM \cite{ouyang2013stochastic}.
When $m=2$, $A=I_{d_1}$, $B_1=-I_{d_1}$, $B_2=-I_{d_1}$ and $c=0$, the problem \eqref{eq:1} will reduce to a risk minimization problem with \emph{two}
penalties such as sparse and low-rank structure in robust PCA \cite{wright2009robust,candes2011robust}, defined as
\begin{align}
\min_{X, Y_1, Y_2} \ &   f(X) + \tau_1\|Y_1\|_1 + \tau_2\|Y_2\|_*,    \\
\mbox{s.t.} \ & X = Y_1+Y_2,  \nonumber
\end{align}
where $X, Y_1, Y_2 \in \mathbb{R}^{d_1\times d_2}$, and $\tau_1,\tau_2 >0$ denote the tuning parameters. Here we can use some nonconvex functions $f(X)$ such as $\|X-D\|_{\gamma} \ (0<\gamma<1)$ instead of $\|X-D\|_2^2$ used in \cite{wright2009robust,candes2011robust}, where $D \in \mathbb{R}^{d_1\times d_2}$ is data matrix.

It is worth highlighting that our \textbf{main contributions}
are four-fold:
\begin{itemize}
\item[1)] We propose a class of new faster zeroth-order stochastic ADMM (\emph{i.e.}, ZO-SPIDER-ADMM) methods to solve the finite-sum problem \eqref{eq:1}, based on variance reduced technique of SPIDER/SARAH and different zeroth-order gradient estimators.
\item[2)] We prove that the ZO-SPIDER-ADMM methods reach a lower function query complexity of $O(dn + dn^{\frac{1}{2}}\epsilon^{-1})$
          for finding an $\epsilon$-stationary point, which improves the existing best nonconvex zeroth-order ADMM methods by a factor of $O(d^{\frac{1}{3}}n^{\frac{1}{6}})$.
\item[3)] We extend the ZO-SPIDER-ADMM methods to the online setting, and propose a class of faster zeroth-order online ADMM methods
          (\emph{i.e.}, ZOO-ADMM+) to solve the online problem \eqref{eq:1}.
\item[4)] We provide the theoretical analysis on the convergence properties of our ZOO-ADMM+ methods, and
          prove that they achieve a lower function query complexity of $O(d\epsilon^{-\frac{3}{2}})$,
          which improves the existing best result by a factor of $O(\epsilon^{-\frac{1}{2}})$.
\end{itemize}

\noindent\textbf{Notations.}
To make the paper easier to follow, we give the following notations:
\begin{itemize}
\item $[m] = \{1,2,\cdots,m\}$ and $[j:m] = \{j,j+1,\cdots,m\}$ for all $1\leq j \leq m$.
\item $\|\cdot\|$ denotes the vector $\ell_2$ norm and the matrix spectral norm, respectively.
\item $\|x\|_G=\sqrt{x^TGx}$, where $G$ is a positive definite matrix.
\item $\sigma^A_{\min}$ and $\sigma^A_{\max}$ denote the minimum and maximum eigenvalues of $A^TA$, respectively.
\item $\sigma^{B_j}_{\max}$ denotes the maximum eigenvalues of $B_j^TB_j$ for all $j\in[m]$, and $\sigma^{B}_{\max} = \max_{1\leq j \leq m} \sigma^{B_j}_{\max}$.
\item $\sigma_{\min}(G)$ and $\sigma_{\max}(G)$ denote the minimum and maximum eigenvalues of matrix $G$, respectively;
the conditional number $\kappa_G = \frac{\sigma_{\max}(G)}{\sigma_{\min}(G)}$.
\item $\sigma_{\min}(H_j)$ and $\sigma_{\max}(H_j)$ denote the minimum and maximum eigenvalues of matrix $H_j$ for all $j\in[m]$, respectively;
$\sigma_{\min}(H) = \min_{1\leq j\leq m}\sigma_{\min}(H_j)$ and $\sigma_{\max}(H) = \max_{1\leq j\leq m}\sigma_{\max}(H_j)$.
\item $\mu, \ \nu$ denote the smoothing parameters of zeroth-order gradient estimators.
\item $\eta$ denotes the step size of updating variable $x$.
\item $L$ denotes the Lipschitz constant of $\nabla f(x)$.
\item $b, \ b_1, \ b_2$ denote the mini-batch sizes of stochastic zeroth-order gradients.
\end{itemize}
\section{Related Work}
In the section, we overview some representative existing ADMM methods and zeroth-order methods, respectively.
\subsection{ ADMM Methods }
ADMM \cite{gabay1976dual,boyd2011distributed} is a popular optimization method
for solving the composite and constrained problems.
For example, due to the flexibility in splitting the objective function into loss and complex penalty,
the ADMM can relatively easily solve some problems with complicated structure penalty
such as the graph-guided fused lasso \cite{kim2009multivariate},
which are too complicated for the other popular optimization methods such as proximal gradient methods \cite{beck2009fast}.
Thus, recently many ADMM methods \cite{he20121,he2016convergence,nishihara2015general,xu2017admm,lu2017unified} have been proposed and studied. For example, the classic convergence analysis of ADMM has been studied in \cite{he20121}.
Subsequently, some stochastic (or online) ADMM methods \cite{suzuki2013dual,ouyang2013stochastic,suzuki2014stochastic,zheng2016fast,liu2017accelerated,liu2020accelerated} have been proposed by visiting only one or mini-batch samples instead of all samples.
In fact, the ADMM method is also successful in solving many nonconvex machine learning problems such as training neural networks \cite{Taylor2016Training}.
Thus, the nonconvex ADMM methods \cite{wang2015convergence,wang2015global,hong2016convergence,jiang2019structured} have been  widely studied. Subsequently, the nonconvex stochastic ADMM methods \cite{huang2016stochastic,huang2019fasters} have been proposed. In particular,  \cite{huang2019fasters} first studied the gradient (or sample) complexities of the nonconvex stochastic ADMM methods.

\subsection{ Zeroth-Order Methods }
Zeroth-order methods have been increasingly  attracted due to effectively solve some complex machine learning problems, whose the explicit gradients are difficult or even infeasible to access.
For example, \cite{ghadimi2013stochastic,duchi2015optimal,Nesterov2017RandomGM} proposed several zeroth-order algorithms (\emph{i.e.}, ZO-SGD and ZO-SCD)
based on the Gaussian smoothing technique.
Subsequently, some accelerated zeroth-order stochastic methods (\emph{i.e.}, ZO-SVRG and ZO-SPIDER)\cite{liu2018zeroth,fang2018spider,Ji2019improved} have been proposed by using the variance reduced techniques of SVRG and SPIDER, respectively.
To solve the nonsmooth optimization, several zeroth-order proximal algorithms \cite{ghadimi2016mini,Huang2019faster} have been proposed.
At the same time, the zeroth-order ADMM-type algorithms \cite{gao2018information,liu2018admm,Huang2019zeroth}
also have been proposed.

\section{ Faster Zeroth-Order ADMMs }
In this section, we propose a class of faster zeroth-order stochastic ADMM methods (\emph{i.e.,} ZO-SPIDER-ADMM with different variants)
to solve the non-convex finite-sum problem \eqref{eq:1}.
At the same time, we extend the ZO-SPIDER-ADMM methods to the online setting and
propose a faster zeroth-order online ADMM methods
(\emph{i.e.}, ZOO-ADMM+) to solve the online problem \eqref{eq:1}.

\subsection{Preliminary}
 We first introduce an augmented Lagrangian function of the problem \eqref{eq:1} as follows:
\begin{align} \label{eq:5}
 \mathcal {L}_{\rho}(x,y_{[m]},\lambda) = & \bigg\{ f(x) \!+\! \sum_{j=1}^m \psi_j(y_j) \!-\! \langle \lambda,
 Ax + \sum_{j=1}^m B_jy_j-c\rangle \nonumber \\
 & + \frac{\rho}{2} \|Ax+ \sum_{j=1}^mB_jy_j-c\|^2 \bigg\},
\end{align}
where $f(x)=\frac{1}{n}\sum_{i=1}^nf_i(x)$ or $f(x)=\mathbb{E}[f(x;\xi)]$, and $\lambda \in \mathbb{R}^l$ denotes the dual variable, and $\rho >0$ denotes the penalty parameter.

In the problem \eqref{eq:1}, the explicit expression of gradients for $f_i(x)$ and $f(x;\xi_i)$ for all $i$ are not available, and
only its function values are available.
We can use the coordinate smoothing gradient estimator (CooGE)
\cite{liu2018zeroth,Ji2019improved} to evaluate gradient:
\begin{align}
 & \hat{\nabla}_{\texttt{coo}} f_i(x) = \sum_{j=1}^d \frac{f_i(x + \mu e_j) - f_i(x - \mu e_j)}{2\mu}e_j, \label{eq:6a} \\
 & \hat{\nabla}_{\texttt{coo}} f(x;\xi_i) = \sum_{j=1}^d \frac{f(x + \mu e_j;\xi_i) - f(x - \mu e_j;\xi_i)}{2\mu}e_j \label{eq:6b}
\end{align}
where $\mu>0$ is a coordinate-wise smoothing parameter, and $e_j$ is
a standard basis vector with 1 at its $j$-th coordinate, and 0 otherwise.
Meanwhile, we can also apply the uniform smoothing gradient estimator (UniGE)\cite{liu2018zeroth,Ji2019improved}
to evaluate gradient:
\begin{align}
 & \hat{\nabla}_{\texttt{uni}} f_i(x) = \frac{ d(f_i(x + \nu u) - f_i(x))}{\nu}u, \label{eq:7a} \\
 & \hat{\nabla}_{\texttt{uni}} f(x;\xi_i) = \frac{ d(f(x + \nu u;\xi_i) - f(x;\xi_i))}{\nu}u, \label{eq:7b}
\end{align}
where $\nu>0$ is a smoothing parameter and $u \in \mathbb{R}^d$ is a vector generated from the uniform distribution over the unit sphere. In addition, we define $f_{\nu}(x)=\mathbb{E}_{u\sim U_B}[f(x+\nu u)]$
be a smooth approximation of $f(x)$, where $U_B$ is the uniform distribution over $d$-dimensional unit Euclidean sphere $B$.

\subsection{ ZO-SPIDER-ADMM Algorithms }
In this subsection, we propose a class of faster zeroth-order SPIDER-ADMM methods (\emph{i.e.}, ZO-SPIDER-ADMM)
to solve the finite-sum problem \eqref{eq:1}. The naive extension
of the multi-block ADMM method may diverge \cite{chen2016direct},
however, our method use the linearized technique as in \cite{jiang2019structured} to guarantee its convergence.
We first propose the ZO-SPIDER-ADMM (CooGE) algorithm based on the CooGE gradient estimator and variance reduced technique of SPIDER/SARAH, which is described in Algorithm \ref{alg:1}.

\begin{algorithm}[htb]
   \caption{ ZO-SPIDER-ADMM (CooGE) Algorithm }
   \label{alg:1}
\begin{algorithmic}[1]
   \STATE {\bfseries Input:} Total iteration $K$, mini-batch size $b$, epoch size $q$, penalty parameter $\rho$ and step size $\eta$;
    \STATE {\bfseries Initialize:} $x_0 \in \mathbb{R}^d$, $y^0_j \in \mathbb{R}^p, \ j\in [m]$ and $\lambda_0 \in \mathbb{R}^l$;
   \FOR{$k=0,1,\cdots,K-1$}
   \IF{$\mod(k,q)=0$}
   \STATE{} Compute $v_k = \frac{1}{n} \sum_{i=1}^n \hat{\nabla}_{\texttt{coo}} f_i(x_k)$;
   \ELSE{}
   \STATE{} Uniformly randomly pick a mini-batch $\mathcal{S}\ (|\mathcal{S}|=b)$ from $\{1,2,\cdots,n\}$ with replacement, \\
            then update $v_k = \frac{1}{b} \sum_{i\in \mathcal{S}} \big( \hat{\nabla}_{\texttt{coo}} f_i(x_k) - \hat{\nabla}_{\texttt{coo}} f_i(x_{k-1}) \big) + v_{k-1}$;
   \ENDIF{}
   \STATE{}  $ y^{k+1}_j= \arg\min_{y_j\in \mathbb{R}^p} \tilde{\mathcal {L}}_{\rho,j} (x_k,y^{k+1}_{[j-1]},y_j,y^{k}_{[j+1:m]},\lambda_k)$ defined in \eqref{eq:y8}, for all $j\in [m]$;
   \STATE{}  $ x_{k+1}= \arg\min_{x\in \mathbb{R}^d} \hat{\mathcal {L}}_{\rho}\big( x,y^{k+1}_{[m]},\lambda_k,v_k \big)$ defined in \eqref{eq:10};
   \STATE{}  $ \lambda_{k+1} = \lambda_k- \rho(Ax_{k+1} + \sum_{j=1}^mB_jy^{k+1}_j - c)$;
   \ENDFOR
   \STATE {\bfseries Output \ (in theory):} $\{x_\zeta,y_{[m]}^{\zeta},\lambda_\zeta\}$ chosen uniformly randomly from $\{x_{k},y_{[m]}^{k},\lambda_k\}_{k=0}^{K-1}$.
   \STATE {\bfseries Output \ (in practice):} $\{x_{K},y_{[m]}^{K},z_K\}$.
\end{algorithmic}
\end{algorithm}

At the step 9 of Algorithm \ref{alg:1}, we update the variables $\{y_j\}_{j=1}^m$ by solving the following sub-problem, for all $j\in [m]$
\begin{align}
& y^{k+1}_j = \arg\min_{y_j\in \mathbb{R}^p} \tilde{\mathcal {L}}_{\rho,j}(x_k,y^{k+1}_{[j-1]},y_j,y^{k}_{[j+1:m]},\lambda_k),   \label{eq:y8} \\
& \tilde{\mathcal {L}}_{\rho,j} (x_k,y^{k+1}_{[j-1]},y_j,y^{k}_{[j+1:m]},\lambda_k) = \bigg\{ f(x_k) + \sum_{i=1}^{j-1} \psi_i(y^{k+1}_{i}) \nonumber  \\
& \quad + \psi_j(y_{j}) \!+\! \sum_{i=j+1}^m \psi_j(y^{k}_{j}) \!-\! \lambda_k^T( B_jy_j +\tilde{c}_j) \!+\! \frac{\rho}{2}\| B_jy_j +\tilde{c} \|^2 \nonumber  \\
& \quad + \frac{1}{2}\|y_j\!-\!y^k_j\|_{H_j}^2 \bigg\}, \nonumber
\end{align}
where $\tilde{c}_j = Ax_k+\sum_{i=1}^{j-1}B_iy_i^{k+1}+ \sum_{i=j+1}^{m}B_iy_i^{k}-c$ and $H_j\succ 0$. Here $\tilde{\mathcal {L}}_{\rho,j}(x_k,y^{k+1}_{[j-1]},y_j,y^{k}_{[j+1:m]},\lambda_k)$ is an approximated function of $\mathcal {L}_{\rho}(x,y_{[m]},\lambda)$ over $x_k$, $y^{k+1}_{[j-1]}$ and  $y^{k}_{[j:m]}$, where the term $\frac{1}{2}\|y_j-y^k_j\|^2_{H_j}$ ensures that the solution $y^{k+1}_j$ is not far from $y^k_j$.
When set $H_j = r_j I_p - \rho B_j^TB_j\succeq I_p$ with $r_j \geq \rho
\sigma_{\max}(B^T_j B_j) + 1$ for all $j\in [m]$ to linearize the term $\frac{\rho}{2}\|B_jy_j + \tilde{c}\|^2$,
then we can use the following proximal operator to update $y_j$, for all $j\in [m]$
\begin{align} \label{eq:s11}
 y^{k+1}_j = \mathop{\arg\min}_{y_j\in \mathbb{R}^p} \frac{1}{2}\|y_j-w^k_j\|^2 +  \frac{1}{r_j}\psi_j(y_j),
\end{align}
where $w^k_j=\frac{1}{r_j}\big( H_jy^k_j - \rho B^T_j\tilde{c} + B^T_j\lambda_k\big)$. For example,  when $\psi_j(y_j)=\|y_j\|_1$, this  subproblem \eqref{eq:s11} has a closed-form solution $y^{k+1}_j=\mbox{sign}(w^k_j)\max\big(|w^k_j|-\frac{1}{r_j},0\big)$.

At the step 10 of Algorithm \ref{alg:1}, we update the variable $x$ by solving the following problem:
\begin{align}
 & x_{k+1}= \arg\min_{x\in \mathbb{R}^d} \hat{\mathcal {L}}_{\rho}\big( x,y^{k+1}_{[m]},\lambda_k,v_k \big)  \label{eq:10} \\
 &\hat{\mathcal {L}}_{\rho} (x,y^{k+1}_{[m]}, \lambda_k, v_k) \!=\! \bigg\{ f(x_k) \!+\! v_k^T(x-x_k) \!+\! \frac{1}{2\eta}\|x-x_k\|^2_G \nonumber \\
 &  \quad + \sum_{j=1}^m \psi_j(y^{k+1}_{j})- \lambda_k^T(Ax+ \sum_{j=1}^m B_jy^{k+1}_j-c) \nonumber \\
 &  \quad + \frac{\rho}{2}\|Ax+\sum_{j=1}^mB_jy_j^{k+1}-c\|^2 \bigg\}, \nonumber
\end{align}
where  $\eta>0$ is a step size and $G\succ 0$. Here $\hat{\mathcal {L}}_{\rho} (x,y^{k+1}_{[m]}, \lambda_k, v_k)$ is an approximated function of $\mathcal {L}_{\rho}(x,y_{[m]},\lambda)$ over $x_k$ and $\{y^{k+1}_j\}_{j=1}^m$, where the term $\frac{1}{2\eta}\|x-x_k\|^2_G$ ensures that the solution $x_{k+1}$
is not far from $x_k$.
By solving the problem \eqref{eq:10},
we can obtain:
\begin{align}
 x_{k+1} = \tilde{A}^{-1}\bigg( \frac{G}{\eta} x_k - v_k - \rho A^T(\sum_{j=1}^m B_jy_j^{k+1}-c-\frac{\lambda_k}{\rho}) \bigg), \nonumber
\end{align}
where $\tilde{A} = \frac{G}{\eta} + \rho A^TA$.
To avoid computing inverse of matrix $\tilde{A}$,
we can set $G = r I_d - \rho \eta A^TA \succeq I_d$ with $r \geq \rho \eta
\sigma^A_{\max} + 1 $ to linearize term $\frac{\rho}{2}
\|Ax+ \sum_{j=1}^mB_jy^{k+1}_{j}-c\|^2$. Then we have
\begin{align}
 x_{k+1} = \frac{Gx_k}{r} \!-\! \frac{\eta v_k}{r} \!-\! \frac{\eta\rho}{r}A^T\big( \sum_{j=1}^m B_jy_j^{k+1}\!-\!c\!-\!\frac{\lambda_k}{\rho} \big).
\end{align}

\begin{algorithm}[htb]
   \caption{ ZO-SPIDER-ADMM (CooGE+UniGE) Algorithm }
   \label{alg:2}
\begin{algorithmic}[1]
   \STATE {\bfseries Input:} Total iteration $K$, mini-batch size $b$, epoch size $q$, penalty parameter $\rho$ and step size $\eta$;
    \STATE {\bfseries Initialize:} $x_0 \in \mathbb{R}^d$, $y^0_j \in \mathbb{R}^p, \ j\in [m]$ and $\lambda_0 \in \mathbb{R}^l$;
   \FOR{$k=0,1,\cdots,K-1$}
   \IF{$\mod(k,q)=0$}
   \STATE{} Compute $v_k = \frac{1}{n} \sum_{i=1}^n \hat{\nabla}_{\texttt{coo}} f_i(x_k)$;
   \ELSE{}
   \STATE{} Uniformly randomly pick a mini-batch $\mathcal{S}\ (|\mathcal{S}|=b)$ from $\{1,2,\cdots,n\}$ with replacement, and draw i.i.d. $\{u_1,\cdots,u_{b}\}$ from uniform distribution over unit sphere, then update $v_k = \frac{1}{b} \sum_{i\in \mathcal{S}} \big( \hat{\nabla}_{\texttt{uni}} f_i(x_k) - \hat{\nabla}_{\texttt{uni}} f_i(x_{k-1}) \big) + v_{k-1}$;
   \ENDIF{}
   \STATE{}  $ y^{k+1}_j= \arg\min_{y_j\in \mathbb{R}^p} \tilde{\mathcal {L}}_{\rho,j} (x_k,y^{k+1}_{[j-1]},y_j,y^{k}_{[j+1:m]},\lambda_k)$ defined in \eqref{eq:y8}, for all $j\in [m]$;
   \STATE{}  $ x_{k+1}= \arg\min_{x\in \mathbb{R}^d} \hat{\mathcal {L}}_{\rho}\big( x,y^{k+1}_{[m]},\lambda_k,v_k \big)$ defined in \eqref{eq:10};
   \STATE{}  $ \lambda_{k+1} = \lambda_k- \rho(Ax_{k+1} + \sum_{j=1}^mB_jy^{k+1}_j - c)$;
   \ENDFOR
  \STATE {\bfseries Output \ (in theory):} $\{x_\zeta,y_{[m]}^{\zeta},\lambda_\zeta\}$ chosen uniformly randomly from $\{x_{k},y_{[m]}^{k},\lambda_k\}_{k=0}^{K-1}$.
  \STATE {\bfseries Output \ (in practice):} $\{x_{K},y_{[m]}^{K},z_K\}$.
\end{algorithmic}
\end{algorithm}

In Algorithm \ref{alg:1}, we use the following semi-stochastic zeroth-order gradient to update $x$:
\begin{align}
  v_k = \left\{ \begin{aligned}
 & \frac{1}{n} \sum_{i=1}^n \hat{\nabla}_{\texttt{coo}} f_i(x_k), \quad \mbox{if} \mod(k,q)=0 \\
 &  \frac{1}{b} \sum_{i\in \mathcal{S}} \big( \hat{\nabla}_{\texttt{coo}} f_i(x_k) - \hat{\nabla}_{\texttt{coo}} f_i(x_{k-1}) \big) + v_{k-1}, \ \mbox{otw.}
 \end{aligned} \right. \nonumber
\end{align}
where \textbf{otw.} denotes otherwise.
In fact, we have $\mathbb{E}_{\mathcal{S}}[v_k] = \hat{\nabla}_{\texttt{coo}} f(x_k) \neq \nabla f(x_k)$, \emph{i.e.}, this stochastic gradient
is a \textbf{biased} estimate of the true full gradient.
Thus, we choose the appropriate step size
$\eta$, penalty parameter $\rho$ and smoothing parameter $\mu$ to guarantee the convergence of our algorithms,
which will be discussed in the following convergence analysis.

From the above \eqref{eq:6a}, the CooGE needs $2d$ function values to estimate one zeroth-order gradient. Clearly, due to relying on the CooGE, Algorithm \ref{alg:1} has an expensive computation cost. To alleviate this issue, we use the UniGE to replace part of the CooGE use in Algorithm \ref{alg:1}, since the UniGE only requires two function values to estimate one zeroth-order gradient in the above \eqref{eq:7a}.
Thus, we propose a ZO-SPIDER-ADMM (CooGE+UniGE) algorithm based on the mixture of CooGE and UniGE gradient estimators.
Algorithm \ref{alg:2} shows the ZO-SPIDER-ADMM (CooGE+UniGE) algorithm. In Algorithm \ref{alg:2},
we use the following zeroth-order gradient estimator
\begin{align}
  v_k = \left\{ \begin{aligned}
 & \frac{1}{n} \sum_{i=1}^n \hat{\nabla}_{\texttt{coo}} f_i(x_k), \quad \mbox{if} \mod(k,q)=0 \\
 &  \frac{1}{b} \sum_{i\in \mathcal{S}} \big( \hat{\nabla}_{\texttt{uni}} f_i(x_k) - \hat{\nabla}_{\texttt{uni}} f_i(x_{k-1}) \big) + v_{k-1}, \ \mbox{otw.}
  \end{aligned} \right. \nonumber
\end{align}
Then updating the parameters $\{x,y_{[m]},\lambda\}$ is the same as the above algorithm \ref{alg:1}.

\begin{algorithm}[htb]
   \caption{ ZOO-ADMM+(CooGE) Algorithm }
   \label{alg:3}
\begin{algorithmic}[1]
   \STATE {\bfseries Input:} Total iteration $K$, mini-batch sizes $b_1, \ b_2$, epoch size $q$, penalty parameter $\rho$ and step size $\eta$;
    \STATE {\bfseries Initialize:} $x_0 \in \mathbb{R}^d$, $y^0_j \in \mathbb{R}^p, \ j\in [m]$ and $\lambda_0 \in \mathbb{R}^l$;
   \FOR{$k=0,1,\cdots,K-1$}
   \IF{$\mod(k,q)=0$}
   \STATE{} Draw independently $\mathcal{S}_1 \ (|\mathcal{S}_1|=b_1)$ samples $\{\xi_i\}_{i=1}^{b_1}$, and compute
            $v_k = \frac{1}{b_1} \sum_{i\in \mathcal{S}_1}\hat{\nabla}_{\texttt{coo}} f(x_k;\xi_i)$;
   \ELSE{}
   \STATE{} Draw independently $\mathcal{S}_2 \ (|\mathcal{S}_2|=b_2)$ samples $\{\xi_i\}_{i=1}^{b_2}$, and compute
            $v_k = \frac{1}{b_2} \sum_{i\in \mathcal{S}_2} \big( \hat{\nabla}_{\texttt{coo}}f(x_k;\xi_i) - \hat{\nabla}_{\texttt{coo}} f_i(x_{k-1};\xi_i) \big)+ v_{k-1}$;
   \ENDIF{}
   \STATE{}  $ y^{k+1}_j= \arg\min_{y_j\in \mathbb{R}^p} \tilde{\mathcal {L}}_{\rho,j} (x_k,y^{k+1}_{[j-1]},y_j,y^{k}_{[j+1:m]},\lambda_k)$ defined in \eqref{eq:y8}, for all $j\in [m]$;
   \STATE{}  $ x_{k+1}= \arg\min_{x\in \mathbb{R}^d} \hat{\mathcal {L}}_{\rho}\big( x,y^{k+1}_{[m]},\lambda_k,v_k \big)$ defined in \eqref{eq:10};
   \STATE{}  $ \lambda_{k+1} = \lambda_k- \rho(Ax_{k+1} + \sum_{j=1}^mB_jy^{k+1}_j - c)$;
   \ENDFOR
\STATE {\bfseries Output \ (in theory):} $\{x_\zeta,y_{[m]}^{\zeta},\lambda_\zeta\}$ chosen uniformly randomly from $\{x_{k},y_{[m]}^{k},\lambda_k\}_{k=0}^{K-1}$.
 \STATE {\bfseries Output \ (in practice):} $\{x_{K},y_{[m]}^{K},z_K\}$.
\end{algorithmic}
\end{algorithm}

\subsection{ZOO-ADMM+ Algorithms}
In this subsection, we extend the above ZO-SPIDER-ADMM methods to the online setting and
propose a class of faster zeroth-order online ADMM (\emph{i.e.}, ZOO-ADMM+) method to solve the online problem \eqref{eq:1}.
We begin with proposing the ZOO-ADMM+(CooGE) algorithm based on the CooGE gradient estimator,
which is described in Algorithm \ref{alg:3}.

Under the online setting, $f(x)=\mathbb{E}_{\xi} [ f(x;\xi) ]$ denotes a population risk
over an underlying data distribution. The online problem \eqref{eq:1} can be viewed as
having infinite samples, so we are not able to estimate zeroth-order gradient of the function $f(x)$.
Thus, we estimate the mini-batch zeroth-order gradient instead of the full zeroth-order gradient.
In Algorithm \ref{alg:3}, we use the zeroth-order stochastic gradient as follows:
\begin{align}
 v_k \!=\! \left\{ \begin{aligned}
 & \frac{1}{b_1} \sum_{i\in \mathcal{S}_1} \hat{\nabla}_{\texttt{coo}} f(x_k;\xi_i), \quad \mbox{if} \mod(k,q)=0 \\
 & \frac{1}{b_2} \sum_{i\in \mathcal{S}_2} \big( \hat{\nabla}_{\texttt{coo}} f(x_k;\xi_i) - \hat{\nabla}_{\texttt{coo}} f(x_{k-1};\xi_i) \big) + v_{k-1}, \ \mbox{otw.}
  \end{aligned} \right. \nonumber
\end{align}

\begin{algorithm}[htb]
   \caption{ ZOO-ADMM+(CooGE+UniGE) Algorithm }
   \label{alg:4}
\begin{algorithmic}[1]
   \STATE {\bfseries Input:} Total iteration $K$, mini-batch sizes $b_1, \ b_2$, epoch size $q$, penalty parameter $\rho$ and step size $\eta$;
    \STATE {\bfseries Initialize:} $x_0 \in \mathbb{R}^d$, $y^0_j \in \mathbb{R}^p, \ j\in [m]$ and $\lambda_0 \in \mathbb{R}^l$;
   \FOR{$k=0,1,\cdots,K-1$}
   \IF{$\mod(k,q)=0$}
   \STATE{} Draw independently $\mathcal{S}_1 \ (|\mathcal{S}_1|=b_1)$ samples $\{\xi_i\}_{i=1}^{b_1}$, and compute
            $v_k = \frac{1}{b_1} \sum_{i\in \mathcal{S}_1}\hat{\nabla}_{\texttt{coo}} f(x_k;\xi_i)$;
   \ELSE{}
   \STATE{} Draw independently $\mathcal{S}_2 \ (|\mathcal{S}_2|=b_2)$ samples $\{\xi_i\}_{i=1}^{b_2}$, and draw \emph{i.i.d.} $\{u_1,\cdots,u_{b_2}\}$ from uniform distribution
            over unit sphere, then compute
            $v_k = \frac{1}{b_2} \sum_{i\in \mathcal{S}_2} \big( \hat{\nabla}_{\texttt{uni}}f(x_k;\xi_i) - \hat{\nabla}_{\texttt{uni}} f(x_{k-1};\xi_i) \big)+ v_{k-1}$;
   \ENDIF{}
   \STATE{}  $y^{k+1}_j= \arg\min_{y_j\in \mathbb{R}^p} \tilde{\mathcal {L}}_{\rho,j} (x_k,y^{k+1}_{[j-1]},y_j,y^{k}_{[j+1:m]},\lambda_k)$ defined in \eqref{eq:y8}, for all $j\in [m]$;
   \STATE{}  $x_{k+1}= \arg\min_{x\in \mathbb{R}^d} \hat{\mathcal {L}}_{\rho}\big( x,y^{k+1}_{[m]},\lambda_k,v_k \big)$ defined in \eqref{eq:10};
   \STATE{}  $\lambda_{k+1} = \lambda_k- \rho(Ax_{k+1} + \sum_{j=1}^mB_jy^{k+1}_j - c)$;
   \ENDFOR
\STATE {\bfseries Output \ (in theory):} $\{x_\zeta,y_{[m]}^{\zeta},\lambda_\zeta\}$ chosen uniformly randomly from $\{x_{k},y_{[m]}^{k},\lambda_k\}_{k=0}^{K-1}$.
\STATE {\bfseries Output \ (in practice):} $\{x_{K},y_{[m]}^{K},z_K\}$.
\end{algorithmic}
\end{algorithm}

Similarly, we propose a ZOO-ADMM+ (CooGE+UniGE) method based on the mixture of CooGE and UniGE gradient estimators.
Algorithm \ref{alg:4} shows the ZOO-ADMM+ (CooGE+UniGE) method. In Algorithm \ref{alg:4},
we use zeroth-order gradient estimator as follows:
\begin{align}
 v_k \! = \! \left\{ \begin{aligned}
 & \frac{1}{b_1} \sum_{i\in \mathcal{S}_1} \hat{\nabla}_{\texttt{coo}} f(x_k;\xi_i), \quad \mbox{if} \mod(k,q)=0 \\
 & \frac{1}{b_2}  \sum_{i\in \mathcal{S}_2} \big( \hat{\nabla}_{\texttt{uni}} f(x_k;\xi_i) - \hat{\nabla}_{\texttt{uni}} f(x_{k-1};\xi_i) \big) + v_{k-1}, \ \mbox{otw.}
  \end{aligned} \right. \nonumber
\end{align}

\section{Theoretical Analysis}
In the section, we study the convergence properties of the above proposed algorithms.
All related proofs are provided in the supplementary document. We first give some standard assumptions regarding the problem \eqref{eq:1}, and define the standard $\epsilon$-stationary point of the problem \eqref{eq:1}, as used in \cite{jiang2019structured,huang2019fasters}.

Throughout the paper, let $c_k = \lfloor k/q\rfloor$ such that $c_kq \leq k \leq (c_k+1)q-1$. Let the sequence $\{x_k,y_{[m]}^k,\lambda_k\}_{k=1}^K$ be generated from our algorithms, we define a useful variable
\begin{align}
\theta_k & = \mathbb{E}\big[ \|x_{k+1}-x_{k}\|^2+\|x_{k}-x_{k-1}\|^2+\frac{1}{q}\sum_{i=c_kq}^k \|x_{i+1}-x_i\|^2 \nonumber \\
& \quad + \sum_{j=1}^m \|y_j^k-y_j^{k+1}\|^2 \big].
\end{align}

\begin{definition} \label{def:1}
Given $\epsilon>0$, the point $(x^*,y_{[m]}^*,\lambda^*)$ is said to be an $\epsilon$-stationary point of the problem \eqref{eq:1},
if it holds that
\begin{align}
 \mathbb{E}\big[ \mbox{dist}(0,\partial \mathcal{L}(x^*,y^*_{[m]},\lambda^*))^2 \big] \leq \epsilon,
\end{align}
where $f(x)=\frac{1}{n}\sum_{i=1}^nf_i(x)$ or $f(x)=\mathbb{E}[f(x;\xi)]$, and $\mathcal{L}(x,y_{[m]},\lambda) = f(x) + \sum_{j=1}^m \psi_j(y_j) - \langle \lambda, Ax + \sum_{j=1}^mB_jy_j-c \rangle$, and
\begin{align}
   \partial \mathcal{L}(x,y_{[m]},\lambda) = \left [ \begin{matrix}
     \nabla_x \mathcal{L}(x,y_{[m]},\lambda) \\
     \partial_{y_1} \mathcal{L}(x,y_{[m]},\lambda) \\
      \cdots \\
     \partial_{y_m} \mathcal{L}(x,y_{[m]},\lambda) \\
     -Ax-\sum_{j=1}^mB_jy_j+c
 \end{matrix}
 \right ], \nonumber
\end{align}
and $\mbox{dist}(0,\partial \mathcal{L})=\inf_{\mathcal{L}'\in \partial \mathcal{L}} \|0-\mathcal{L}'\|.$
\end{definition}

\begin{assumption} \label{ass:1}
Each loss function $f_i(x)$ and $f(x,\xi_i)$ is $L$-smooth such that, for any $x,y \in \mathbb{R}^d$
\begin{align}
&\|\nabla f_i(x)-\nabla f_i(y)\| \leq L \|x - y\|, \nonumber \\
&\|\nabla f(x;\xi_i)-\nabla f(y;\xi_i)\| \leq L \|x - y\|.  \nonumber
\end{align}
\end{assumption}

\begin{assumption} \label{ass:2}
The functions $\psi_j(y_j)$ for all $j\in [m]$ are convex and possibly nonsmooth.
\end{assumption}

\begin{assumption} \label{ass:3}
$f(x) $ and $\psi_j(y_j)$ for all $j\in [m]$ are all lower bounded, and let $f^*=\inf_x f(x) > - \infty$ and $\psi_j^*=\inf_{y_j} \psi_j(y_j) > - \infty$.
\end{assumption}

\begin{assumption} \label{ass:4}
Matrix $A$ is full row or column rank.
\end{assumption}

\begin{assumption} \label{ass:5}
For the online setting, the variance of unbiased stochastic gradient is bounded, \emph{i.e.},
$\mathbb{E}[\nabla f(x;\xi)]=\nabla f(x)$ and $\mathbb{E} \|\nabla f(x;\xi)-\nabla f(x)\|^2 \leq \sigma^2$ for all $x \in \mathbb{R}^d$.
\end{assumption}

Assumption 1 imposes the smoothness on individual loss functions, which is commonly used in the convergence analysis of variance-reduced algorithms \cite{allen2016variance,reddi2016stochastic,fang2018spider,wang2019spiderboost}.
For any $x_1,x_2 \in \mathbb{R}^d$, we have $\|\nabla f(x_1) - \nabla f(x_2)\| \leq \|\mathbb{E}[\nabla f(x_1;\xi) - \nabla f(x_2;\xi)]\| \leq \mathbb{E}\|\nabla f(x_1;\xi) - \nabla f(x_2;\xi)\| \leq L\|x_1 - x_2\|$. Assumption 1 also implies the full gradient is $L$-Lipschitz continuous.
Assumption 2 is widely used in ADMM methods \cite{huang2019fasters} and proximal methods \cite{ghadimi2016mini,Huang2019faster}. For example, the sparse regularization $\psi_1(y_1)=\|y_1\|_1$.
Assumptions 3 guarantees the feasibility of the optimization problem, which has been used in the study of nonconvex ADMMs
\cite{hong2016convergence,jiang2019structured}.
For example, given a set of training samples $(a_i,b_i)_{i=1}^n$,
where $a_i\in \mathbb{R}^d$, $b_i \in \{-1,1\}$, we consider the nonconvex smooth sigmoid loss function $f(x)=\frac{1}{n}\sum_{i=1}^n \frac{1}{1+\exp(b_i a_i^Tx)}$ used in \cite{allen2016variance,huang2019fasters}. Clearly, we have $ f(x)\geq 0$. Let $\psi_1(y_1)=\|y_1\|_1$, we have $
\psi_1(y_1)\geq 0 $ for any $y_1\in \mathbb{R}^p$.

Assumption 4 guarantees the matrix $A^TA$ or $AA^T$ is non-singular, which is commonly used in the convergence analysis of nonconvex ADMM algorithms \cite{hong2016convergence,jiang2019structured}.
Without loss of generality, we will use the full column rank matrix $A$ as in \cite{huang2019fasters}.
Assumption 5 shows that the variance of stochastic gradient is bounded in norm, which is widely used in the study of
nonconvex methods \cite{ghadimi2016mini} and zeroth-order methods \cite{gao2018information}.
Considering the nonconvex cases where there are multiple global minima, we admit that Assumptions~\ref{ass:1} and~\ref{ass:4} are not mild. Thus, we hope to relax these assumptions in the future work. For example, we do not require full row or column rank of matrix $A$. Next, we give a useful lemma.

\begin{lemma} \label{lem:1}
Suppose that  $\inf_x f(x) = f^* > -\infty$, then
\begin{align}
\inf_x\left\{f(x) - \frac{\beta}{2L}\|\nabla f(x)\|^2\right\} \geq f^* > -\infty,
\end{align}
where $\beta \in [0,1]$.
\end{lemma}

 In the convergence analysis, we admit that the main proofs in our paper can follow the proofs in \cite{huang2019fasters}. However, the main differences between our convergence analysis and the convergence analysis in \cite{huang2019fasters} come from the variances of (zeroth-order) stochastic gradient estimators. For example, in our convergence analysis, the Lyapunov (or potential) functions defined in the proofs are not monotonous, while the Lyapunov (or potential) functions defined in \cite{huang2019fasters} are monotonous. Clearly, our proofs are more complex than the proofs in \cite{huang2019fasters}. For example, we need to choose some appropriate parameters (such as smoothing parameter $\mu$ in ZO-SIPDER-ADMM (Coord) algorithm) control the convergence of our algorithms.

In fact, Lemma~\ref{lem:1} instead of the bounded norm of objective function, \emph{i.e.}, $\|\nabla f(x)\|^2\leq \delta^2$ or $\|\nabla f(x;\xi)\|^2\leq \delta^2$, which is widely used in the nonconvex (stochastic) ADMM methods~\cite{huang2019fasters}. Thus, our convergence analysis uses some milder conditions than the existing convergence analysis in \cite{huang2019fasters}.

\subsection{ Convergence Analysis of ZO-SPIDER-ADMM (CooGE) Algorithm }
In the subsection, we study the convergence properties of the ZO-SPIDER-ADMM (CooGE) algorithm.

\begin{theorem} \label{th:1}
Suppose the sequence $\{x_k,y_{[m]}^k,\lambda_k)_{k=1}^K$ be generated from Algorithm \ref{alg:1}. Let
 $b=q$, $\eta = \frac{\alpha\sigma_{\min}(G)}{4L} \ (0<\alpha \leq 1)$, $\rho = \frac{2\sqrt{237} \kappa_G L}{\sigma^A_{\min}\alpha}$,
 we have
 \begin{align}
 \mathbb{E}\big[ \mbox{dist}(0,\partial L(x_\zeta,y_{[m]}^\zeta,\lambda_\zeta))^2\big]\leq  O(\frac{1}{K}) + O(d\mu^2). \nonumber
\end{align}
It implies that the iteration number $K$ and the smoothing parameter $\mu$ satisfy $K = O(\frac{1}{\epsilon}), \ \mu = O(\sqrt{\frac{\epsilon}{d}})$,
then $(x_{k^*},y^{k^*}_{[m]},\lambda_{k^*})$ is an $\epsilon$-approximate stationary point of the problem \eqref{eq:1}, where $k^* = \mathop{\arg\min}_{k}\theta_{k}$.
\end{theorem}

\begin{remark}
Theorem \ref{th:1} shows that given $b=q$ and $\mu = \frac{1}{\sqrt{dK}}$, the Algorithm \ref{alg:1} has a convergence rate of $O(\frac{1}{K})$.
It implies that given $b=q=\sqrt{n}$, $K = O(\frac{1}{\epsilon}), \ \mu = O(\sqrt{\frac{\epsilon}{d}})$, the ZO-SPIDER-ADMM (CooGE) algorithm
reach a lower function query complexity of $O(nd + n^{\frac{1}{2}}d\epsilon^{-1})$ for finding an $\epsilon$-stationary point
of the \textbf{finite-sum} problem \eqref{eq:1}.
Specifically, when $\mod(k,q)=0$, our algorithm needs $2nd$ function values to compute a zeroth-order full gradient at each iteration, and needs $K/q= n^{-\frac{1}{2}}\epsilon^{-1}$ iterations.
When $\mod(k,q)\neq0$, our algorithm needs $4bd$ function values to compute the zeroth-order stochastic gradient at each iteration, and needs $K=\frac{1}{\epsilon}$ iterations.
Meanwhile, our algorithm needs to calculate a zeroth-order full gradient at least.
Thus, the ZO-SPIDER-ADMM (CooGE) requires the function query complexity $nd+ 2ndK/q + 4bdK=nd+ 2ndn^{-\frac{1}{2}}\epsilon^{-1}+ 4n^{-\frac{1}{2}}d\epsilon^{-1} =O(nd+n^{\frac{1}{2}}d\epsilon^{-1})$
for obtaining an $\epsilon$-approximate stationary point.
\end{remark}

\subsection{ Convergence Analysis of ZO-SPIDER-ADMM (CooGE+UniGE) Algorithm}
In the subsection, we study the convergence properties of the ZO-SPIDER-ADMM (CooGE+UniGE) algorithm.

\begin{theorem} \label{th:2}
Suppose the sequence $\{x_k,y_{[m]}^k,\lambda_k)_{k=1}^K$ be generated from the Algorithm \ref{alg:2}.
Further, let $b=qd$, $\eta = \frac{\alpha\sigma_{\min}(G)}{4L} \ (0<\alpha \leq 1)$ and $\rho = \frac{2\sqrt{237} \kappa_G L}{\sigma^A_{\min}\alpha}$,
we have
\begin{align}
\mathbb{E}\big[ \mbox{dist}(0,\partial L(x_\zeta,y_{[m]}^\zeta,\lambda_\zeta))^2\big]
\leq  O(\frac{1}{K})+ O(d^2\nu^2) + O(d\mu^2), \nonumber
\end{align}
where $\{x_\zeta,y_{[m]}^{\zeta},\lambda_\zeta\}$ is chosen uniformly randomly from $\{x_{k},y_{[m]}^{k},\lambda_k\}_{k=0}^{K-1}$.
It implies that the iteration number $K$, the smoothing parameter $\nu$ and the mini-batch size $b_1$ satisfy
$K = O(\frac{1}{\epsilon}), \ \nu = O(\frac{\sqrt{\epsilon}}{d}), \ \mu = O(\frac{\sqrt{\epsilon}}{\sqrt{d}})$
then $(x_{k^*},y^{k^*}_{[m]},\lambda_{k^*})$ is an $\epsilon$-stationary point of the problem \eqref{eq:1}, where $k^* = \mathop{\arg\min}_{k}\theta_{k}$.
\end{theorem}

\begin{remark}
Theorem \ref{th:2} shows that given $b=dq$, $\mu = \frac{1}{\sqrt{dK}}$ and $\nu=\frac{1}{d\sqrt{K}}$,
the Algorithm \ref{alg:2} has a convergence rate of $O(\frac{1}{K})$.
It implies that given $q=\sqrt{n}$, $K = O(\frac{1}{\epsilon}), \ \mu = O(\sqrt{\frac{\epsilon}{d}}), \ \nu=O(\frac{\sqrt{\epsilon}}{d})$,
the ZO-SPIDER-ADMM (CooGE+UniGE) algorithm
reach a lower function query complexity of $O(nd + n^{\frac{1}{2}}d\epsilon^{-1})$ for finding an $\epsilon$-stationary point
of the \textbf{finite-sum} problem \eqref{eq:1}.
Specifically, when $\mod(k,q)=0$, our algorithm needs $2nd$ function values to compute a zeroth-order full gradient at each iteration, and needs $K/q= n^{-\frac{1}{2}}\epsilon^{-1}$ iterations.
When $\mod(k,q)\neq0$, our algorithm needs $4b$ function values to compute the zeroth-order stochastic gradient at each iteration, and needs $K=\frac{1}{\epsilon}$ iterations.
Meanwhile, our algorithm needs to calculate a zeroth-order full gradient at least.
Thus, the ZO-SPIDER-ADMM (CooGE+UniGE) requires the function query complexity $nd+ 2ndK/q + 4bK=nd+ 2ndn^{-\frac{1}{2}}\epsilon^{-1}+ 4n^{-\frac{1}{2}}d\epsilon^{-1} =O(nd+n^{\frac{1}{2}}d\epsilon^{-1})$
for obtaining an $\epsilon$-stationary point.
\end{remark}

\subsection{ Convergence Analysis of ZOO-ADMM+(CooGE) Algorithm}
In this subsection, we study the convergence properties of the ZOO-ADMM+(CooGE) algorithm.

\begin{theorem}\label{th:3}
Suppose the sequence $\{x_k,y_{[m]}^k,\lambda_k)_{k=1}^K$ be generated from the Algorithm \ref{alg:3}.
Further, let $b_2=q$, $\eta = \frac{\alpha\sigma_{\min}(G)}{4L} \ (0<\alpha \leq 1)$ and $\rho = \frac{2\sqrt{237} \kappa_G L}{\sigma^A_{\min}\alpha}$,
we have
\begin{align}
\mathbb{E}\big[ \mbox{dist}(0,\partial L(x_\zeta,y_{[m]}^\zeta,\lambda_\zeta))^2\big] \leq  O(\frac{1}{K}) + O(d\mu^2)+ O(\frac{1}{b_1}). \nonumber
\end{align}
It implies that the parameters $K$, $b_1$ and $\mu$ satisfy
$K = O(\frac{1}{\epsilon})$, $b_1 = O(\frac{1}{\epsilon})$ and $\mu=\sqrt{\frac{\epsilon}{d}}$,
then $(x_{k^*},y^{k^*}_{[m]},\lambda_{k^*})$ is an $\epsilon$-approximate stationary point of the problem \eqref{eq:1},
where $k^* = \mathop{\arg\min}_{k}\theta_{k}$.
\end{theorem}
\begin{remark}
 Theorem \ref{th:3} shows that given $b_2=q$, $b_1= K$, $\mu=\frac{1}{\sqrt{dK}}$, the ZOO-ADMM+(CooGE) algorithm has a convergence rate of $O(\frac{1}{K})$.
 Let $K=\epsilon^{-1}$, $b_2=q=\epsilon^{-1/2}$, $b_1= \epsilon^{-1}$ and $\mu=\sqrt{\frac{\epsilon}{d}}$,
 the ZOO-ADMM+(CooGE) algorithm reaches a function query complexity of $2db_1K/q+4db_2K=O(d\epsilon^{-\frac{3}{2}})$
 for finding an $\epsilon$-stationary point
 of the \textbf{online} problem \eqref{eq:1}.
\end{remark}

\subsection{ Convergence Analysis of ZOO-ADMM+ (CooGE+UniGE) Algorithm }
In this subsection, we study the convergence properties of the ZOO-ADMM+(CooGE+UniGE) algorithm.

\begin{theorem} \label{th:4}
Suppose the sequence $\{x_k,y_{[m]}^k,\lambda_k)_{k=1}^K$ be generated from Algorithm \ref{alg:4}.
Further, let $b_2=qd$, $\eta = \frac{\alpha\sigma_{\min}(G)}{4L} \ (0<\alpha \leq 1)$
and $\rho = \frac{2\sqrt{237} \kappa_G L}{\sigma^A_{\min}\alpha}$,
we have
\begin{align}
\mathbb{E}\big[ \mbox{dist}(0,\partial L(x_\zeta,y_{[m]}^\zeta,\lambda_\zeta))^2\big] & \leq  O(\frac{1}{K})+ O(d^2\nu^2)  \nonumber \\
& \quad + O(d\mu^2) + O(\frac{1}{b_1}),
\end{align}
where $\{x_\zeta,y_{[m]}^{\zeta},\lambda_\zeta\}$ is chosen uniformly randomly from $\{x_{k},y_{[m]}^{k},\lambda_k\}_{k=0}^{K-1}$.
It implies that given
$K = O(\frac{1}{\epsilon}), \ \nu = O(\frac{\sqrt{\epsilon}}{d}), \ \mu=O(\sqrt{\frac{\epsilon}{d}}), \ b_1 = O(\frac{1}{\epsilon})$
then $(x_{k^*},y^{k^*}_{[m]},\lambda_{k^*})$ is an $\epsilon$-stationary point of the problem \eqref{eq:1},
where $k^* = \mathop{\arg\min}_{k}\theta_{k}$.
\end{theorem}
\begin{remark}
 Theorem \ref{th:4} shows that given $b_2=qd$, $b_1= K$ and $\nu=\frac{1}{d\sqrt{K}}$, the Algorithm \ref{alg:4} has $O(\frac{1}{K})$ convergence rate.
 Let $K=\epsilon^{-1}$, $q=\epsilon^{-1/2}$, $b_2=dq=d\epsilon^{-1/2}$, $b_1= \epsilon^{-1}$, $\nu=\frac{\sqrt{\epsilon}}{d}$ and $\mu=O(\sqrt{\frac{\epsilon}{d}})$, the ZOO-ADMM+(CooGE+UniGE) algorithm reaches the function query complexity of $2b_1K/q+4b_2K=O(d\epsilon^{-\frac{3}{2}})$ for finding an $\epsilon$-stationary point of the \textbf{online} problem \eqref{eq:1}.
\end{remark}

\begin{table}[!h]
  \centering
  \caption{ Four Benchmark Datasets for Attacking Black-Box DNNs }
  \label{tab:2}
  \resizebox{0.49\textwidth}{!}{
  \begin{tabular}{c|c|c|c}
  \hline
  datasets & \#test samples & \#dimension & \#classes \\ \hline
  \emph{MNIST} & 10,000 &  $28\!\times \!28$ & 10 \\
  \emph{Fashion-MNIST} & 10,000 &  $28\!\times\!28$ & 10 \\
  \emph{SVHN} & 26,032 & $32\!\times\!32\times\!3$ & 10 \\
  \emph{CIFAR-10}  & 10,000& $32\!\times32\!\times\!3$ & 10 \\
  \emph{STL-10} (resized)  & 8,000& $32\!\times32\!\times\!3$ & 10 \\
  \emph{CINIC-10}  & 90,000& $32\!\times32\!\times\!3$ & 10
  \\
  \hline
  \end{tabular}
  }
\end{table}

\begin{figure*}[htbp]
\centering
\subfigure[\emph{MNIST}]{\includegraphics[width=0.3\textwidth]{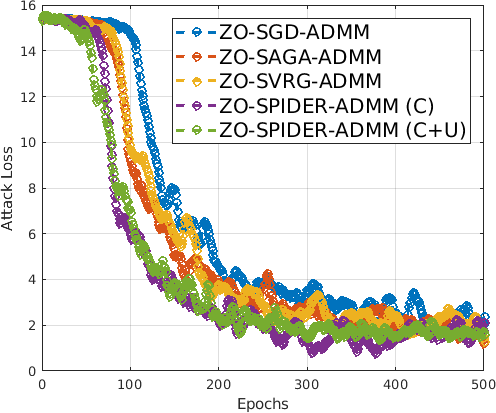}}
\hspace{0.5em}
\subfigure[\emph{Fashion-MNIST}]{\includegraphics[width=0.3\textwidth]{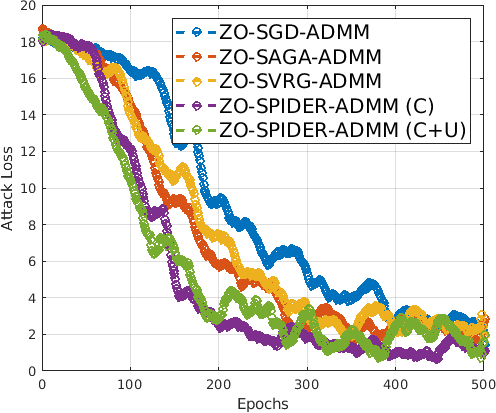}}
\hspace{0.5em}
\subfigure[\emph{SVHN}]{\includegraphics[width=0.3\textwidth]{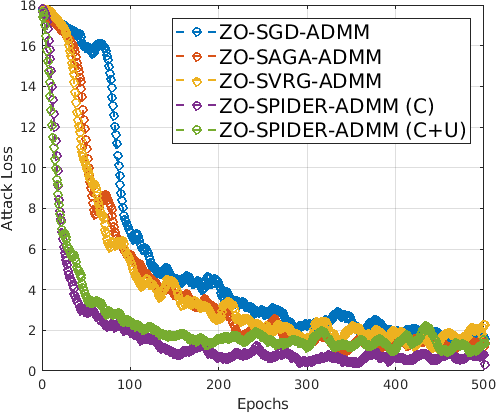}}\\

\subfigure[\emph{CIFAR-10}]{\includegraphics[width=0.3\textwidth]{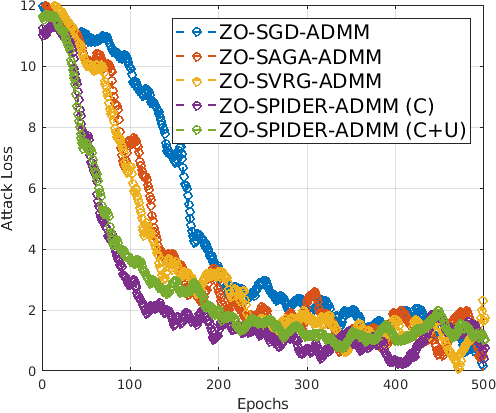}}
\hspace{0.5em}
\subfigure[\emph{STL-10}]{\includegraphics[width=0.3\textwidth]{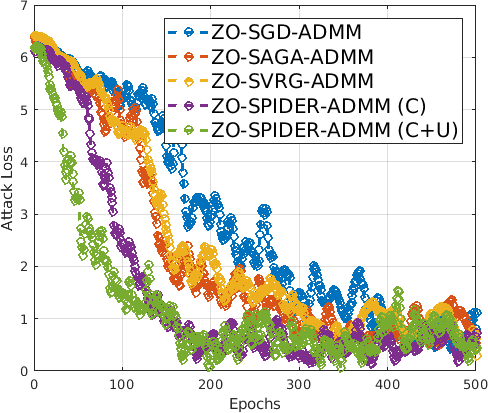}}
\hspace{0.5em}
\subfigure[\emph{CINIC-10}]{\includegraphics[width=0.3\textwidth]{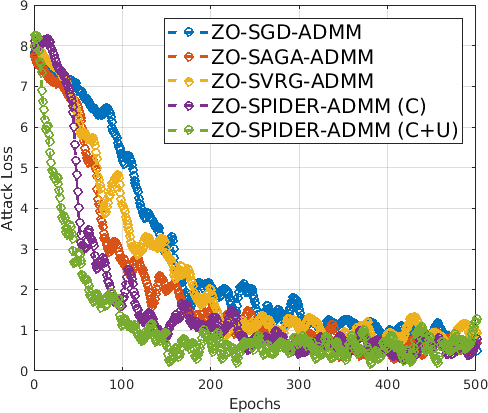}}
\caption{ Attack loss vs. epochs on adversarial attacks black-box DNNs in the finite-sum setting. ZO-SPIDER-ADMM (C) and ZO-SPIDER-ADMM (C+U) represent ZO-SPIDER-ADMM (CooGE) and ZO-SPIDER-ADMM (CooGE+UniGE), respectively.}
\label{fig:2}
\end{figure*}

\begin{figure*}[htbp]
\centering
\subfigure[\emph{CIFAR-10}]{\includegraphics[width=0.3\textwidth]{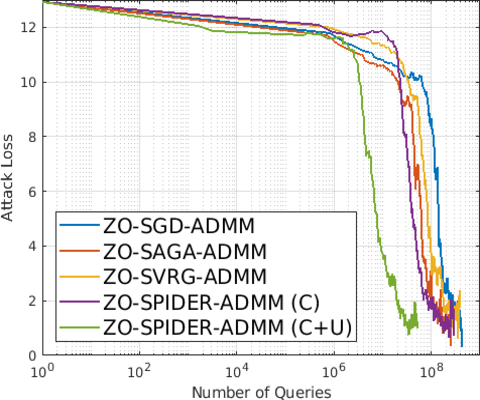}}
\hspace{0.5em}
\subfigure[\emph{STL-10}]{\includegraphics[width=0.3\textwidth]{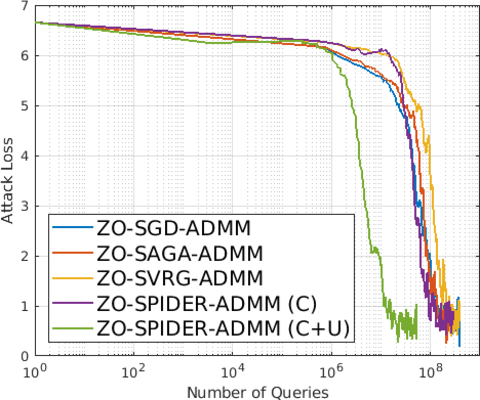}}
\hspace{0.5em}
\subfigure[\emph{CINIC-10}]{\includegraphics[width=0.3\textwidth]{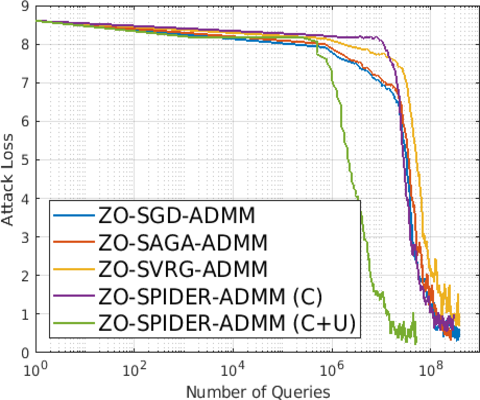}}
\caption{ Attack loss vs. number of queries on adversarial attacks black-box DNNs in the finite-sum setting. }
\label{fig:fs}
\end{figure*}

\begin{figure*}[htbp]
\centering
\subfigure[\emph{MNIST}]{\includegraphics[width=0.3\textwidth]{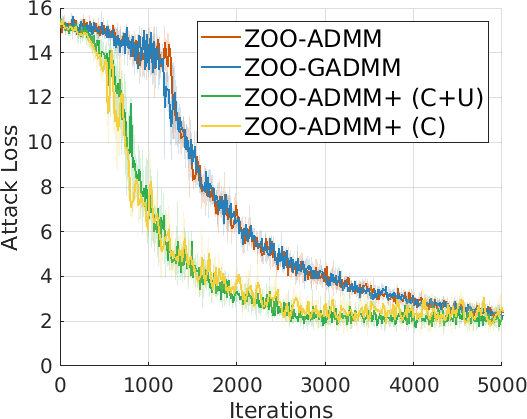}}
\hspace{0.5em}
\subfigure[\emph{Fashion-MNIST}]{\includegraphics[width=0.3\textwidth]{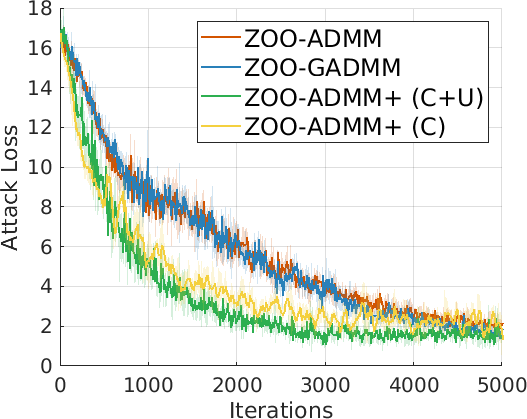}}
\hspace{0.5em}
\subfigure[\emph{SVHN}]{\includegraphics[width=0.3\textwidth]{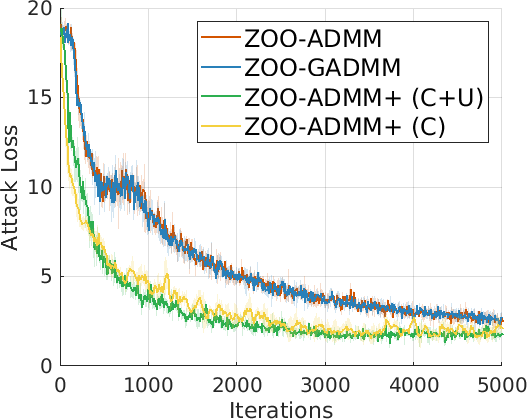}}\\
\hspace{0.5em}
\subfigure[\emph{CIFAR-10}]{\includegraphics[width=0.3\textwidth]{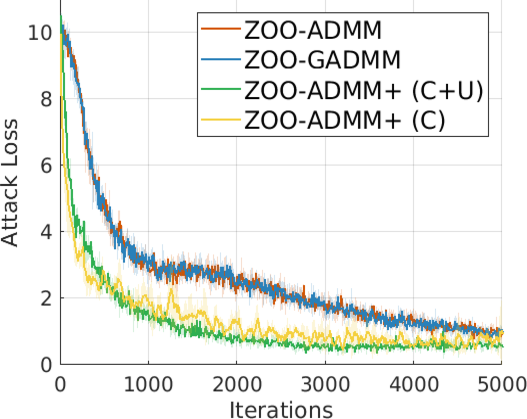}}
\hspace{0.5em}
\subfigure[\emph{STL-10}]{\includegraphics[width=0.3\textwidth]{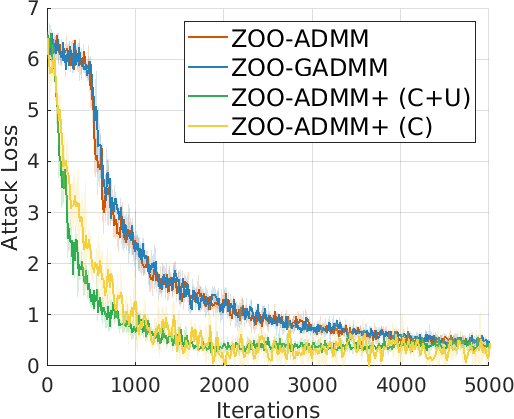}}
\hspace{0.5em}
\subfigure[\emph{CINIC-10}]{\includegraphics[width=0.3\textwidth]{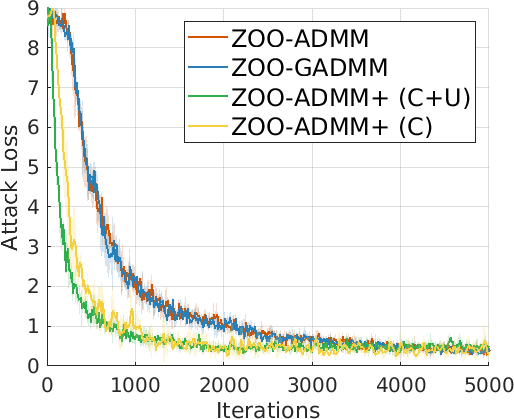}}

\caption{ Attack loss vs. epochs on adversarial attacks black-box DNNs in the online setting. ZOO-ADMM+ (C) and ZOO-ADMM+ (C+U) represent ZOO-ADMM+ (CooGE) and ZOO-ADMM+ (CooGE+UniGE), respectively.}
\label{fig:3}
\end{figure*}

\begin{figure*}[htbp]
\centering
\subfigure[\emph{CIFAR-10}]{\includegraphics[width=0.3\textwidth]{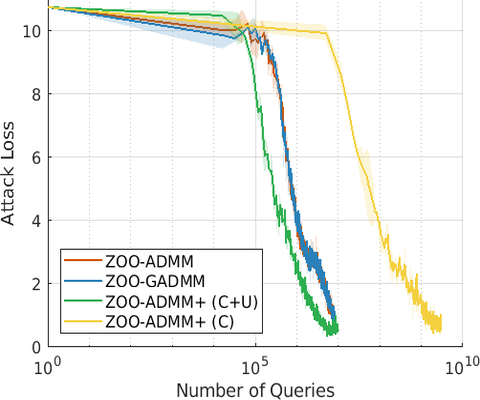}}
\hspace{0.5em}
\subfigure[\emph{STL-10}]{\includegraphics[width=0.3\textwidth]{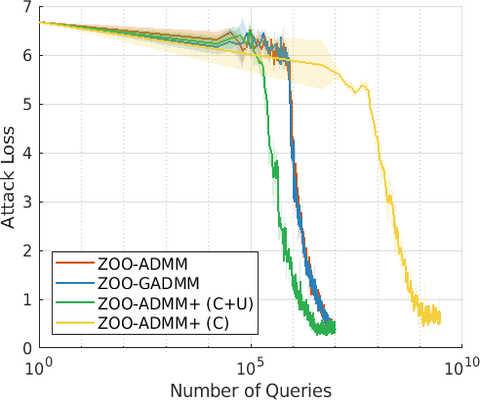}}
\hspace{0.5em}
\subfigure[\emph{CINIC-10}]{\includegraphics[width=0.3\textwidth]{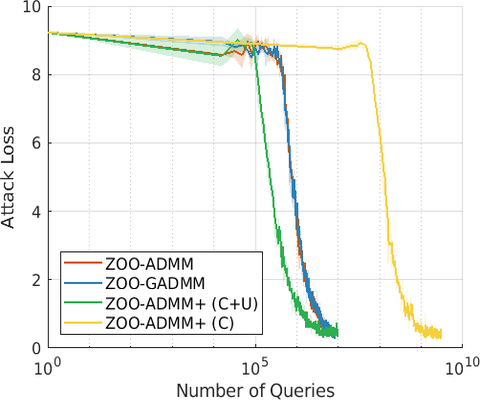}}

\caption{ Attack loss vs. number of queries on adversarial attacks black-box DNNs in the online setting.}
\label{fig:ol}
\end{figure*}

\begin{figure*}[htbp]
\centering
\subfigure[\emph{MNIST}]{\includegraphics[width=0.31\textwidth]{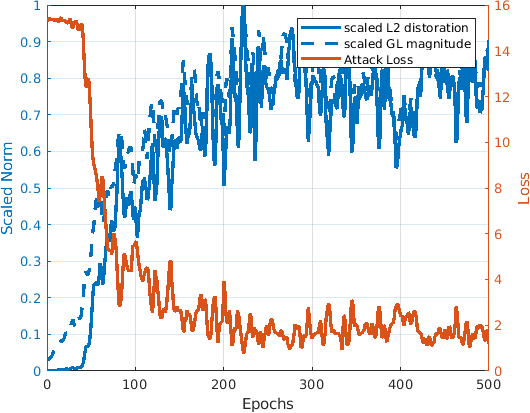}}
\hspace{0.5em}
\subfigure[\emph{Fashion-MNIST}]{\includegraphics[width=0.31\textwidth]{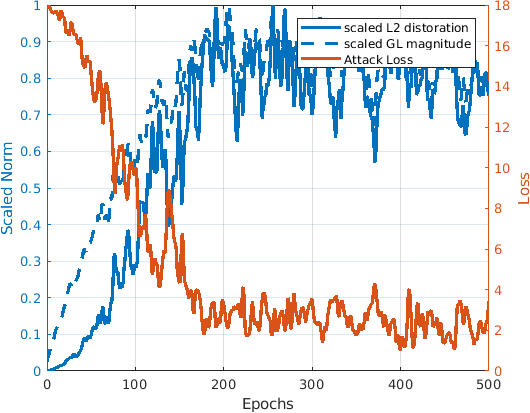}}
\hspace{0.5em}
\subfigure[\emph{SVHN}]{\includegraphics[width=0.31\textwidth]{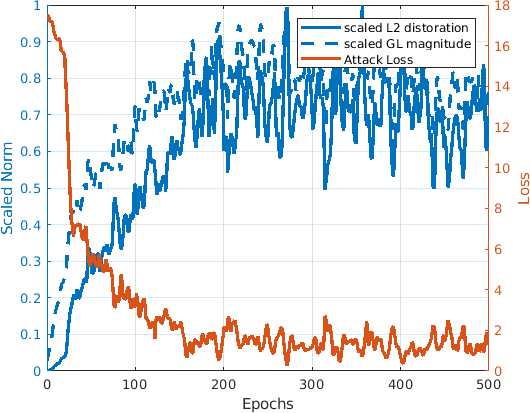}}\\

\subfigure[\emph{CIFAR-10}]{\includegraphics[width=0.31\textwidth]{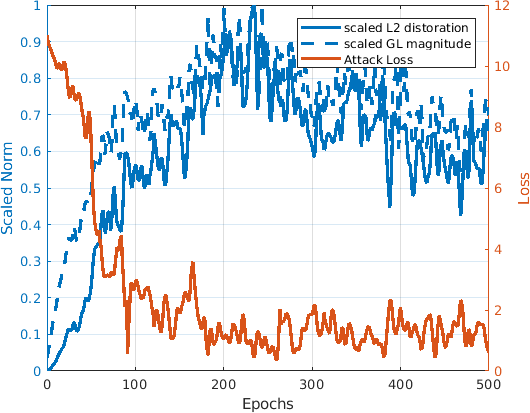}}
\hspace{0.5em}
\subfigure[\emph{STL-10}]{\includegraphics[width=0.31\textwidth]{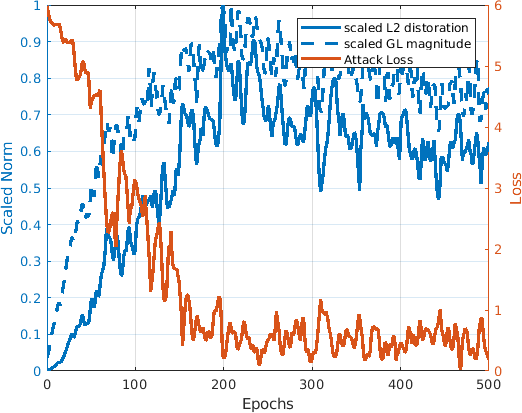}}
\hspace{0.5em}
\subfigure[\emph{CINIC-10}]{\includegraphics[width=0.31\textwidth]{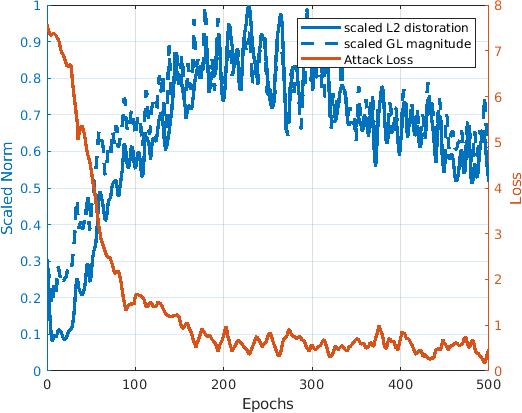}}\\
\caption{  Relationship between attack loss and norms of perturbation ($L_2$ norm distoration and Group lasso magnitude) on adversarial attacks black-box DNNs. }
\label{fig:4}
\end{figure*}

\section{Experiments}
In this section, we apply the universal structured adversarial attack on black-box DNNs to demonstrate efficiency of our algorithms.
In the finite-sum setting, we compare the proposed ZO-SPIDER-ADMM methods with the existing
zeroth-order stochastic ADMM methods (i.e., ZO-SAGA-ADMM and ZO-SVRG-ADMM \cite{Huang2019zeroth}),
and the zeroth-order stochastic ADMM without variance reduction (ZO-SGD-ADMM).
In the online setting, we compare the proposed ZOO-ADMM+ methods with
the ZOO-ADMM \cite{liu2018admm} and ZO-GADM \cite{gao2018information} methods.

\begin{figure}[htbp]
\centering
\subfigure[\emph{ResNet-20}]{\includegraphics[width=0.235\textwidth]{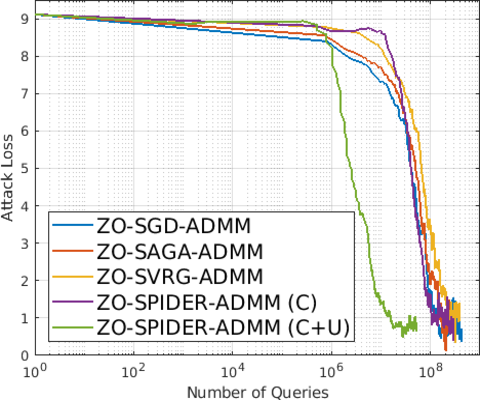}}
\hspace{0.2em}
\subfigure[\emph{7-layer DNN}]{\includegraphics[width=0.235\textwidth]{figs/cifar10_oldq.png}}
\caption{  The attack loss for different models on the CIFAR-10 dataset. }
\label{fig:diff-dnn}
\end{figure}

\subsection{ Experimental Setups }
In this subsection, we apply our algorithms to generate adversarial examples to attack the pre-trained black-box DNNs,
whose parameters are hidden from us and only its outputs are accessible.
In the experiment, we mainly focus on finding a universal structured perturbation to fool the DNNs from multiple images as in \cite{moosavi2017universal,Huang2019zeroth}.
Given the samples $\big\{ a_i\in \mathbb{R}^d , \ l_i \in \{1,2,\cdots,c\} \ i=1,2,\dots,n \big\}$, as in \cite{Huang2019zeroth},
we solve the following problem
\begin{align} \label{eq:15}
\min_{x \in \mathbb{R}^d} & \ \bigg\{\frac{1}{n}\sum\limits_{i=1}^n \max\big\{F_{l_i}(a_i+x) - \max\limits_{j\neq l_i} F_j(a_{i}+x), 0 \big\} \\
 &  + \tau_1 \sum_{p=1}^P\sum_{q=1}^Q \|x_{\mathcal{G}_{p,q}}\| + \tau_2 \|x \|^2 + \tau_3 h(x) \bigg\},\nonumber
\end{align}
where $x\in \mathbb{R}^d$ denotes a universal structured perturbation, and $\{a_i,l_i\}_{i=1}^n$ are the input images and the corresponding labels. Here $F(a)$ represents the final layer output before softmax of neural network,
and $\tau_i \ (i=1,2,3)$ are non-negative tuning parameters.
In the above problem \eqref{eq:15}, the first penalty $\sum_{p,q}\|x_{G_{p,q}}\|$ is to penalize structured permutations as in \cite{xu2018structured}, where the overlapping groups $\{\{\mathcal{G}_{p,q}\}_{p=1}^P\}_{q=1}^Q$
generate from dividing an image into sub-groups of pixels. The second penalty $\|x\|^2$ is to ensure the permutation is small, so that it become harder to detect. The third penalty $h(x)$ is to ensure the attacked image $a_i+x\in [0,1]^d$ is a valid image, defined as
\begin{align}
 h(x) = \left\{ \begin{aligned}
 & 0, \quad \mbox{if} \ \{a_i+x \in [0,1]^d\}_{i=1}^n \ \mbox{and} \ \|x\|_\infty \leq \varepsilon \\
 &  \infty,
 \quad \mbox{otherwise}
  \end{aligned} \right. \nonumber
\end{align}

Let $f_i(x)=\max\big\{F_{l_i}(a_i+x) - \max\limits_{j\neq l_i} F_j(a_{i}+x), 0 \big\}$, then we rewrite the problem \eqref{eq:15} in the following form:
\begin{align} \label{eq:16}
 \min & \ \frac{1}{n}\sum\limits_{i=1}^nf_i(x) + \tau_1 \sum_{j=1}^{PQ}\|y_{j,\mathcal{G}_j}\| + \tau_2\|z\|^2 + \tau_3h(w) ,  \\
 \mbox{s.t.} & \ z = x, \ w = x,\ y_j = x, \ \mbox{for} \ j=1,\cdots,PQ \nonumber
\end{align}
where $y_{j,\mathcal{G}_j}$ denotes the sub-vector of $y_j\in \mathbb{R}^d$ with indices given by group $\mathcal{G}_j$.
Clearly, the above problem \eqref{eq:16} can easy be rewritten as the form of the problem \eqref{eq:1}.

\begin{figure*}[!t]
  \centering
  \includegraphics[width=0.9\textwidth]{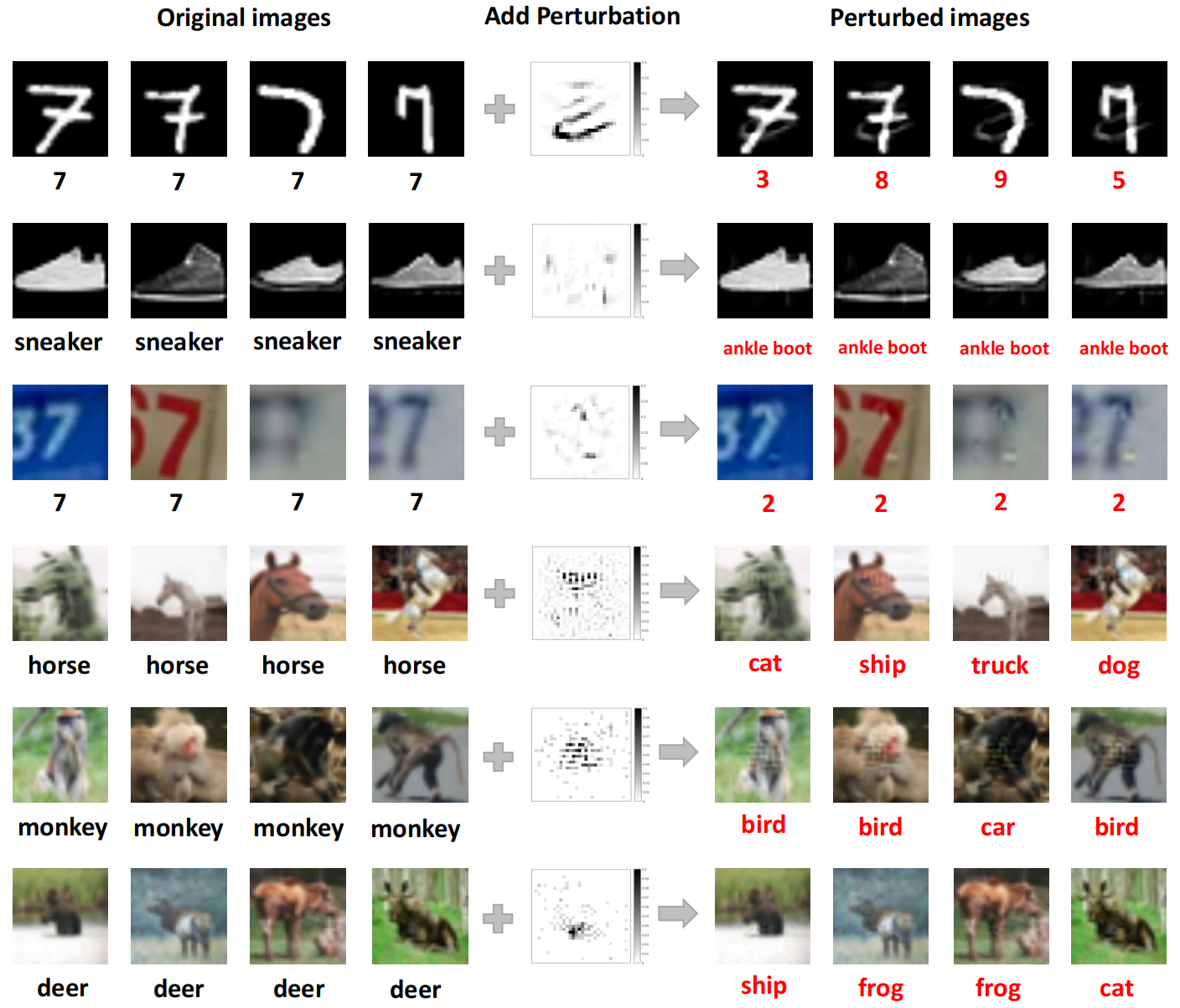}\\
  \caption{ Group-sparsity perturbations are learned from \emph{MNIST}, \emph{Fashion-MNIST}, \emph{SVHN}, \emph{CIFAR-10}, \emph{STL-10} and \emph{CINIC-10} (from top to bottom) datasets.
  Black and red labels denote the initial label and the label after attack, respectively. }
  \label{fig:5}
\end{figure*}

In the experiment, we use the pre-trained DNN models on four benchmark datasets in Table \ref{tab:2} as the target black-box models.
Specifically, the pre-trained DNNs on MNIST, Fashion-MNIST, SVHN, CIFAR-10, STL-10~\cite{coates2011analysis} and CINIC-10~\cite{darlow2018cinic} can attain $99.4\%$, $91.8\%$, $93.2\%$, $80.8\%$, $71.8\%$ and $83.3\%$ test accuracy, respectively.
For MNIST and FashionMNIST, a 5-layer DNN (3 conv layers and 2 dense layers) is used. For the CINIC-10 dataset, ResNet-32 is used. For the rest dataset, an 7-layer (5 conv layers and 2 dense layers) DNN is used.

In the experiment, we set the parameters $\varepsilon=0.4$, $\tau_1 = 1$, $\tau_2=2$ and $\tau_3=1$.
In the zeroth-order gradient estimator, we choose the smoothing parameters $\mu=\frac{1}{\sqrt{d k}}$ and $\nu=\frac{1}{d \sqrt{k}}$.
For all datasets, the kernel size for overlapping group lasso is set to $3 \times 3$ and the stride is one.
In the \textbf{finite-sum} setting, we select 40 samples from each class, and choose 4 as the batch size.
In the \textbf{online} setting, we select 400 samples from each class, and choose 10 as the batch size.

\subsection{ Experimental Results }
Figure \ref{fig:2} illustrates that, in the finite-sum setting, the attack losses (\emph{i.e.} the first term of the problem \eqref{eq:15}) of
our ZO-SPIDER-ADMM algorithms faster decrease than the other algorithms, as the epoch (\emph{i.e.}, 10 iterations) increases.
Figure \ref{fig:3} shows that, in the online setting, the attack losses of our ZOO-ADMM+ algorithms
also faster decrease than the ZOO-ADMM and ZO-GADM algorithms, as the number of iteration increases.
These results demonstrate that our algorithms have faster convergence than the existing zeroth-order algorithms in solving the above complex problem \eqref{eq:15}.

Figs.~\ref{fig:fs} and~\ref{fig:ol} show the results on attack losses vs number of queries in the finite-sum and stochastic cases, respectively. In the finite-sum setting, our ZO-SPIDER-ADMM (C+U) has the lowest query costs. In the online setting, our ZOO-ADMM+ (C+U) has the lowest query costs.  ZO-SPIDER-ADMM (C) has similar query costs compared to other methods in the finite sum setting, and ZOO-ADMM+ (C) has the largest query costs in the online learning setting. It is because that all comparison methods use cheaper zeroth order gradient estimators, while our ZO-SPIDER-ADMM (C) and ZOO-ADMM+ (C) methods use the expensive CooGE gradient estimator.
Fig.~\ref{fig:diff-dnn} shows the attack losses vs number of queries in different models (i.e., ResNet-20  and 7-layer DNN) on the CIFAR-10 dataset.
From Fig.~\ref{fig:diff-dnn}, we can find that different models only change the numeric results, but the overall tendency is similar.

Fig. \ref{fig:4} shows that relationship between the attack loss and the norms of perturbation ($L_2$ norm distoration and Group lasso magnitude),
which conforms to the phenomenon of adversarial attack.
From Fig. \ref{fig:5}, we can find that our algorithms can learn some interpretable structured perturbations,
which can successfully attack the corresponding DNNs.

\section{ Conclusions }
In this paper, we proposed a class of faster zeroth-order stochastic ADMM (ZO-SPIDER-ADMM) methods to solve the problem \eqref{eq:1}.
We proved that the ZO-SPIDER-ADMM achieves a lower function query complexity of $O(dn + dn^{\frac{1}{2}}\epsilon^{-1})$ for
finding an $\epsilon$-stationary point, which improves the existing zeroth-order ADMM methods by a factor $O(d^{\frac{1}{3}}n^{\frac{1}{6}})$.
At the same time, we extend the ZO-SPIDER-ADMM method to the online setting, and propose a class of faster zeroth-order online ADMM (ZOO-ADMM+)
methods. Moreover, we proved that the ZOO-ADMM+ reaches a lower function query complexity of $O(d\epsilon^{-\frac{3}{2}})$,
which improves the existing best result by a factor of $O(\epsilon^{-\frac{1}{2}})$.


%

%

\ifCLASSOPTIONcompsoc
  \section*{Acknowledgments}
\else
  \section*{Acknowledgment}
\fi
We would like to express our gratitude to both the editor and the anonymous reviewers for their
valuable comments that significantly improved the quality of our paper. F. Huang was a postdoctoral researcher at the University of Pittsburgh. S.Q. Gao and H. Huang were partially supported by U.S. NSF IIS 1838627, 1837956, 1956002, 2211492, CNS 2213701, CCF 2217003, DBI 2225775.

\ifCLASSOPTIONcaptionsoff
  \newpage
\fi



%
%
%

\bibliographystyle{IEEEtran}
\bibliography{IEEEabrv,ZO-ADMM}

\vfill


\newpage

\begin{onecolumn}

\setcounter{page}{1}

\begin{appendices}

\section{ Supplementary Materials }
In this section, we provide the detail proofs of the above theoretical results.
First, we restate and give some useful lemmas.
Throughout the paper, let $c_k = \lfloor k/q\rfloor$ such that $c_kq \leq k \leq (c_k+1)q-1$.

\begin{lemma}  \label{lem:A0}
(Lemma 7 of \cite{Reddi2016Prox}) For random variables $\xi_1, \cdots, \xi_n$ are independent and mean $0$, we have
\begin{align}
 \mathbb{E}[ \|\xi_1 + \cdots +\xi_n\|^2 ] = \mathbb{E}[ \|\xi_1\|^2 + \cdots + \|\xi_n\|^2].
\end{align}
\end{lemma}

\begin{lemma} \label{lem:A1}
(Lemma 3 of \cite{Ji2019improved}) Let $\hat{\nabla}_{\texttt{coo}}f(x) = \sum_{j=1}^d \frac{f(x + \mu e_j)-f(x - \mu e_j)}{2\mu}e_j$, for
any $x\in \mathbb{R}^d$, we have
\begin{align}
 \mathbb{E}\|\hat{\nabla}_{\texttt{coo}}f(x) - \nabla f(x)\|^2 \leq L^2d\mu^2.
\end{align}
\end{lemma}

\begin{lemma} \label{lem:A2}
(Lemma 6 of \cite{Ji2019improved}) Under Assumption 1, and given $v_k = \frac{1}{|\mathcal{S}_2|} \sum_{i\in \mathcal{S}_2}\big ( \hat{\nabla}_{\texttt{coo}}f_i(x_k) - \hat{\nabla}_{\texttt{coo}}f_i(x_{k-1}) \big) + v_{k-1}$ with $\mathcal{S}_2\subseteq \{1,2,\cdots,n\}$,
then we have, for any $c_kq \leq k \leq \min\{(c_k+1)q-1,K\}$,
\begin{align}
 \mathbb{E}\| v_k - \hat{\nabla}_{\texttt{coo}}f(x_k) \|^2 \leq \frac{3L^2}{|S_2|} \sum_{i=c_kq}^{k-1} \mathbb{E}\|x_{i+1}-x_i\|^2 + (k-c_k q)\frac{6L^2d\mu^2}{|S_2|}
 + \frac{3\mathcal{I}(|S_1|<n)}{|S_1|} (2L^2d\mu^2 + \sigma^2),
\end{align}
where $\mathcal{I}(\cdot)$ is the indicator function, and $\sum_{i=c_kq}^{c_kq-1} \mathbb{E}\|x_{i+1}-x_i\|^2=0$.
\end{lemma}

\begin{lemma} \label{lem:A3}
(Lemma 5 of \cite{Ji2019improved}) Let $f_{\nu}(x)=\mathbb{E}_{u\sim U_B}[f(x+\nu u)]$ be a smooth approximation of $f(x)$,
where $U_B$ is the uniform distribution over the $d$-dimensional unit Euclidean ball $B$. Then we have
\begin{itemize}
\item[(1)] $|f_{\nu}(x) - f(x)| \leq \frac{\nu^2 L}{2}$ and $\|\nabla f_{\nu}(x) - \nabla f(x)\| \leq \frac{\nu Ld}{2}$ for any $x\in \mathbb{R}^d$;
\item[(2)] $\mathbb{E}[\frac{1}{|\mathcal{S}|}\sum_{i\in \mathcal{S}}\hat{\nabla}_{\texttt{uni}}f_i(x)] = \nabla f_{\nu}(x)$ for any $x\in \mathbb{R}^d$;
\item[(3)] $\mathbb{E}\|\hat{\nabla}_{\texttt{uni}}f_i(x_1)-\hat{\nabla}_{\texttt{uni}}f_i(x_2)\|^2 \leq 3dL^2\|x_1-x_2\|^2 + \frac{3L^2d^2\nu^2}{2}$ for any $i$ and any $x_1, x_2 \in \mathbb{R}^d$ .
\end{itemize}
\end{lemma}

\begin{lemma} \label{lem:A4}
(Restatement of Lemma \ref{lem:1})
Suppose that  $\inf_x f(x) = f^* > -\infty$, then
\begin{align}
\inf_x\left\{f(x) - \frac{\beta}{2L}\|\nabla f(x)\|^2\right\} \geq f^* > -\infty,
\end{align}
where $\beta \in [0,1]$.
\end{lemma}

\begin{proof}
\begin{eqnarray*}
f^* &\leq& \inf_x f\big( x - \frac{\beta}{L}\nabla f(x) \big) \\
&\leq & \inf_x \left\{f(x) + \langle \nabla f(x), x -\frac{\beta}{L}\nabla f(x) - x\rangle + \frac{L}{2}\|x - \frac{\beta}{L}\nabla f(x) -x \|^2\right\} \\
& \leq & \inf_x \left\{f(x) - \frac{\beta}{2L}\|\nabla f(x)\|^2 \right\},
\end{eqnarray*}
where the second inequality is due to Assumption 1, and the last inequality holds by $\beta \in [0,1]$.

\end{proof}

\subsection{ Convergence Analysis of ZO-SPIDER-ADMM (CooGE) Algorithm }
\label{Appendix:A1}
In the subsection, we study the convergence properties of the ZO-SPIDER-ADMM (CooGE) algorithm.
We begin with giving some useful lemmas.

\begin{lemma} \label{lem:B1}
Suppose the zeroth-order gradient $v_k$ be generated from Algorithm \ref{alg:1}, we have
\begin{align}
 \mathbb{E}\| v_k - \nabla f(x_k) \|^2 \leq \frac{6L^2}{b} \sum_{i=c_kq}^{k-1} \mathbb{E}\|x_{i+1}-x_i\|^2 + \frac{12qL^2d\mu^2}{b} + 2 L^2d\mu^2.
\end{align}
\end{lemma}
\begin{proof}
In Algorithm \ref{alg:1}, we have $|S_2|=b$, and $|S_1|=n$. By Lemma \ref{lem:A2}, we have
\begin{align} \label{eq:B1}
 \mathbb{E}\| v_k - \hat{\nabla}_{\texttt{coo}}f(x_k) \|^2 & \leq \frac{3L^2}{b} \sum_{i=c_kq}^{k-1} \mathbb{E}\|x_{i+1}-x_i\|^2 + (k-c_k q)\frac{6L^2d\mu^2}{b}
 + \frac{3\mathcal{I}(|S_1|<n)}{|S_1|} (2L^2d\mu^2 + \sigma^2) \nonumber \\
 & \leq \frac{3L^2}{b} \sum_{i=c_kq}^{k-1} \mathbb{E}\|x_{i+1}-x_i\|^2 + \frac{6qL^2d\mu^2}{b},
\end{align}
where the second inequality holds by $\mathcal{I}(|S_1|<n)=0$, and the inequality $k-c_kq \leq q$ from $c_kq \leq k \leq (c_k+1)q-1$.

It follows that
\begin{align}
\mathbb{E}\| v_k - \nabla f(x_k) \|^2 & = \mathbb{E}\| v_k - \hat{\nabla}_{\texttt{coo}}f(x_k) + \hat{\nabla}_{\texttt{coo}}f(x_k) -  \nabla f(x_k) \|^2 \nonumber \\
 & \leq  2\mathbb{E}\| v_k - \hat{\nabla}_{\texttt{coo}}f(x_k)\|^2 + 2\mathbb{E}\| \hat{\nabla}_{\texttt{coo}}f(x_k) - \nabla f(x_k) \|^2 \nonumber \\
 & \leq \frac{6L^2}{b} \sum_{i=c_kq}^{k-1} \mathbb{E}\|x_{i+1}-x_i\|^2 + \frac{12qL^2d\mu^2}{b} + 2 L^2d\mu^2,
\end{align}
where the second inequality holds by the above inequality \eqref{eq:B1}, and Lemmas \ref{lem:A1}.
\end{proof}

\begin{lemma} \label{lem:B2}
Suppose the sequence $\{x_k,y_{[m]}^k,\lambda_k\}_{k=1}^K$ be generated from Algorithm \ref{alg:1}, it holds that
\begin{align}
\mathbb{E}\|\lambda_{k+1} - \lambda_k\|^2 \leq & \frac{60L^2}{b\sigma^A_{\min}} \sum_{i=c_kq}^{k-1} \mathbb{E}\|x_{i+1}-x_i\|^2
+ (\frac{5L^2}{\sigma^A_{\min}}+\frac{5\sigma^2_{\max}(G)}{\sigma^A_{\min}\eta^2})\|x_k-x_{k-1}\|^2 \nonumber \\
& + \frac{5\sigma^2_{\max}(G)}{\sigma^A_{\min}\eta^2}\|x_{k+1}-x_{k}\|^2  + \frac{120qL^2d\mu^2}{b\sigma^A_{\min}} + \frac{20 L^2d\mu^2}{\sigma^A_{\min}}.
\end{align}
\end{lemma}
\begin{proof}
By using the optimal condition of the step 10 in Algorithm \ref{alg:1},
we have
\begin{align}
  v_k + \frac{G}{\eta}(x_{k+1}-x_k) - A^T\lambda_k + \rho A^T(Ax_{k+1}+\sum_{j=1}^mB_jy_j^{k+1}-c) = 0.
\end{align}
By using the step 11 of Algorithm \ref{alg:1}, we have
\begin{align}
 A^T\lambda_{k+1} = v_k + \frac{G}{\eta}(x_{k+1}-x_k).
\end{align}
It follows that
\begin{align} \label{eq:B2}
 \lambda_{k+1} = (A^T)^+ \big( v_k + \frac{G}{\eta}(x_{k+1}-x_k) \big),
\end{align}
where $(A^T)^+$ is the pseudoinverse of $A^T$. By Assumption 4, without loss of generality, we use the full column rank matrix $A$.
It is easily verified that $(A^T)^+=A(A^TA)^{-1}$.
By using the above formula \eqref{eq:B2}, we have
\begin{align}
\mathbb{E}\|\lambda_{k+1}-\lambda_k\|^2 & = \mathbb{E} \|(A^T)^{+}\big( v_k + \frac{G}{\eta}(x_{k+1}-x_k) - v_{k-1} - \frac{G}{\eta}(x_{k}-x_{k-1})\big)\|^2  \nonumber \\
& \leq \frac{1}{\sigma^A_{\min}} \mathbb{E} \|v_k + \frac{G}{\eta}(x_{k+1}-x_k) - v_{k-1} - \frac{G}{\eta}(x_{k}-x_{k-1})\|^2  \nonumber \\
& = \frac{1}{\sigma^A_{\min}} \mathbb{E} \|v_k - \nabla f(x_k) + \nabla f(x_k) -  \nabla f(x_{k-1}) + \nabla f(x_{k-1}) - v_{k-1} \nonumber \\
& \qquad + \frac{G}{\eta}(x_{k+1}-x_k) - \frac{G}{\eta}(x_{k}-x_{k-1})\|^2  \nonumber \\
& \leq \frac{1}{\sigma^A_{\min}}\big[5 \mathbb{E}\|v_k - \nabla f(x_k)\|^2 +  5\mathbb{E}\|\nabla f(x_k) - \nabla f(x_{k-1})\|^2 + 5 \mathbb{E}\|v_{k-1} - \nabla f(x_{k-1})\|^2 \nonumber\\
& \qquad +\frac{5\sigma^2_{\max}(G)}{\eta^2}\mathbb{E}\|x_{k+1}-x_k\|^2 + \frac{5\sigma^2_{\max}(G)}{\eta^2}\mathbb{E}\|x_{k}-x_{k-1}\|^2 \big] \nonumber \\
& \leq \frac{1}{\sigma^A_{\min}}\big[\frac{30L^2}{b} \sum_{i=c_kq}^{k-1} \mathbb{E}\|x_{i+1}-x_i\|^2 + \frac{120qL^2d\mu^2}{b} + 20 L^2d\mu^2 + \frac{30L^2}{b} \sum_{i=c_kq}^{k-2} \mathbb{E}\|x_{i+1}-x_i\|^2  \nonumber\\
& \qquad  +  5L^2\mathbb{E}\|x_k - x_{k-1}\|^2+\frac{5\sigma^2_{\max}(G)}{\eta^2}\mathbb{E}\|x_{k+1}-x_k\|^2 + \frac{5\sigma^2_{\max}(G)}{\eta^2}\mathbb{E}\|x_{k}-x_{k-1}\|^2 \big] \nonumber \\
& \leq  \frac{1}{\sigma^A_{\min}}\big[\frac{60L^2}{b} \sum_{i=c_kq}^{k-1} \mathbb{E}\|x_{i+1}-x_i\|^2 + \frac{120qL^2d\mu^2}{b} + 20 L^2d\mu^2 +  5L^2\mathbb{E}\|x_k - x_{k-1}\|^2  \nonumber\\
& \qquad +\frac{5\sigma^2_{\max}(G)}{\eta^2}\mathbb{E}\|x_{k+1}-x_k\|^2 + \frac{5\sigma^2_{\max}(G)}{\eta^2}\mathbb{E}\|x_{k}-x_{k-1}\|^2 \big],
\end{align}
where the first inequality follows by $((A^T)^+)^T(A^T)^+=(A(A^TA)^{-1})^TA(A^TA)^{-1}=(A^TA)^{-1}$, and the third inequality follows by Lemma \ref{lem:B1}.

\end{proof}

\begin{lemma} \label{lem:B3}
Suppose the sequence $\{x_k,y_{[m]}^k,\lambda_k\}_{k=1}^K$ be generated from Algorithm \ref{alg:1},
and define a \emph{Lyapunov} function $\Omega_k$ as follows:
\begin{align}
 \Omega_k = \mathbb{E}\big[\mathcal{L}_{\rho} (x_k,y_{[m]}^k,\lambda_k) + (\frac{5L^2}{\sigma^A_{\min}\rho}+\frac{5\sigma^2_{\max}(G)}{\sigma^A_{\min}\eta^2\rho})\|x_k-x_{k-1}\|^2
 + \frac{12L^2}{\sigma^A_{\min}\rho b} \sum_{i=c_kq}^{k-1}\|x_{i+1}-x_i\|^2 \big]. \nonumber
\end{align}
Let $b=q$, $\eta=\frac{\alpha\sigma_{\min}(G)}{4L} \ (0<\alpha \leq 1)$ and $\rho = \frac{2\sqrt{237}\kappa_GL}{\sigma^A_{\min}\alpha}$,
we have
\begin{align}
\frac{1}{K}\sum_{k=0}^{K-1} (\|x_{k+1} - x_{k}\|^2 + \sum_{j=1}^m \|y_j^k-y_j^{k+1}\|^2) \leq \frac{\Omega_{0} - \Omega^*}{K\gamma}
+ \frac{7\vartheta_1L^2d\mu^2}{\gamma}, \nonumber
\end{align}
where $\vartheta_1 = \frac{1}{L} + \frac{20}{\sigma^A_{\min}\rho}$, $\gamma = \min(\chi,\sigma_{\min}^H)$ with $\chi\geq \frac{\sqrt{237}\kappa_GL}{2\alpha}$
and $\Omega^*$ is a lower bound of the function $\Omega_k$.
\end{lemma}
\begin{proof}
By the optimal condition of step 9 in Algorithm \ref{alg:1},
we have, for $j\in [m]$
\begin{align}
& 0 =(y_j^k-y_j^{k+1})^T\big(\partial \psi_j(y_j^{k+1}) - B_j^T\lambda_k + \rho B_j^T(Ax_k + \sum_{i=1}^jB_iy_i^{k+1} + \sum_{i=j+1}^mB_iy_i^{k}-c) + H_j(y_j^{k+1}-y_j^k)\big) \nonumber \\
& \leq \psi_j(y_j^k) \!-\! \psi_j(y_j^{k+1}) \!-\! \lambda_k^T(B_jy_j^k-B_jy_j^{k+1}) \!+\! \rho(B_jy_j^k-B_jy_j^{k+1})^T(Ax_k + \sum_{i=1}^jB_iy_i^{k+1} + \sum_{i=j+1}^mB_iy_i^{k}-c) \!-\! \|y_j^{k+1}-y_j^k\|^2_{H_j} \nonumber \\
& = \psi_j(y_j^k)- \psi_j(y_j^{k+1}) - (\lambda_k)^T(Ax_k+\sum_{i=1}^{j-1}B_iy_i^{k+1} + \sum_{i=j}^mB_iy_i^{k}-c) + (\lambda_k)^T(Ax_k+\sum_{i=1}^jB_iy_i^{k+1}+ \sum_{i=j+1}^mB_iy_i^{k}-c) \nonumber \\
&  + \frac{\rho}{2}\|Ax_k \!+\! \sum_{i=1}^{j-1}B_iy_i^{k+1} \!+\! \sum_{i=j}^mB_iy_i^{k}-c\|^2 \!-\! \frac{\rho}{2}\|Ax_k \!+\!\sum_{i=1}^jB_iy_i^{k+1}\!+\! \sum_{i=j+1}^mB_iy_i^{k}-c\|^2 \!-\!\frac{\rho}{2}\|B_jy_j^k-B_jy_j^{k+1}\|^2 \!-\! \|y_j^{k+1}-y_j^k\|^2_{H_j} \nonumber \\
& =  \big( \underbrace{ f(x_k)\! +\! \sum_{i=1}^{j-1}\psi_i(y_i^{k+1}) \!+\! \sum_{i=j}^{m}\psi_i(y_i^{k}) \!-\! \lambda_k^T(Ax_k \!+\! \sum_{i=1}^{j-1}B_iy_i^{k+1} \!+\! \sum_{i=j}^mB_iy_i^{k}-c) \!+\! \frac{\rho}{2}\|Ax_k \!+\!\sum_{i=1}^{j-1}B_iy_i^{k+1} \!+\! \sum_{i=j}^mB_iy_i^{k}-c\|^2}_{=\mathcal{L}_{\rho} (x_k,y_{[j-1]}^{k+1},y_{[j:m]}^k,\lambda_k)} \big)\nonumber \\
& - \big( \underbrace{ f(x_k) \!+\! \sum_{i=1}^{j}\psi_i(y_i^{k+1}) \!+\! \sum_{i=j+1}^{m}\psi_i(y_i^{k}) \!-\! \lambda_k^T(Ax_k\!+\!\sum_{i=1}^jB_iy_i^{k+1}\!+\! \sum_{i=j+1}^mB_iy_i^{k}-c) \!+\! \frac{\rho}{2}\|Ax_k\!+\!\sum_{i=1}^jB_iy_i^{k+1} \!+\! \sum_{i=j+1}^mB_iy_i^{k}-c\|^2}_{=\mathcal{L}_{\rho} (x_k,y_{[j]}^{k+1},y_{[j+1:m]}^k,\lambda_k)} \big) \nonumber \\
&  -\frac{\rho}{2}\|B_jy_j^k-B_jy_j^{k+1}\|^2 - \|y_j^{k+1}-y_j^k\|^2_{H_j} \nonumber \\
& \leq \mathcal{L}_{\rho} (x_k,y_{[j-1]}^{k+1},y_{[j:m]}^k,\lambda_k) - \mathcal{L}_{\rho} (x_k,y_{[j]}^{k+1},y_{[j+1:m]}^k,\lambda_k)
- \sigma_{\min}(H_j)\|y_j^k-y_j^{k+1}\|^2,
\end{align}
where the first inequality holds by the convexity of function $\psi_j(y)$,
and the second equality follows by applying the equality
$(a-b)^Tb = \frac{1}{2}(\|a\|^2-\|b\|^2-\|a-b\|^2)$ on the term $(B_jy_j^k-B_jy_j^{k+1})^T(Ax_k + \sum_{i=1}^jB_iy_i^{k+1} + \sum_{i=j+1}^mB_iy_i^{k}-c)$.
Thus, we have, for all $j\in[m]$
\begin{align} \label{eq:C1}
\mathcal{L}_{\rho} (x_k,y_{[j]}^{k+1},y_{[j+1:m]}^k,\lambda_k) \leq \mathcal{L}_{\rho} (x_k,y_{[j-1]}^{k+1},y_{[j:m]}^k,\lambda_k)
- \sigma_{\min}(H_j)\|y_j^k-y_j^{k+1}\|^2.
\end{align}
Telescoping \eqref{eq:C1} over $j$ from $1$ to $m$, we obtain
\begin{align} \label{eq:C2}
\mathcal{L}_{\rho} (x_k,y^{k+1}_{[m]},\lambda_k) \leq \mathcal{L}_{\rho} (x_k,y^k_{[m]},\lambda_k)
- \sigma_{\min}^H\sum_{j=1}^m \|y_j^k-y_j^{k+1}\|^2,
\end{align}
where $\sigma_{\min}^H=\min_{j\in[m]}\sigma_{\min}(H_j)$.

By Assumption 1, we have
\begin{align} \label{eq:C3}
0 \leq f(x_k) - f(x_{k+1}) + \nabla f(x_k)^T(x_{k+1}-x_k) + \frac{L}{2}\|x_{k+1}-x_k\|^2.
\end{align}
By using the optimal condition of step 10 in Algorithm \ref{alg:1},
we have
\begin{align} \label{eq:C4}
0 = (x_k-x_{k+1})^T \big( v_k - A^T\lambda_k + \rho A^T(Ax_{k+1} + \sum_{j=1}^mB_jy_j^{k+1}-c) + \frac{G}{\eta}(x_{k+1}-x_k) \big).
\end{align}
Combining \eqref{eq:C3} and \eqref{eq:C4}, we have
\begin{align}
0 & \leq f(x_k) - f(x_{k+1}) + \nabla f(x_k)^T(x_{k+1}-x_k) + \frac{L}{2}\|x_{k+1}-x_k\|^2 \nonumber \\
& \quad + (x_k-x_{k+1})^T \big( v_k - A^T\lambda_k + \rho A^T(Ax_{k+1} + \sum_{j=1}^mB_jy_j^{k+1}-c) + \frac{G}{\eta}(x_{k+1}-x_k) \big)  \nonumber \\
& = f(x_k) - f(x_{k+1}) + \frac{L}{2}\|x_k-x_{k+1}\|^2 - \frac{1}{\eta}\|x_k - x_{k+1}\|^2_G + (x_k-x_{k+1})^T(v_k-\nabla f(x_k)) \nonumber \\
& \quad -(\lambda_k)^T(Ax_k-Ax_{k+1}) + \rho(Ax_k - Ax_{k+1})^T(Ax_{k+1} + \sum_{j=1}^mB_jy_j^{k+1}-c) \nonumber \\
& = f(x_k) \!-\! f(x_{k+1}) \!+\! \frac{L}{2}\|x_k-x_{k+1}\|^2 \!-\! \frac{1}{\eta}\|x_k - x_{k+1}\|^2_G
\!+\! (x_k-x_{k+1})^T\big(v_k-\nabla f(x_k)\big) \!-\!\lambda_k^T(Ax_k + \sum_{j=1}^mB_jy_j^{k+1}-c)\nonumber \\
& \quad  + \lambda_k^T(Ax_{k+1} \!+\! \sum_{j=1}^mB_jy_j^{k+1}\!-\!c) +\frac{\rho}{2}\big(\|Ax_k \!+\! \sum_{j=1}^mB_jy_j^{k+1}\!-\!c\|^2
\!-\! \|Ax_{k+1} \!+\! \sum_{j=1}^mB_jy_j^{k+1}\!-\!c\|^2 \!- \!\|Ax_k \!-\! Ax_{k+1}\|^2 \big) \nonumber \\
& = \underbrace{f(x_k) + \sum_{j=1}^{m}\psi_j(y_j^{k+1}) -\lambda_k^T(Ax_k + \sum_{j=1}^mB_jy_j^{k+1}-c) + \frac{\rho}{2}\|Ax_{k} + \sum_{j=1}^mB_jy_j^{k+1}-c\|^2}_{=\mathcal{L}_{\rho}(x_k,y_{[m]}^{k+1},\lambda_k)}  \nonumber \\
& \quad - \underbrace{\big(f(x_{k+1}) + \sum_{j=1}^{m}\psi_j(y_j^{k+1})-  \lambda_k^T(Ax_{k+1}+ \sum_{j=1}^mB_jy_j^{k+1}-c) + \frac{\rho}{2} \|Ax_{k+1} + \sum_{j=1}^mB_jy_j^{k+1}-c\|^2
\big)}_{=\mathcal{L}_{\rho}(x_{k+1},y^{k+1}_{[m]},z_{k})} \nonumber \\
& \quad  + \frac{L}{2}\|x_k-x_{k+1}\|^2 + (x_k-x_{k+1})^T\big(v_k-\nabla f(x_k)\big)  - \frac{1}{\eta}\|x_k - x_{k+1}\|^2_G - \frac{\rho}{2}\|Ax_k - Ax_{k+1}\|^2 \nonumber \\
& \leq \mathcal{L}_{\rho} (x_k,y_{[m]}^{k+1},\lambda_k) -  \mathcal{L}_{\rho} (x_{k+1},y_{[m]}^{k+1},\lambda_k)
- (\frac{\sigma_{\min}(G)}{\eta}+\frac{\rho \sigma^A_{\min}}{2}-\frac{L}{2}) \|x_{k+1} - x_{k}\|^2 + (x_k-x_{k+1})^T(v_k-\nabla f(x_k)) \nonumber \\
& \leq \mathcal{L}_{\rho} (x_k,y_{[m]}^{k+1},\lambda_k) -  \mathcal{L}_{\rho} (x_{k+1},y_{[m]}^{k+1},\lambda_k)
- (\frac{\sigma_{\min}(G)}{\eta}+\frac{\rho \sigma^A_{\min}}{2}-L ) \|x_{k+1} - x_{k}\|^2 + \frac{1}{2L}\|v_k - \nabla f(x_k)\|^2 \nonumber \\
& \leq \mathcal{L}_{\rho} (x_k,y_{[m]}^{k+1},\lambda_k) -  \mathcal{L}_{\rho} (x_{k+1},y_{[m]}^{k+1},\lambda_k)
- (\frac{\sigma_{\min}(G)}{\eta}+\frac{\rho \sigma^A_{\min}}{2}-L ) \|x_{k+1} - x_{k}\|^2 + \frac{3L}{b} \sum_{i=c_kq}^{k-1} \mathbb{E}\|x_{i+1}-x_i\|^2 \nonumber \\
& \quad + \frac{6qLd\mu^2}{b} + Ld\mu^2, \nonumber
\end{align}
where the second equality follows by applying the equality
$(a-b)^Tb = \frac{1}{2}(\|a\|^2-\|b\|^2-\|a-b\|^2)$ over the term $(Ax_k - Ax_{k+1})^T(Ax_{k+1}+\sum_{j=1}^mB_jy_j^{k+1}-c)$; the third inequality
follows by the inequality $a^Tb \leq \frac{1}{2L}\|a\|^2 + \frac{L}{2}\|b\|^2$, and the forth inequality holds by Lemma \ref{lem:B1}.
It follows that
\begin{align} \label{eq:C5}
\mathcal{L}_{\rho} (x_{k+1},y_{[m]}^{k+1},\lambda_k) & \leq \mathcal{L}_{\rho} (x_k,y_{[m]}^{k+1},\lambda_k)
- (\frac{\sigma_{\min}(G)}{\eta}+\frac{\rho \sigma^A_{\min}}{2}-L ) \|x_{k+1} - x_{k}\|^2 \nonumber \\
& \quad + \frac{3L}{b} \sum_{i=c_kq}^{k-1} \mathbb{E}\|x_{i+1}-x_i\|^2 + \frac{6qLd\mu^2}{b} + Ld\mu^2.
\end{align}

By using the step 11 in Algorithm \ref{alg:1}, we have
\begin{align} \label{eq:C6}
\mathcal{L}_{\rho} (x_{k+1},y_{[m]}^{k+1},\lambda_{k+1}) -\mathcal{L}_{\rho} (x_{k+1},y_{[m]}^{k+1},\lambda_k)
& = \frac{1}{\rho}\|\lambda_{k+1}-\lambda_k\|^2 \nonumber \\
& \leq \frac{60L^2}{b\sigma^A_{\min}\rho} \sum_{i=c_kq}^{k-1} \mathbb{E}\|x_{i+1}-x_i\|^2
+ (\frac{5L^2}{\sigma^A_{\min}\rho}+\frac{5\sigma^2_{\max}(G)}{\sigma^A_{\min}\eta^2\rho})\|x_k-x_{k-1}\|^2 \nonumber \\
& \quad + \frac{5\sigma^2_{\max}(G)}{\sigma^A_{\min}\eta^2\rho}\|x_{k+1}-x_{k}\|^2  + \frac{120qL^2d\mu^2}{b\sigma^A_{\min}\rho} + \frac{20 L^2d\mu^2}{\sigma^A_{\min}\rho},
\end{align}
where the above inequality holds by Lemma \ref{lem:B2}.

Combining \eqref{eq:C2}, \eqref{eq:C5} and \eqref{eq:C6}, we have
\begin{align} \label{eq:C7}
\mathcal{L}_{\rho} (x_{k+1},y_{[m]}^{k+1},\lambda_{k+1}) & \leq \mathcal{L}_{\rho} (x_k,y_{[m]}^k,\lambda_k)
 - \sigma_{\min}^H\sum_{j=1}^m \|y_j^k-y_j^{k+1}\|^2 + (\frac{3L}{b} + \frac{60L^2}{\sigma^A_{\min}b\rho})\sum_{i=c_kq}^{k-1}\mathbb{E}\|x_{i+1}-x_i\|^2 \nonumber \\
& \quad + (\frac{5L^2}{\sigma^A_{\min}\rho}+\frac{5\sigma^2_{\max}(G)}{\sigma^A_{\min}\eta^2\rho})\|x_k-x_{k-1}\|^2 + \big( \frac{6qL}{b} + L + \frac{120qL^2}{b\sigma^A_{\min}\rho}
+ \frac{20 L^2}{\sigma^A_{\min}\rho} \big) d\mu^2\nonumber \\
& \quad - (\frac{\sigma_{\min}(G)}{\eta}+\frac{\rho \sigma^A_{\min}}{2}-L - \frac{5\sigma^2_{\max}(G)}{\sigma^A_{\min}\eta^2\rho}) \|x_{k+1} - x_{k}\|^2.
\end{align}

Next, we define a useful \emph{Lyapunov} function $\Omega_k$ as follows:
\begin{align}
\Omega_k = \mathbb{E}\big[\mathcal{L}_{\rho} (x_k,y_{[m]}^k,\lambda_k) + (\frac{5L^2}{\sigma^A_{\min}\rho}+\frac{5\sigma^2_{\max}(G)}{\sigma^A_{\min}\eta^2\rho})\|x_k-x_{k-1}\|^2 + \frac{12L^2}{\sigma^A_{\min}\rho b} \sum_{i=c_kq}^{k-1}\|x_{i+1}-x_i\|^2 \big].
\end{align}
It follows that
\begin{align} \label{eq:C8}
\Omega_{k+1}  &= \mathbb{E} \big[ \mathcal{L}_{\rho} (x_{k+1},y_{[m]}^{k+1},\lambda_{k+1}) + (\frac{5L^2}{\sigma^A_{\min}\rho} + \frac{5\sigma^2_{\max}(G)}{\sigma^A_{\min}\eta^2\rho})\|x_{k+1}-x_{k}\|^2
+ \frac{12L^2}{\sigma^A_{\min}\rho b} \sum_{i=c_kq}^{k}\|x_{i+1}-x_i\|^2 \big]\nonumber \\
& \leq \mathbb{E} \big[ \mathcal{L}_{\rho} (x_k,y_{[m]}^k,\lambda_k) + (\frac{5L^2}{\sigma^A_{\min}\rho}+\frac{5\sigma^2_{\max}(G)}{\sigma^A_{\min}\eta^2\rho})\|x_{k}-x_{k-1}\|^2
+ \frac{12L^2}{\sigma^A_{\min}\rho b} \sum_{i=c_kq}^{k-1}\|x_{i+1}-x_i\|^2 \big] \nonumber \\
& \quad - \sigma_{\min}^H\sum_{j=1}^m \|y_j^k-y_j^{k+1}\|^2 - (\frac{\sigma_{\min}(G)}{\eta}+\frac{\rho \sigma^A_{\min}}{2}-L - \frac{10\sigma^2_{\max}(G)}{\sigma^A_{\min}\eta^2\rho}-\frac{5L^2}{\sigma^A_{\min}\rho}
 - \frac{12L^2}{\sigma^A_{\min}\rho b})\mathbb{E} \|x_{k+1} - x_{k}\|^2 \nonumber \\
& \quad + (\frac{3L}{b} + \frac{60L^2}{\sigma^A_{\min}b\rho})\sum_{i=c_kq}^{k-1}\mathbb{E}\|x_{i+1}-x_i\|^2 + \big( \frac{6qL}{b} + L + \frac{120qL^2}{b\sigma^A_{\min}\rho}
+ \frac{20 L^2}{\sigma^A_{\min}\rho} \big) d\mu^2 \nonumber \\
& \leq \Omega_k - \sigma_{\min}^H\sum_{j=1}^m \|y_j^k-y_j^{k+1}\|^2 - \big(\frac{\sigma_{\min}(G)}{\eta}+\frac{\rho \sigma^A_{\min}}{2}- L - \frac{10\sigma^2_{\max}(G)}{\sigma^A_{\min}\eta^2\rho}-\frac{5L^2}{\sigma^A_{\min}\rho}
 - \frac{12L^2}{\sigma^A_{\min}\rho b}\big)\mathbb{E}\|x_{k+1} - x_{k}\|^2 \nonumber \\
& \quad+ (\frac{3L}{b} + \frac{60L^2}{\sigma^A_{\min}b\rho})\sum_{i=c_kq}^{k-1}\mathbb{E}\|x_{i+1}-x_i\|^2 + \big( \frac{6qL}{b} + L + \frac{120qL^2}{b\sigma^A_{\min}\rho}
+ \frac{20L^2}{\sigma^A_{\min}\rho} \big) d\mu^2 ,
\end{align}
where the first inequality holds by the inequality \eqref{eq:C7} and the equality
$$\sum_{i=c_kq}^{k}\mathbb{E}\|x_{i+1}-x_i\|^2 = \sum_{i=c_kq}^{k-1}\mathbb{E}\|x_{i+1}-x_i\|^2 + \mathbb{E}\|x_{k+1}-x_k\|^2.$$

Since $c_kq \leq k \leq (c_k+1)q -1$, and let $c_kq \leq t \leq (c_k+1)q-1$, then telescoping inequality \eqref{eq:C8} over $k$ from $c_kq$ to $k$, we have
\begin{align} \label{eq:C9}
\Omega_{k+1} &\leq \Omega_{c_kq}
 - (\frac{\sigma_{\min}(G)}{\eta}+\frac{\rho \sigma^A_{\min}}{2}-L - \frac{10\sigma^2_{\max}(G)}{\sigma^A_{\min}\eta^2\rho}-\frac{5L^2}{\sigma^A_{\min}\rho}
  -\frac{12L^2}{\sigma^A_{\min}\rho b} ) \sum_{t=c_kq}^k \mathbb{E}\|x_{t+1} - x_{t}\|^2 \nonumber \\
& \quad- \sigma_{\min}^H\sum_{t=c_kq}^{k}\sum_{j=1}^m \|y_j^t\!-\!y_j^{t+1}\|^2 + (\frac{3L}{b} \!+\! \frac{60L^2}{\sigma^A_{\min}b\rho})\sum_{t=c_kq}^{k}\sum_{i=c_kq}^{k-1}\mathbb{E}\|x_{i+1}-x_i\|^2
+ \sum_{t=c_kq}^{k}\big( \frac{6qL}{b} \!+\! L \!+\! \frac{120qL^2}{b\sigma^A_{\min}\rho}
\!+\! \frac{20 L^2}{\sigma^A_{\min}\rho} \big) d\mu^2 \nonumber \\
& \leq \Omega_{c_kq} - (\frac{\sigma_{\min}(G)}{\eta}+\frac{\rho \sigma^A_{\min}}{2}-L - \frac{10\sigma^2_{\max}(G)}{\sigma^A_{\min}\eta^2\rho}-\frac{5L^2}{\sigma^A_{\min}\rho}
 -\frac{12L^2}{\sigma^A_{\min}\rho b}) \sum_{i = c_kq}^k \mathbb{E}\|x_{i+1} - x_{i}\|^2 \nonumber \\
& \quad - \sigma_{\min}^H\sum_{i=c_kq}^{k-1}\sum_{j=1}^m \|y_j^i-y_j^{i+1}\|^2 + (\frac{3Lq}{b} + \frac{60L^2q}{\sigma^A_{\min}b\rho})\sum_{i=c_kq}^{k}\mathbb{E}\|x_{i+1}-x_i\|^2 + \big( \frac{6qL}{b}
+ L + \frac{120qL^2}{b\sigma^A_{\min}\rho} + \frac{20 L^2}{\sigma^A_{\min}\rho} \big) qd\mu^2 \nonumber \\
& = \Omega_{c_kq} - \underbrace{ (\frac{\sigma_{\min}(G)}{\eta}+\frac{\rho \sigma^A_{\min}}{2}-L - \frac{10\sigma^2_{\max}(G)}{\sigma^A_{\min}\eta^2\rho}-\frac{5L^2}{\sigma^A_{\min}\rho}
- \frac{12L^2}{\sigma^A_{\min}\rho b}-\frac{3Lq}{b} - \frac{60L^2q}{\sigma^A_{\min}b\rho} )}_{\chi}  \sum_{i=c_kq}^k \|x_{i+1} - x_{i}\|^2 \nonumber \\
& \quad - \sigma_{\min}^H\sum_{i=c_kq}^{k-1}\sum_{j=1}^m \|y_j^i-y_j^{i+1}\|^2 + \big( \frac{6qL}{b} + L + \frac{120qL^2}{b\sigma^A_{\min}\rho}
+ \frac{20L^2}{\sigma^A_{\min}\rho} \big) qd\mu^2,
\end{align}
where the second inequality holds by the fact that
\begin{align}
\sum_{j=c_kq}^{k}\sum_{i=c_kq}^{k-1}\mathbb{E}\|x_{i+1}-x_i\|^2 \leq \sum_{j=c_kq}^{k}\sum_{i=c_kq}^{k}\mathbb{E}\|x_{i+1}-x_i\|^2 \leq q\sum_{i=c_kq}^{k}\mathbb{E}\|x_{i+1}-x_i\|^2. \nonumber
\end{align}

Since $b=q$, we have
\begin{align}
\chi & = \frac{\sigma_{\min}(G)}{\eta}+\frac{\rho \sigma^A_{\min}}{2}-L - \frac{10\sigma^2_{\max}(G)}{\sigma^A_{\min}\eta^2\rho}
-\frac{5L^2}{\sigma^A_{\min}\rho} -  \frac{12L^2}{\sigma^A_{\min}\rho b} -\frac{3Lq}{b} - \frac{60L^2q}{\sigma^A_{\min}b\rho} \nonumber \\
& = \underbrace{\frac{\sigma_{\min}(G)}{\eta} - 4L }_{=T_1} + \underbrace{\frac{\rho\sigma^A_{\min}}{2} - \frac{10\sigma^2_{\max}(G)}{\sigma^A_{\min}\eta^2\rho}
-\frac{65L^2}{\sigma^A_{\min}\rho}-  \frac{12L^2}{\sigma^A_{\min}\rho b}}_{=T_2}.
\end{align}
Given $0< \eta \leq \frac{\sigma_{\min}(G)}{4L}$, we have $T_1\geq 0$. Further, let $\eta = \frac{\alpha\sigma_{\min}(G)}{4L} \ (0<\alpha \leq 1)$ and
$\rho = \frac{2\sqrt{237}\kappa_GL}{\sigma^A_{\min}\alpha}$,
we have
\begin{align}
T_2 & = \frac{\rho \sigma^A_{\min}}{2} - \frac{10\sigma^2_{\max}(G)}{\sigma^A_{\min}\eta^2\rho}
- \frac{65L^2}{\sigma^A_{\min}\rho} - \frac{12L^2}{\sigma^A_{\min}\rho b}\nonumber \\
& =  \frac{\rho \sigma^A_{\min}}{2} - \frac{160L^2\kappa^2_G}{\sigma^A_{\min}\rho\alpha^2}
-\frac{65L^2}{\sigma^A_{\min}\rho} - \frac{12L^2}{\sigma^A_{\min}\rho b}\nonumber \\
& \geq \frac{\rho \sigma^A_{\min}}{2} - \frac{160L^2\kappa^2_G}{\sigma^A_{\min}\rho\alpha^2}
-\frac{65L^2\kappa^2_G}{\sigma^A_{\min}\rho\alpha^2} -  \frac{12\kappa^2_GL^2}{\sigma^A_{\min}\rho \alpha^2} \nonumber \\
& = \frac{\rho \sigma^A_{\min}}{4} + \underbrace{\frac{\rho \sigma^A_{\min}}{4} - \frac{237L^2\kappa^2_G }{\sigma^A_{\min}\rho\alpha^2}}_{\geq 0} \nonumber \\
& \geq \frac{\sqrt{237}\kappa_GL}{2\alpha},
\end{align}
where the first inequality holds by $\kappa_G \geq 1$ and
$b \geq 1 \geq \alpha^2$ and the second inequality holds by $\rho = \frac{2\sqrt{237}\kappa_GL}{\sigma^A_{\min}\alpha}$.
Thus, we obtain $\chi \geq \frac{\sqrt{237}\kappa_GL}{2\alpha}$.

By Assumption 4, i.e., $A$ is a full column rank matrix,
we have $(A^T)^+ = A(A^T A)^{-1}$.
It follows that $\sigma_{\max}((A^T)^+)^T(A^T)^+) = \sigma_{\max}((A^TA)^{-1}) = \frac{1}{\sigma_{\min}^A}$.
By using the above formula \eqref{eq:B2}, then we have
\begin{align}
&\mathcal{L}_{\rho} (x_{k+1},y_{[m]}^{k+1},\lambda_{k+1}) \nonumber \\
& = f(x_{k+1}) + \sum_{j=1}^m \psi_j(y_j^{k+1}) - \lambda_{k+1}^T(Ax_{k+1} + \sum_{j=1}^mB_jy_j^{k+1} -c) + \frac{\rho}{2}\|Ax_{k+1} + \sum_{j=1}^mB_jy_j^{k+1} -c\|^2 \nonumber \\
& = f(x_{k+1}) + \sum_{j=1}^m \psi_j(y_j^{k+1}) -  \langle(A^T)^+(v_k + \frac{G}{\eta}(x_{k+1}-x_k)), Ax_{k+1} + \sum_{j=1}^mB_jy_j^{k+1} -c\rangle + \frac{\rho}{2}\|Ax_{k+1} + \sum_{j=1}^mB_jy_j^{k+1} -c\|^2 \nonumber \\
& = f(x_{k+1}) + \sum_{j=1}^m \psi_j(y_j^{k+1}) - \langle(A^T)^+(v_{k} - \nabla f(x_{k}) + \nabla f(x_{k})+ \frac{G}{\eta}(x_{k+1}-x_k)), Ax_{k+1} + \sum_{j=1}^mB_jy_j^{k+1} -c\rangle  \nonumber \\
& \quad +  \frac{\rho}{2}\|Ax_{k+1} + \sum_{j=1}^mB_jy_j^{k+1} -c\|^2 \nonumber \\
& \geq f(x_{k+1}) + \sum_{j=1}^m \psi_j(y_j^{k+1}) - \frac{2}{\sigma^A_{\min}\rho}\|v_k - \nabla f(x_{k})\|^2 - \frac{2}{\sigma^A_{\min}\rho}\|\nabla f(x_{k})\|^2
- \frac{2\sigma^2_{\max}(G)}{\sigma^A_{\min}\eta^2\rho}\|x_{k+1}-x_k\|^2 \nonumber \\
& \quad + \frac{\rho}{8}\|Ax_{k+1} + \sum_{j=1}^mB_jy_j^{k+1} -c\|^2 \nonumber\\
& \geq f(x_{k+1}) - \frac{4}{\sigma^A_{\min}\rho}\|\nabla f(x_{k+1})\|^2 + \sum_{j=1}^m \psi_j(y_j^{k+1})  - \frac{4}{\sigma^A_{\min}\rho}\|\nabla f(x_{k})-\nabla f(x_{k+1})\|^2  \nonumber \\
& \quad - \frac{12L^2}{\sigma^A_{\min}\rho b} \sum_{i=c_kq}^{k-1} \mathbb{E}\|x_{i+1}-x_i\|^2 - \frac{24qL^2d\mu^2}{\sigma^A_{\min}\rho b}
- \frac{4 L^2d\mu^2}{\sigma^A_{\min}\rho} - \frac{2\sigma^2_{\max}(G)}{\sigma^A_{\min}\eta^2\rho}\|x_{k+1}-x_k\|^2 \nonumber \\
& \geq f^* + \sum_{j=1}^m \psi_j^*  - \frac{12L^2}{\sigma^A_{\min}\rho b} \sum_{i=c_kq}^{k-1} \mathbb{E}\|x_{i+1}-x_i\|^2 - \frac{24qL^2d\mu^2}{\sigma^A_{\min}\rho b}
- \frac{4 L^2d\mu^2}{\sigma^A_{\min}\rho} -  \big(\frac{2\sigma^2_{\max}(G)}{\sigma^A_{\min}\eta^2\rho} + \frac{4L^2}{\sigma^A_{\min}\rho} \big)\|x_{k+1}-x_k\|^2,
\end{align}
where the first inequality is obtained by applying $ \langle a, b\rangle \leq \frac{1}{2\tau}\|a\|^2 + \frac{\tau}{2}\|b\|^2$ to the terms
$\langle(A^T)^+(v_{k} - \nabla f(x_{k})), Ax_{k+1} + \sum_{j=1}^mB_jy_j^{k+1} -c\rangle$, $\langle(A^T)^+v_{k}, Ax_{k+1} + \sum_{j=1}^mB_jy_j^{k+1} -c\rangle $ and
$\langle(A^T)^+\frac{G}{\eta}(x_{k+1}-x_k), Ax_{k+1} + \sum_{j=1}^mB_jy_j^{k+1} -c\rangle$ with $\tau = \frac{\rho}{4}$, respectively;
the second inequality follows by Lemma \ref{lem:B1}; the last inequality holds by Assumption 1, Assumption 3 and the above lemma \ref{lem:A4}, i.e., it follows the following fact: According to $\rho = \frac{2\sqrt{237}\kappa_GL}{\sigma^A_{\min}\alpha}$, we have
$\sigma^A_{\min}\rho = \frac{2\sqrt{237}\kappa_GL}{\alpha}$. Let $\beta=\frac{4\alpha}{\sqrt{237}\kappa_G}$, then we have $\frac{4}{\sigma^A_{\min}\rho} = \frac{2\alpha}{\sqrt{237}\kappa_G L}=\frac{\beta}{2L}$. Due to $\kappa_G\geq 1$ and $\alpha \in (0,1]$, we have $0<\beta=\frac{4\alpha}{\sqrt{237}\kappa_G} \leq 1$. According to lemma \ref{lem:A4}, then we have
\begin{align}
    f(x_{k+1}) - \frac{4}{\sigma^A_{\min}\rho}\|\nabla f(x_{k+1})\|^2 \geq f^* .
\end{align}
Then by using the definition of $\Omega_k$, we have, for all $k=0,1,2,\cdots$
\begin{align}
\Omega_{k+1}\geq f^* + \sum_{j=1}^m \psi_j^* - \frac{24qL^2d\mu^2}{\sigma^A_{\min}\rho b}
- \frac{4 L^2d\mu^2}{\sigma^A_{\min}\rho}.
\end{align}
It follows that the function $\Omega_k$ is bounded from below. Let $\Omega^*$ denotes a lower bound of function $\Omega_k$.

Telescoping inequality \eqref{eq:C9} over $k$ from $0$ to $K$, we have
\begin{align}
\Omega_{K} - \Omega_{0}  &= \Omega_{q} - \Omega_{0} + \Omega_{2q} - \Omega_{q} + \cdots + \Omega_{K} - \Omega_{c_Kq}\nonumber \\
& \leq - \sum_{i=0}^{q-1} (\chi\|x_{i+1} - x_{i}\|^2 + \sigma_{\min}^H\sum_{j=1}^m \|y_j^i-y_j^{i+1}\|^2) - \sum_{i=q}^{2q-1} ( \chi\|x_{i+1} - x_{i}\|^2 + \sigma_{\min}^H\sum_{j=1}^m \|y_j^i-y_j^{i+1}\|^2) \nonumber \\
& \quad - \cdots - \sum_{i=c_Kq}^{K-1} ( \chi \|x_{i+1} - x_{i}\|^2  + \sigma_{\min}^H\sum_{j=1}^m \|y_j^i-y_j^{i+1}\|^2) + \big( \frac{6qL}{b} + L + \frac{120qL^2}{b\sigma^A_{\min}\rho}
+ \frac{20 L^2}{\sigma^A_{\min}\rho} \big) K d\mu^2 \nonumber \\
& = - \sum_{i=0}^{K-1} (\chi \|x_{i+1} - x_{i}\|^2+ \sigma_{\min}^H\sum_{j=1}^m \|y_j^i-y_j^{i+1}\|^2) + \big( \frac{6qL}{b} + L + \frac{120qL^2}{b\sigma^A_{\min}\rho}
+ \frac{20 L^2}{\sigma^A_{\min}\rho} \big) K d\mu^2.
\end{align}
Since $b=q$, we obtain
\begin{align} \label{eq:C10}
\frac{1}{K}\sum_{k=0}^{K-1} \big(\mathbb{E}\|x_{k+1} - x_{k}\|^2 + \sum_{j=1}^m \|y_j^k-y_j^{k+1}\|^2\big) \leq \frac{\Omega_{0} - \Omega^*}{K\gamma}
+ \frac{7\vartheta_1 L^2d\mu^2}{\gamma},
\end{align}
where $\vartheta_1 = \frac{1}{L} + \frac{20}{\sigma^A_{\min}\rho}$, $\gamma = \min(\chi,\sigma_{\min}^H)$ with $\chi\geq \frac{\sqrt{237}\kappa_GL}{2\alpha}$.

\end{proof}

Next, based on the above lemmas, we study the convergence properties of ZO-SPIDER-ADMM (CooGE) algorithm.
First, we define a useful variable $\theta_k = \mathbb{E}\big[ \|x_{k+1}-x_{k}\|^2+\|x_{k}-x_{k-1}\|^2+\frac{1}{q}\sum_{i=c_kq}^k \|x_{i+1}-x_i\|^2 + \sum_{j=1}^m \|y_j^k-y_j^{k+1}\|^2 \big] $.

\begin{theorem} \label{th:A1}
(Restatement of Theorem \ref{th:1})
Suppose the sequence $\{x_k,y_{[m]}^k,\lambda_k)_{k=1}^K$ be generated from Algorithm \ref{alg:1}. Under the above Assumptions,
the function $\Omega_k$ has a lower bound.
Further, let $b=q$, $\eta = \frac{\alpha\sigma_{\min}(G)}{4L} \ (0<\alpha \leq 1)$ and $\rho = \frac{2\sqrt{237} \kappa_G L}{\sigma^A_{\min}\alpha}$,
we have
\begin{align}
\mathbb{E}\big[ \mbox{dist}(0,\partial L(x_\zeta,y_{[m]}^\zeta,\lambda_\zeta))^2\big] =\frac{1}{K}\sum_{k=1}^K \mathbb{E}\big[ \mbox{dist}(0,\partial L(x_k,y_{[m]}^k,\lambda_k))^2\big] \leq O(\frac{1}{K}) + O(d\mu^2),
\end{align}
where $\{x_\zeta,y_{[m]}^{\zeta},\lambda_\zeta\}$ chosen uniformly randomly from $\{x_{k},y_{[m]}^{k},\lambda_k\}_{k=0}^{K-1}$.
It implies that given
$K = O(\frac{1}{\epsilon}), \ \mu = O(\sqrt{\frac{\epsilon}{d}})$,
then $(x_{k^*},y^{k^*}_{[m]},\lambda_{k^*})$ is an $\epsilon$-approximate stationary point of the problem \eqref{eq:1}, where $k^* = \mathop{\arg\min}_{k}\theta_{k}$.
\end{theorem}

\begin{proof}
By the optimal condition of the step 9 in Algorithm \ref{alg:1}, we
have, for all $j\in [m]$
\begin{align} \label{eq:C11}
\mathbb{E}\big[\mbox{dist}(0,\partial_{y_j} L(x,y_{[m]},\lambda))^2\big]_{k+1} &= \mathbb{E}\big[\mbox{dist} (0, \partial \psi_j(y_j^{k+1})-B_j^T\lambda_{k+1})^2\big] \nonumber \\
& = \|B_j^T\lambda_k -\rho B_j^T(Ax_k + \sum_{i=1}^jB_iy_i^{k+1} + \sum_{i=j+1}^mB_iy_i^{k} -c) - H_j(y_j^{k+1}-y_j^k) -B_j^T\lambda_{k+1}\|^2 \nonumber \\
& = \|\rho B_j^TA(x_{k+1}-x_{k}) + \rho B_j^T \sum_{i=j+1}^m B_i (y_i^{k+1}-y_i^{k})- H_j(y_j^{k+1}-y_j^k) \|^2 \nonumber \\
& \leq m\rho^2\sigma^{B_j}_{\max}\sigma^A_{\max}\|x_{k+1}\!-\!x_k\|^2 \!+\! m\rho^2\sigma^{B_j}_{\max}\!\sum_{i=j+1}^m \!\sigma^{B_i}_{\max}\|y_i^{k+1}-y_i^{k}\|^2
\!+\! m\sigma^2_{\max}(H_j)\|y_j^{k+1}\!-\!y_j^k\|^2\nonumber \\
& \leq m\big(\rho^2\sigma^B_{\max}\sigma^A_{\max} + \rho^2(\sigma^B_{\max})^2 + \sigma^2_{\max}(H)\big) \theta_{k},
\end{align}
where the first inequality follows by the inequality $\|\sum_{i=1}^r \alpha_i\|^2 \leq r\sum_{i=1}^r \|\alpha_i\|^2$.

By the step 10 of Algorithm \ref{alg:1}, we have
\begin{align} \label{eq:C12}
\mathbb{E}\big[\mbox{dist}(0,\nabla_x L(x,y_{[m]},\lambda))^2\big]_{k+1} & = \mathbb{E}\|A^T\lambda_{k+1}-\nabla f(x_{k+1})\|^2  \nonumber \\
& = \mathbb{E}\|v_k - \nabla f(x_{k+1}) - \frac{G}{\eta} (x_k - x_{k+1})\|^2 \nonumber \\
& = \mathbb{E}\|v_k - \nabla f(x_{k}) + \nabla f(x_{k})- \nabla f(x_{k+1}) - \frac{G}{\eta}(x_k-x_{k+1})\|^2  \nonumber \\
& \leq \frac{18L^2}{b} \sum_{i=c_kq}^{k-1} \mathbb{E}\|x_{i+1}-x_i\|^2 \!+\! \frac{36qL^2d\mu^2}{b} \!+\! 6L^2d\mu^2 \!+\! 3(L^2+ \frac{\sigma^2_{\max}(G)}{\eta^2})\|x_k-x_{k+1}\|^2  \nonumber \\
& \leq 18(L^2+ \frac{\sigma^2_{\max}(G)}{\eta^2})\theta_{k} + 42L^2d\mu^2,
\end{align}
where where the first inequality follows the above Lemma \ref{lem:B1}, and the second inequality holds by $b=q$.

By the step 11 of Algorithm \ref{alg:1}, we have
\begin{align} \label{eq:C13}
\mathbb{E}\big[\mbox{dist}(0,\nabla_\lambda L(x,y_{[m]},\lambda))^2\big]_{k+1} & = \mathbb{E}\|Ax_{k+1}+ \sum_{j=1}^m B_jy_j^{k+1}-c\|^2 \nonumber \\
&= \frac{1}{\rho^2} \mathbb{E} \|\lambda_{k+1}-\lambda_k\|^2  \nonumber \\
& \leq \frac{60L^2}{b\sigma^A_{\min}\rho^2} \sum_{i=c_kq}^{k-1} \mathbb{E}\|x_{i+1}-x_i\|^2
+ (\frac{5L^2}{\sigma^A_{\min}\rho^2}+\frac{5\sigma^2_{\max}(G)}{\sigma^A_{\min}\eta^2\rho^2})\|x_k-x_{k-1}\|^2 \nonumber \\
& \quad + \frac{5\sigma^2_{\max}(G)}{\sigma^A_{\min}\eta^2\rho^2}\|x_{k+1}-x_{k}\|^2  + \frac{120qL^2d\mu^2}{b\sigma^A_{\min}\rho^2} + \frac{20 L^2d\mu^2}{\sigma^A_{\min}\rho^2} \nonumber \\
& \leq \big( \frac{60 L^2 }{\sigma^A_{\min} \rho^2} + \frac{5\sigma^2_{\max}(G) }{\sigma^A_{\min}\eta^2\rho^2} \big) \theta_{k}
+ \frac{140L^2d\mu^2}{\rho^2\sigma^A_{\min}},
\end{align}
where the first inequality follows the above Lemma \ref{lem:B2}, and the last inequality holds by $b=q$.

Let $\beta_{\max}=\max\{\beta_1,\beta_2,\beta_3\}$ with
\begin{align}
  \beta_1 \!=\! m\big(\rho^2\sigma^B_{\max}\sigma^A_{\max} \!+\! \rho^2(\sigma^B_{\max})^2 \!+\! \sigma^2_{\max}(H)\big), \ \beta_2 \!=\! 18(L^2 \!+\! \frac{\sigma^2_{\max}(G)}{\eta^2}), \
  \beta_3 \!=\! \frac{35L^2 }{\sigma^A_{\min} \rho^2} \!+\! \frac{5\sigma^2_{\max}(G) }{\sigma^A_{\min}\eta^2\rho^2}. \nonumber
\end{align}
Combining the above inequalities \eqref{eq:C11}, \eqref{eq:C12} and \eqref{eq:C13},
we have
\begin{align}
\frac{1}{K}\sum_{k=1}^K \mathbb{E}\big[ \mbox{dist}(0,\partial L(x_k,y_{[m]}^k,\lambda_k))^2\big] & \leq \frac{\beta_{\max}}{K}\sum_{k=1}^{K-1} \theta_k  + 14\vartheta_2L^2d\mu^2 \nonumber \\
& \leq  \frac{3\beta_{\max}(\Omega_{0} - \Omega^*)}{K \gamma} + \frac{21\vartheta_1L^2 d\mu^2}{\gamma} + 14\vartheta_2L^2d\mu^2,
\end{align}
where the second inequality holds by the above Lemma \ref{lem:B3} and $\sum_{k=0}^{K-1}\sum_{i=c_kq}^k \|x_{i+1}-x_i\|^2 \leq q\sum_{k=0}^{K-1} \|x_{k+1}-x_k\|^2$, $\vartheta_1 = \frac{1}{L} + \frac{20}{\sigma^A_{\min}\rho}$, $\vartheta_2 =\max\big\{3, \frac{10}{\rho^2\sigma^A_{\min}}\big\}$,
and $\gamma \geq \frac{\sqrt{237}\kappa_GL}{2\alpha}$.

Let $\eta = \frac{\alpha\sigma_{\min}(G)}{4L} \ (0<\alpha \leq 1)$ and $\rho = \frac{2\sqrt{237} \kappa_G L}{\sigma^A_{\min}\alpha}$.
Since $m$ denotes the number of nonsmooth regularization functions, it is relatively small.
It is easily verified that $\beta_{\max} = O(1)$ and $\gamma=O(1)$, which are independent on $n$ and $K$. Thus, we obtain
\begin{align}
\mathbb{E}\big[ \mbox{dist}(0,\partial L(x_\zeta,y_{[m]}^\zeta,\lambda_\zeta))^2\big] =\frac{1}{K}\sum_{k=1}^K \mathbb{E}\big[ \mbox{dist}(0,\partial L(x_k,y_{[m]}^k,\lambda_k))^2\big] \leq  O(\frac{1}{K}) + O(d\mu^2),
\end{align}
where $\{x_\zeta,y_{[m]}^{\zeta},\lambda_\zeta\}$ chosen uniformly randomly from $\{x_{k},y_{[m]}^{k},\lambda_k\}_{k=0}^{K-1}$.

\end{proof}

\subsection{ Convergence Analysis of ZO-SPIDER-ADMM (CooGE+UniGE) Algorithm }
\label{Appendix:A2}
In the subsection, we analyze the convergence of the ZO-SPIDER-ADMM (CooGE+UniGE) algorithm.
We begin with giving an upper bound of variance of stochastic zeroth-order gradient $v_k$.

\begin{lemma} \label{lem:M1}
 Suppose the zeroth-order stochastic gradient $v_t$ be generated from the Algorithm \ref{alg:2}, we have
 \begin{align}
  \mathbb{E} \|\nabla f(x_k)-v_k\|^2 \leq \frac{6dL^2}{b}\sum_{i=c_kq}^{k-1}\mathbb{E}\|x_{i+1}-x_{i}\|^2 + \frac{3\nu^2 L^2 d^2}{2} + \frac{3qL^2d^2\nu^2}{b} + 4L^2d\mu^2.
 \end{align}
\end{lemma}

\begin{proof}
 We first define $f_{\nu}(x)=\mathbb{E}_{u\sim U_B}[f(x+\nu u)]$ be a smooth approximation of $f(x)$,
 where $U_B$ is the uniform distribution over the $d$-dimensional unit Euclidean ball $B$. By the Lemma \ref{lem:A3},
 we have $\mathbb{E}_{(u,\xi)}[\hat{\nabla}_{\texttt{uni}} f_{\xi}(x)]=f_{\nu}(x)$.

 Next, we give an upper bound of $\mathbb{E}\|v_k - \nabla f_{\nu}(x_{k})\|^2$. By the definition of $v_k$,
 we have
 \begin{align} \label{eq:M1}
  &\mathbb{E}\| \nabla f_{\nu}(x_{k})-v_k\|^2= \mathbb{E}\|\underbrace{\nabla f_{\nu}(x_{k}) -\nabla f_{\nu}(x_{k-1}) -\frac{1}{b} \sum_{j \in \mathcal{S}} [\hat{\nabla}_{\texttt{uni}} f_{j}(x_k) - \hat{\nabla}_{\texttt{uni}} f_{j}(x_{k-1})] }_{=T_1} + \underbrace{\nabla f_{\nu}(x_{k-1})- v_{k-1}}_{=T_2}\|^2 \nonumber \\
  & = \mathbb{E}\|\nabla f_{\nu}(x_{k}) -\nabla f_{\nu}(x_{k-1}) -\frac{1}{b} \sum_{j \in \mathcal{S}} [\hat{\nabla}_{\texttt{uni}} f_{j}(x_k) - \hat{\nabla}_{\texttt{uni}} f_{j}(x_{k-1})] \|^2
  + \mathbb{E}\| \nabla f_{\nu}(x_{k-1})- v_{k-1}\|^2 \nonumber \\
  & = \frac{1}{b^2} \sum_{j \in \mathcal{S}}\mathbb{E}\| \nabla f_{\nu}(x_{k}) -\nabla f_{\nu}(x_{k-1}) - \hat{\nabla}_{\texttt{uni}} f_{j}(x_k) + \hat{\nabla}_{\texttt{uni}} f_{j}(x_{k-1}) \|^2
  + \mathbb{E}\| \nabla f_{\nu}(x_{k-1}) - v_{k-1}\|^2 \nonumber \\
  & \leq\frac{1}{b^2} \sum_{j \in \mathcal{S}}\mathbb{E}\| \hat{\nabla}_{\texttt{uni}} f_{j}(x_k) - \hat{\nabla}_{\texttt{uni}} f_{j}(x_{k-1}) \|^2
  + \mathbb{E}\| \nabla f_{\nu}(x_{k-1}) - v_{k-1}\|^2 \nonumber \\
  & \leq \frac{1}{b} \big( 3dL^2\|x_{k}-x_{k-1}\|^2 + \frac{3L^2d^2\nu^2}{2} \big) +  \mathbb{E}\| \nabla f_{\nu}(x_{k-1})- v_{k-1}\|^2,
 \end{align}
 where the second equality follows by $\mathbb{E}[T_1]=0$ and $T_2$ is independent to $\mathcal{S}$, and
 the third equality holds by the above Lemma \ref{lem:A0}, and the second inequality holds by the Lemma \ref{lem:A3}.

 When $k=c_kq$, we have $v_k = \frac{1}{n}\sum_{j=1}^n \hat{\nabla}_{\texttt{coo}} f_{j}(x_{k})=\hat{\nabla}_{\texttt{coo}} f(x_{k})$.
 We have
 \begin{align} \label{eq:M2}
  \mathbb{E}\| \nabla f_{\nu}(x_{c_kq})- v_{c_kq}\|^2 & = \mathbb{E}\| \nabla f_{\nu}(x_{c_kq})- \nabla f(x_{c_kq}) + \nabla f(x_{c_kq}) -\hat{\nabla}_{\texttt{coo}} f(x_{c_kq})\|^2 \nonumber \\
  & \leq 2\mathbb{E}\| \nabla f_{\nu}(x_{c_kq})- \nabla f(x_{c_kq})\|^2 + 2\mathbb{E}\|\nabla f(x_{c_kq}) - \hat{\nabla}_{\texttt{coo}} f(x_{c_kq})\|^2 \nonumber \\
  & \leq \frac{\nu^2L^2d^2}{2} + 2L^2d\mu^2,
 \end{align}
 where the final inequalities holds by Lemmas \ref{lem:A1} and \ref{lem:A3}.

 By recursion to \eqref{eq:M1}, we have
 \begin{align}
 \mathbb{E}\|\nabla f_{\nu}(x_{k})-v_k\|^2 & \leq (k-c_kq)\frac{3L^2d^2\nu^2}{2b} + \frac{3dL^2}{b}\sum_{i=c_kq}^{k-1}\mathbb{E}\|x_{i+1}-x_{i}\|^2
 + \mathbb{E}\| \nabla f_{\nu}(x_{c_kq})- v_{c_kq}\|^2 \nonumber \\
 & \leq \frac{3qL^2d^2\nu^2}{2b} + \frac{3dL^2}{b}\sum_{i=c_kq}^{k-1}\mathbb{E}\|x_{i+1}-x_{i}\|^2 + \frac{\nu^2L^2d^2}{2} + 2L^2d\mu^2,
 \end{align}
 where the last inequality holds by the above inequality \eqref{eq:M2}
 and $k-c_kq \leq q$.

 Finally, we have
 \begin{align}
  \mathbb{E} \|\nabla f(x_k)-v_k\|^2 &=\mathbb{E} \|\nabla f(x_k)- \nabla f_{\nu}(x_{k}) + \nabla f_{\nu}(x_{k})- v_k\|^2 \nonumber\\
  & \leq 2\mathbb{E} \|\nabla f(x_k)- \nabla f_{\nu}(x_{k})\|^2 + 2\mathbb{E} \| \nabla f_{\nu}(x_{k})- v_k\|^2 \nonumber\\
  & \leq \frac{6dL^2}{b}\sum_{i=c_kq}^{k-1}\mathbb{E}\|x_{i+1}-x_{i}\|^2 + \frac{3\nu^2 L^2 d^2}{2} + \frac{3qL^2d^2\nu^2}{b} + 4L^2d\mu^2,
 \end{align}
 where the last inequality holds by the above Lemma \ref{lem:A3}.
\end{proof}

\begin{lemma} \label{lem:M2}
Suppose the sequence $\{x_k,y_{[m]}^k,\lambda_k\}_{k=1}^K$ be generated from Algorithm \ref{alg:2}, it holds that
\begin{align}
\mathbb{E}\|\lambda_{k+1} - \lambda_k\|^2 \leq & \frac{60dL^2}{b\sigma^A_{\min}} \sum_{i=c_kq}^{k-1} \mathbb{E}\|x_{i+1}-x_i\|^2
+ (\frac{5L^2}{\sigma^A_{\min}}+\frac{5\sigma^2_{\max}(G)}{\sigma^A_{\min}\eta^2})\|x_k-x_{k-1}\|^2 \nonumber \\
& + \frac{5\sigma^2_{\max}(G)}{\sigma^A_{\min}\eta^2}\|x_{k+1}-x_{k}\|^2  + \frac{15\nu^2 L^2 d^2 }{\sigma^A_{\min}}
+ \frac{30qL^2d^2\nu^2}{b\sigma^A_{\min}} + \frac{40L^2d\mu^2}{\sigma^A_{\min}}.
\end{align}
\end{lemma}
\begin{proof}
 This proof can follow the proof of the Lemma \ref{lem:B2}.
\end{proof}

\begin{lemma} \label{lem:M3}
Suppose the sequence $\{x_k,y_{[m]}^k,\lambda_k\}_{k=1}^K$ be generated from Algorithm \ref{alg:2},
and define a \emph{Lyapunov} function $\Phi_k$ as follows:
\begin{align}
\Phi_k = \mathbb{E} \big[ \mathcal{L}_{\rho} (x_k,y_{[m]}^k,\lambda_k) + (\frac{5L^2}{\sigma^A_{\min}\rho}+\frac{5\sigma^2_{\max}(G)}{\sigma^A_{\min}\eta^2\rho})\|x_k-x_{k-1}\|^2 + \frac{12dL^2}{b\sigma^A_{\min}\rho}\sum_{i=c_kq}^{k-1}\mathbb{E}\|x_{i+1}-x_i\|^2 \big]. \nonumber
\end{align}
Let $b=qd$, $\eta=\frac{\alpha\sigma_{\min}(G)}{4L} \ (0<\alpha \leq 1)$ and $\rho = \frac{2\sqrt{237}\kappa_GL}{\sigma^A_{\min}\alpha}$,
we have
\begin{align}
\frac{1}{K}\sum_{k=0}^{K-1} \big(\mathbb{E}\|x_{k+1} - x_{k}\|^2 + \sum_{j=1}^m \|y_j^k-y_j^{k+1}\|^2\big) \leq \frac{\Phi_{0} - \Phi^*}{K\gamma}
+ \big(\frac{3L}{4} + \frac{3qL}{2b}\big)\frac{\vartheta_1d^2\nu^2}{\gamma} + \frac{2\vartheta_1Ld\mu^2}{\gamma},
\end{align}
where $\vartheta_1 = 1 + \frac{20L}{\sigma^A_{\min}\rho}$, $\gamma = \min(\chi,\sigma_{\min}^H)$ with $\chi\geq \frac{\sqrt{237}\kappa_GL}{2\alpha}$,
and $\Phi^*$ is a lower bound of the function $\Phi_k$.
\end{lemma}

\begin{proof}
The proof can follow the proof of Lemma \ref{lem:B3}.
Similarly, we can obtain
\begin{align} \label{eq:M3}
\mathcal{L}_{\rho} (x_k,y^{k+1}_{[m]},\lambda_k) \leq \mathcal{L}_{\rho} (x_k,y^k_{[m]},\lambda_k)
- \sigma_{\min}^H\sum_{j=1}^m \|y_j^k-y_j^{k+1}\|^2,
\end{align}
where $\sigma_{\min}^H=\min_{j\in[m]}\sigma_{\min}(H_j)$.

Similarly, we have
\begin{align}
0 & \leq \mathcal{L}_{\rho} (x_k,y_{[m]}^{k+1},\lambda_k) -  \mathcal{L}_{\rho} (x_{k+1},y_{[m]}^{k+1},\lambda_k)
- (\frac{\sigma_{\min}(G)}{\eta}+\frac{\rho \sigma^A_{\min}}{2}-L ) \|x_{k+1} - x_{k}\|^2 + \frac{1}{2L}\|v_k - \nabla f(x_k)\|^2 \nonumber \\
& \leq \mathcal{L}_{\rho} (x_k,y_{[m]}^{k+1},\lambda_k) -  \mathcal{L}_{\rho} (x_{k+1},y_{[m]}^{k+1},\lambda_k)
- (\frac{\sigma_{\min}(G)}{\eta}+\frac{\rho \sigma^A_{\min}}{2}-L ) \|x_{k+1} - x_{k}\|^2 \nonumber \\
& \quad + \frac{3dL}{b}\sum_{i=c_kq}^{k-1}\mathbb{E}\|x_{i+1}-x_{i}\|^2 + \frac{3\nu^2 L d^2}{4} + \frac{3qLd^2\nu^2}{2b} + 2Ld\mu^2, \nonumber
\end{align}
where the final inequality holds by Lemma \ref{lem:M1}.
It follows that
\begin{align} \label{eq:M4}
\mathcal{L}_{\rho} (x_{k+1},y_{[m]}^{k+1},\lambda_k) & \leq \mathcal{L}_{\rho} (x_k,y_{[m]}^{k+1},\lambda_k)
- (\frac{\sigma_{\min}(G)}{\eta}+\frac{\rho \sigma^A_{\min}}{2}-L )\mathbb{E}\|x_{k+1} - x_{k}\|^2 \nonumber \\
& \quad + \frac{3dL}{b}\sum_{i=c_kq}^{k-1}\mathbb{E}\|x_{i+1}-x_{i}\|^2 + \frac{3\nu^2 L d^2}{4} + \frac{3qLd^2\nu^2}{2b} + 2Ld\mu^2.
\end{align}

By using the step 11 in Algorithm \ref{alg:2}, we have
\begin{align} \label{eq:M5}
\mathcal{L}_{\rho} (x_{k+1},y_{[m]}^{k+1},\lambda_{k+1}) -\mathcal{L}_{\rho} (x_{k+1},y_{[m]}^{k+1},\lambda_k)
& = \frac{1}{\rho}\mathbb{E}\|\lambda_{k+1}-\lambda_k\|^2  \\
& \leq \frac{60dL^2}{\sigma^A_{\min}b\rho}\sum_{i=c_kq}^{k-1}\mathbb{E}\|x_{i+1}-x_i\|^2 + \big(\frac{5\sigma^2_{\max}(G)}{\sigma^A_{\min}\eta^2\rho}+\frac{5L^2}{\sigma^A_{\min}\rho}\big)\mathbb{E}\|x_k-x_{k-1}\|^2\nonumber \\
& \quad + \frac{5\sigma^2_{\max}(G)}{\sigma^A_{\min}\eta^2\rho}\mathbb{E}\|x_{k+1}-x_{k}\|^2 + \frac{15\nu^2 L^2 d^2 }{\sigma^A_{\min}\rho}
+ \frac{30qL^2d^2\nu^2}{b\sigma^A_{\min}\rho} + \frac{40L^2d\mu^2}{\sigma^A_{\min}\rho}, \nonumber
\end{align}
where the above inequality holds by Lemma \ref{lem:M2}.

Combining \eqref{eq:M3}, \eqref{eq:M4} and \eqref{eq:M5}, we have
\begin{align} \label{eq:M6}
\mathcal{L}_{\rho} (x_{k+1},y_{[m]}^{k+1},\lambda_{k+1}) & \leq \mathcal{L}_{\rho} (x_k,y_{[m]}^k,\lambda_k)
 \!-\! \sigma_{\min}^H\sum_{j=1}^m \|y_j^k-y_j^{k+1}\|^2 \!-\! \big(\frac{\sigma_{\min}(G)}{\eta}+\frac{\rho \sigma^A_{\min}}{2}-L - \frac{5\sigma^2_{\max}(G)}{\sigma^A_{\min}\eta^2\rho}\big)
 \mathbb{E}\|x_{k+1} - x_{k}\|^2 \nonumber \\
& \quad + \big(\frac{5L^2}{\sigma^A_{\min}\rho}+\frac{5\sigma^2_{\max}(G)}{\sigma^A_{\min}\eta^2\rho}\big)\mathbb{E}\|x_k-x_{k-1}\|^2 + \big(\frac{60dL^2}{\sigma^A_{\min}b\rho} + \frac{3dL}{b}  \big)\sum_{i=c_kq}^{k-1}\mathbb{E}\|x_{i+1}-x_i\|^2 \nonumber \\
& \quad + \big( \frac{3L}{4} + \frac{3qL}{2b} + \frac{15 L^2}{\sigma^A_{\min}\rho}+ \frac{30qL^2}{b\sigma^A_{\min}\rho}\big) d^2\nu^2
+ (2L + \frac{40L^2}{\sigma^A_{\min}\rho})d\mu^2.
\end{align}

We define a useful \emph{Lyapunov} function $\Phi_k$ as follows:
\begin{align}
\Phi_k = \mathbb{E} \big[ \mathcal{L}_{\rho} (x_k,y_{[m]}^k,\lambda_k) + (\frac{5L^2}{\sigma^A_{\min}\rho}+\frac{5\sigma^2_{\max}(G)}{\sigma^A_{\min}\eta^2\rho})\|x_k-x_{k-1}\|^2 + \frac{12dL^2}{b\sigma^A_{\min}\rho}\sum_{i=c_kq}^{k-1}\mathbb{E}\|x_{i+1}-x_i\|^2 \big].
\end{align}
Following the above lemma \ref{lem:B3}, it is easily verified that the function $\Phi_k$ is bounded from below.
Let $\Phi^*$ denotes a lower bound of function $\Phi_k$.
It follows that
\begin{align} \label{eq:M7}
\Phi_{k+1} & = \mathbb{E} \big[ \mathcal{L}_{\rho} (x_{k+1},y_{[m]}^{k+1},\lambda_{k+1}) + (\frac{5L^2}{\sigma^A_{\min}\rho} + \frac{5\sigma^2_{\max}(G)}{\sigma^A_{\min}\eta^2\rho})\|x_{k+1}-x_{k}\|^2
+ \frac{12dL^2}{b\sigma^A_{\min}\rho}\sum_{i=c_kq}^{k}\mathbb{E}\|x_{i+1}-x_i\|^2 \big]\nonumber \\
& \leq \mathbb{E} \big[ \mathcal{L}_{\rho} (x_k,y_{[m]}^k,\lambda_k) + (\frac{5L^2}{\sigma^A_{\min}\rho}+\frac{5\sigma^2_{\max}(G)}{\sigma^A_{\min}\eta^2\rho})\|x_{k}-x_{k-1}\|^2
+ \frac{12dL^2}{b\sigma^A_{\min}\rho}\sum_{i=c_kq}^{k-1}\mathbb{E}\|x_{i+1}-x_i\|^2 \big]   \nonumber \\
& \quad - \sigma_{\min}^H\sum_{j=1}^m \|y_j^k-y_j^{k+1}\|^2 - \big(\frac{\sigma_{\min}(G)}{\eta}+\frac{\rho \sigma^A_{\min}}{2}-L - \frac{10\sigma^2_{\max}(G)}{\sigma^A_{\min}\eta^2\rho}-\frac{5L^2}{\sigma^A_{\min}\rho}
 - \frac{12dL^2}{b\sigma^A_{\min}\rho} \big) \mathbb{E}\|x_{k+1} - x_{k}\|^2 \nonumber \\
& \quad + \big(\frac{60dL^2}{\sigma^A_{\min}b\rho} + \frac{3dL}{b}  \big)\sum_{i=c_kq}^{k-1}\mathbb{E}\|x_{i+1}-x_i\|^2  + \big( \frac{3L}{4} + \frac{3qL}{2b} + \frac{15 L^2}{\sigma^A_{\min}\rho}+ \frac{30qL^2}{b\sigma^A_{\min}\rho}\big) d^2\nu^2
+ (2L + \frac{40L^2}{\sigma^A_{\min}\rho})d\mu^2 \nonumber \\
& = \Phi_k - \sigma_{\min}^H\sum_{j=1}^m \|y_j^k-y_j^{k+1}\|^2 - \big(\frac{\sigma_{\min}(G)}{\eta}+\frac{\rho \sigma^A_{\min}}{2}-L - \frac{10\sigma^2_{\max}(G)}{\sigma^A_{\min}\eta^2\rho}-\frac{5L^2}{\sigma^A_{\min}\rho}
 - \frac{12dL^2}{b\sigma^A_{\min}\rho} \big) \mathbb{E}\|x_{k+1} - x_{k}\|^2 \nonumber \\
& \quad + \big(\frac{60dL^2}{\sigma^A_{\min}b\rho} + \frac{3dL}{b} \big)\sum_{i=c_kq}^{k-1}\mathbb{E}\|x_{i+1}-x_i\|^2 + \big( \frac{3L}{4} + \frac{3qL}{2b} + \frac{15 L^2}{\sigma^A_{\min}\rho}+ \frac{30qL^2}{b\sigma^A_{\min}\rho}\big) d^2\nu^2
+ (2L + \frac{40L^2}{\sigma^A_{\min}\rho})d\mu^2,
\end{align}
where the first inequality holds by the inequality \eqref{eq:M6} and the equality
$$\sum_{i=c_kq}^{k}\mathbb{E}\|x_{i+1}-x_i\|^2 = \sum_{i=c_kq}^{k-1}\mathbb{E}\|x_{i+1}-x_i\|^2 + \mathbb{E}\|x_{k+1}-x_k\|^2.$$

Since $c_kq \leq k \leq (c_k+1)q -1$, and let $c_kq \leq t \leq (c_k+1)q-1$, then telescoping inequality \eqref{eq:M7} over $k$ from $c_kq$ to $k$, we have
\begin{align} \label{eq:M8}
\Phi_{k+1} & \leq \Phi_{c_kq}
 - \big(\frac{\sigma_{\min}(G)}{\eta}+\frac{\rho \sigma^A_{\min}}{2}-L - \frac{10\sigma^2_{\max}(G)}{\sigma^A_{\min}\eta^2\rho}-\frac{5L^2}{\sigma^A_{\min}\rho}
 - \frac{12dL^2}{b\sigma^A_{\min}\rho} \big) \sum_{t=c_kq}^k \mathbb{E}\|x_{t+1} - x_{t}\|^2 \nonumber \\
& \quad - \sigma_{\min}^H\sum_{t=c_kq}^{k}\sum_{j=1}^m \|y_j^t\!-\!y_j^{t+1}\|^2 + \big(\frac{60dL^2}{\sigma^A_{\min}b\rho} + \frac{3dL}{b} \big)
\sum_{t=c_kq}^{k}\sum_{i=c_kq}^{k-1}\mathbb{E}\|x_{i+1}-x_i\|^2 \nonumber \\
& \quad + \sum_{t=c_kq}^{k}\big( \frac{3L}{4} + \frac{3qL}{2b} + \frac{15 L^2}{\sigma^A_{\min}\rho}+ \frac{30qL^2}{b\sigma^A_{\min}\rho}\big) d^2\nu^2 + \sum_{t=c_kq}^{k}\big(2L + \frac{40L^2}{\sigma^A_{\min}\rho}\big)d\mu^2 \nonumber \\
& \leq \Phi_{c_kq} - (\frac{\sigma_{\min}(G)}{\eta}+\frac{\rho \sigma^A_{\min}}{2}-L - \frac{10\sigma^2_{\max}(G)}{\sigma^A_{\min}\eta^2\rho}-\frac{5L^2}{\sigma^A_{\min}\rho}
 -\frac{12dL^2}{\sigma^A_{\min}\rho b}) \sum_{i = c_kq}^k \mathbb{E}\|x_{i+1} - x_{i}\|^2 \nonumber \\
& \quad - \sigma_{\min}^H\sum_{i=c_kq}^{k-1}\sum_{j=1}^m \|y_j^i-y_j^{i+1}\|^2 + \big(\frac{60qdL^2}{\sigma^A_{\min}b\rho} + \frac{3qdL}{b} \big)\sum_{i=c_kq}^{k}\mathbb{E}\|x_{i+1}-x_i\|^2 \nonumber \\
& \quad + \big( \frac{3L}{4} + \frac{3qL}{2b} + \frac{15 L^2}{\sigma^A_{\min}\rho}+ \frac{30qL^2}{b\sigma^A_{\min}\rho}\big) qd^2\nu^2
+ (2L + \frac{40L^2}{\sigma^A_{\min}\rho})qd\mu^2 \nonumber \\
& = \Phi_{c_kq} - \underbrace{ \big(\frac{\sigma_{\min}(G)}{\eta}+\frac{\rho \sigma^A_{\min}}{2}-L - \frac{10\sigma^2_{\max}(G)}{\sigma^A_{\min}\eta^2\rho}-\frac{5L^2}{\sigma^A_{\min}\rho}
- \frac{12dL^2}{b\sigma^A_{\min}\rho}-\frac{60qdL^2}{\sigma^A_{\min}b\rho}
- \frac{3qdL}{b} \big)}_{\chi}  \sum_{i=c_kq}^k \|x_{i+1} - x_{i}\|^2 \nonumber \\
& \quad - \sigma_{\min}^H\sum_{i=c_kq}^{k-1}\sum_{j=1}^m \|y_j^i-y_j^{i+1}\|^2  + \big( \frac{3L}{4} + \frac{3qL}{2b} + \frac{15 L^2}{\sigma^A_{\min}\rho}+ \frac{30qL^2}{b\sigma^A_{\min}\rho}\big) qd^2\nu^2
+ (2L + \frac{40L^2}{\sigma^A_{\min}\rho})qd\mu^2,
\end{align}
where the second inequality holds by the fact that
\begin{align}
\sum_{j=c_kq}^{k}\sum_{i=c_kq}^{k-1}\mathbb{E}\|x_{i+1}-x_i\|^2 \leq \sum_{j=c_kq}^{k}\sum_{i=c_kq}^{k}\mathbb{E}\|x_{i+1}-x_i\|^2 \leq q\sum_{i=c_kq}^{k}\mathbb{E}\|x_{i+1}-x_i\|^2. \nonumber
\end{align}

Since $b=dq$, we have
\begin{align}
\chi & = \frac{\sigma_{\min}(G)}{\eta}+\frac{\rho \sigma^A_{\min}}{2}-L - \frac{10\sigma^2_{\max}(G)}{\sigma^A_{\min}\eta^2\rho}
-\frac{5L^2}{\sigma^A_{\min}\rho}- \frac{12dL^2}{b\sigma^A_{\min}\rho}-\frac{60qdL^2}{\sigma^A_{\min}b\rho} - \frac{3qdL}{b} \nonumber \\
& \geq \underbrace{\frac{\sigma_{\min}(G)}{\eta} - 4L }_{=T_3} + \underbrace{\frac{\rho\sigma^A_{\min}}{2} - \frac{10\sigma^2_{\max}(G)}{\sigma^A_{\min}\eta^2\rho}
-\frac{77L^2}{\sigma^A_{\min}\rho}}_{=T_4},
\end{align}
where the above inequality holds by $q\geq 1$.
Let $0< \eta \leq \frac{\sigma_{\min}(G)}{4L}$, we have $T_3\geq 0$. Further, let $\eta = \frac{\alpha\sigma_{\min}(G)}{4L} \ (0<\alpha \leq 1)$ and
$\rho = \frac{2\sqrt{237}\kappa_GL}{\sigma^A_{\min}\alpha}$,
we have
\begin{align}
T_4 & = \frac{\rho \sigma^A_{\min}}{2} - \frac{10\sigma^2_{\max}(G)}{\sigma^A_{\min}\eta^2\rho}
- \frac{77dL^2}{\sigma^A_{\min}\rho} =  \frac{\rho \sigma^A_{\min}}{2} - \frac{160L^2\kappa^2_G}{\sigma^A_{\min}\rho\alpha^2}
-\frac{77L^2}{\sigma^A_{\min}\rho} \nonumber \\
& \geq \frac{\rho \sigma^A_{\min}}{2} - \frac{160L^2\kappa^2_G}{\sigma^A_{\min}\rho\alpha^2}-\frac{77L^2\kappa^2_G}{\sigma^A_{\min}\rho\alpha^2}  \nonumber \\
& = \frac{\rho \sigma^A_{\min}}{4} + \underbrace{\frac{\rho \sigma^A_{\min}}{4} - \frac{237L^2\kappa^2_G }{\sigma^A_{\min}\rho\alpha^2}}_{\geq 0} \nonumber \\
& \geq \frac{\sqrt{237}\kappa_GL}{2\alpha},
\end{align}
where the first inequality holds by $\kappa_G \geq 1$,
and the second inequality holds by $\rho = \frac{2\sqrt{237}\kappa_GL}{\sigma^A_{\min}\alpha}$.
Thus, we obtain $\chi \geq \frac{\sqrt{237}\kappa_GL}{2\alpha}$.

Telescoping inequality \eqref{eq:M8} over $k$ from $0$ to $K$, we have
\begin{align}
&\Phi_{K} - \Phi_{0}  = \Phi_{q} - \Phi_{0} + \Phi_{2q} - \Phi_{q} + \cdots + \Phi_{K} - \Phi_{c_Kq}\nonumber \\
& \leq - \sum_{i=0}^{q-1} (\chi\|x_{i+1} - x_{i}\|^2 + \sigma_{\min}^H\sum_{j=1}^m \|y_j^i-y_j^{i+1}\|^2) - \sum_{i=q}^{2q-1} ( \chi\|x_{i+1} - x_{i}\|^2
+ \sigma_{\min}^H\sum_{j=1}^m \|y_j^i-y_j^{i+1}\|^2) - \cdots \nonumber \\
& \quad - \sum_{i=c_Kq}^{K-1} ( \chi \|x_{i+1} - x_{i}\|^2  + \sigma_{\min}^H\sum_{j=1}^m \|y_j^i-y_j^{i+1}\|^2) + K\big( \frac{3L}{4} + \frac{3qL}{2b} + \frac{15 L^2}{\sigma^A_{\min}\rho}+ \frac{30qL^2}{b\sigma^A_{\min}\rho}\big)d^2\nu^2
+ K(2L + \frac{40L^2}{\sigma^A_{\min}\rho})d\mu^2 \nonumber \\
& = - \sum_{i=0}^{K-1} (\chi \|x_{i+1} - x_{i}\|^2+ \sigma_{\min}^H\sum_{j=1}^m \|y_j^i-y_j^{i+1}\|^2)
+ K\big( \frac{3L}{4} + \frac{3qL}{2b} + \frac{15 L^2}{\sigma^A_{\min}\rho}+ \frac{30qL^2}{b\sigma^A_{\min}\rho}\big)d^2\nu^2
+ K(2L + \frac{40L^2}{\sigma^A_{\min}\rho})d\mu^2.
\end{align}
Thus, we have
\begin{align} \label{eq:M9}
\frac{1}{K}\sum_{k=0}^{K-1} \big(\mathbb{E}\|x_{k+1} - x_{k}\|^2 + \sum_{j=1}^m \|y_j^k-y_j^{k+1}\|^2\big) \leq \frac{\Phi_{0} - \Phi^*}{K\gamma}
+ \big(\frac{3L}{4} + \frac{3qL}{2b}\big)\frac{\vartheta_1d^2\nu^2}{\gamma} + \frac{2\vartheta_1Ld\mu^2}{\gamma},
\end{align}
where $\vartheta_1 = 1 + \frac{20L}{\sigma^A_{\min}\rho}$, $\gamma = \min(\chi,\sigma_{\min}^H)$ with $\chi\geq \frac{\sqrt{237}\kappa_GL}{2\alpha}$.

\end{proof}

Next, based on the above lemmas, we study the convergence properties of the ZO-SPIDER-ADMM (CooGE+UniGE) algorithm.
First, we define a useful variable $\theta_k = \mathbb{E}\big[ \|x_{k+1}-x_{k}\|^2 + \|x_{k}-x_{k-1}\|^2 + \frac{1}{q}\sum_{i={c_kq}^k}\|x_{i+1}-x_{i}\|^2 + \sum_{j=1}^m \|y_j^k-y_j^{k+1}\|^2 \big]$.

\begin{theorem} \label{th:A2}
(Restatement of Theorem \ref{th:2})
Suppose the sequence $\{x_k,y_{[m]}^k,\lambda_k)_{k=1}^K$ be generated from the Algorithm \ref{alg:2}. Under the above Assumptions,
the function $\Phi_k$ has a lower bound.
Further, let $b=qd$, $\eta = \frac{\alpha\sigma_{\min}(G)}{4L} \ (0<\alpha \leq 1)$ and $\rho = \frac{2\sqrt{237} \kappa_G L}{\sigma^A_{\min}\alpha}$,
we have
\begin{align}
\mathbb{E}\big[ \mbox{dist}(0,\partial L(x_\zeta,y_{[m]}^\zeta,\lambda_\zeta))^2\big] =\frac{1}{K}\sum_{k=1}^K \mathbb{E}\big[ \mbox{dist}(0,\partial L(x_k,y_{[m]}^k,\lambda_k))^2\big] \leq  O(\frac{1}{K})+ O(d^2\nu^2) + O(d\mu^2),
\end{align}
where $\{x_\zeta,y_{[m]}^{\zeta},\lambda_\zeta\}$ is chosen uniformly randomly from $\{x_{k},y_{[m]}^{k},\lambda_k\}_{k=0}^{K-1}$.
It implies that given
$K = O(\frac{1}{\epsilon}), \ \nu = O(\frac{\sqrt{\epsilon}}{d}), \ \mu = O(\frac{\sqrt{\epsilon}}{\sqrt{d}})$,
then $(x_{k^*},y^{k^*}_{[m]},\lambda_{k^*})$ is an $\epsilon$-stationary point of the problem \eqref{eq:1}, where $k^* = \mathop{\arg\min}_{k}\theta_{k}$.
\end{theorem}
\begin{proof}
By the optimal condition of the step 9 in Algorithm \ref{alg:2}, we
have, for all $j\in [m]$
\begin{align} \label{eq:N1}
\mathbb{E}\big[\mbox{dist}(0,\partial_{y_j} L(x,y_{[m]},\lambda))^2\big]_{k+1} & = \mathbb{E}\big[\mbox{dist} (0, \partial \psi_j(y_j^{k+1})-B_j^T\lambda_{k+1})^2\big] \nonumber \\
& = \|B_j^T\lambda_k -\rho B_j^T(Ax_k + \sum_{i=1}^jB_iy_i^{k+1} + \sum_{i=j+1}^mB_iy_i^{k} -c) - H_j(y_j^{k+1}-y_j^k) -B_j^T\lambda_{k+1}\|^2 \nonumber \\
& = \|\rho B_j^TA(x_{k+1}-x_{k}) + \rho B_j^T \sum_{i=j+1}^m B_i (y_i^{k+1}-y_i^{k})- H_j(y_j^{k+1}-y_j^k) \|^2 \nonumber \\
& \leq m\rho^2\sigma^{B_j}_{\max}\sigma^A_{\max}\|x_{k+1}\!-\!x_k\|^2 + m\rho^2\sigma^{B_j}_{\max}\!\sum_{i=j+1}^m \!\sigma^{B_i}_{\max}\|y_i^{k+1}\!-\!y_i^{k}\|^2
+ m\sigma^2_{\max}(H_j)\|y_j^{k+1}\!-\!y_j^k\|^2\nonumber \\
& \leq m\big(\rho^2\sigma^B_{\max}\sigma^A_{\max} + \rho^2(\sigma^B_{\max})^2 + \sigma^2_{\max}(H)\big) \theta_{k},
\end{align}
where the first inequality follows by the inequality $\|\sum_{i=1}^r \alpha_i\|^2 \leq r\sum_{i=1}^r \|\alpha_i\|^2$.

By the step 10 of Algorithm \ref{alg:2}, we have
\begin{align} \label{eq:N2}
\mathbb{E}\big[\mbox{dist}(0,\nabla_x L(x,y_{[m]},\lambda))^2\big]_{k+1} & = \mathbb{E}\|A^T\lambda_{k+1}-\nabla f(x_{k+1})\|^2  \nonumber \\
& = \mathbb{E}\|v_k - \nabla f(x_{k+1}) - \frac{G}{\eta} (x_k - x_{k+1})\|^2 \nonumber \\
& = \mathbb{E}\|v_k - \nabla f(x_{k}) + \nabla f(x_{k})- \nabla f(x_{k+1}) - \frac{G}{\eta}(x_k-x_{k+1})\|^2  \nonumber \\
& \leq 3\mathbb{E}\|v_k - \nabla f(x_{k})\|^2 +  3\mathbb{E}\|\nabla f(x_{k})- \nabla f(x_{k+1})\|^2 + \frac{3\sigma^2_{\max}(G)}{\eta^2}\|x_k-x_{k+1}\|^2  \nonumber \\
& \leq \frac{18dL^2}{b}\sum_{i=c_kq}^{k-1}\mathbb{E}\|x_{i+1}-x_{i}\|^2 + \frac{9\nu^2 L^2 d^2}{2} + \frac{9qL^2d^2\nu^2}{b} + 12L^2d\mu^2 \nonumber \\
& \quad + 3(L^2+ \frac{\sigma^2_{\max}(G)}{\eta^2})\|x_k-x_{k+1}\|^2 \nonumber \\
& \leq 18(L^2+ \frac{\sigma^2_{\max}(G)}{\eta^2})\theta_{k} + \frac{9\nu^2 L^2 d^2}{2} + \frac{9qL^2d^2\nu^2}{b} + 12L^2d\mu^2,
\end{align}
where the second inequality holds by Lemma \ref{lem:M1}, and the last inequality is due to the definition of $\theta_k$ and $b=qb$.

By the step 11 of Algorithm \ref{alg:2}, we have
\begin{align} \label{eq:N3}
\mathbb{E}\big[\mbox{dist}(0,\nabla_\lambda L(x,y_{[m]},\lambda))^2\big]_{k+1} & = \mathbb{E}\|Ax_{k+1}+ \sum_{j=1}^m B_jy_j^{k+1}-c\|^2 \nonumber \\
&= \frac{1}{\rho^2} \mathbb{E} \|\lambda_{k+1}-\lambda_k\|^2  \nonumber \\
& \leq \frac{60dL^2}{\rho^2\sigma^A_{\min}b}\sum_{i=c_kq}^{k-1}\mathbb{E}\|x_{i+1}-x_i\|^2 + \big(\frac{5\sigma^2_{\max}(G)}{\rho^2\sigma^A_{\min}\eta^2}+\frac{5L^2}{\rho^2\sigma^A_{\min}}\big)\mathbb{E}\|x_k-x_{k-1}\|^2\nonumber \\
& \quad + \frac{5\sigma^2_{\max}(G)}{\rho^2\sigma^A_{\min}\eta^2}\|x_{k+1}-x_{k}\|^2 + \frac{15\nu^2 L^2 d^2 }{\sigma^A_{\min}\rho^2}
+ \frac{30qL^2d^2\nu^2}{b\sigma^A_{\min}\rho^2} + \frac{40L^2d\mu^2}{\sigma^A_{\min}\rho^2} \nonumber \\
& \leq \big( \frac{60L^2}{\sigma^A_{\min} \rho^2} + \frac{5\sigma^2_{\max}(G) }{\sigma^A_{\min}\eta^2\rho^2} \big) \theta_{k}
+ \frac{15\nu^2 L^2 d^2 }{\sigma^A_{\min}\rho^2}
+ \frac{30qL^2d^2\nu^2}{b\sigma^A_{\min}\rho^2} + \frac{40L^2d\mu^2}{\sigma^A_{\min}\rho^2},
\end{align}
where the first inequality follows by Lemma \ref{lem:M2}, and the last inequality is due to the definition of $\theta_k$ and $b=qb$.

Let $\beta_{\max}=\max\{\beta_1,\beta_2,\beta_3\}$ with
\begin{align}
  \beta_1 = m\big(\rho^2\sigma^B_{\max}\sigma^A_{\max} + \rho^2(\sigma^B_{\max})^2 + \sigma^2_{\max}(H)\big), \ \beta_2 = 18(L^2 + \frac{\sigma^2_{\max}(G)}{\eta^2}), \
  \beta_3 = \frac{60L^2 }{\sigma^A_{\min} \rho^2} + \frac{5\sigma^2_{\max}(G) }{\sigma^A_{\min}\eta^2\rho^2}. \nonumber
\end{align}
Combining the above inequalities \eqref{eq:N1}, \eqref{eq:N2} and \eqref{eq:N3}, we have
\begin{align}
\frac{1}{K}\sum_{k=1}^K \mathbb{E}\big[ \mbox{dist}(0,\partial L(x_k,y_{[m]}^k,\lambda_k))^2\big] & \leq \frac{\beta_{\max}}{K}\sum_{k=1}^{K-1} \theta_k
 + 15\vartheta_2L^2d^2\nu^2 + \frac{40L^2d\mu^2}{\sigma^A_{\min}\rho^2} \nonumber \\
& \leq  \frac{3\beta_{\max}(\Phi_{0} - \Phi^*)}{K \gamma} + \big(\frac{3L}{4} + \frac{3qL}{2b}\big)\frac{3\beta_{\max}\vartheta_1d^2\nu^2}{\gamma} + \frac{6\beta_{\max}\vartheta_1Ld\mu^2}{\gamma} \nonumber \\
& \quad  + 15\vartheta_2L^2d^2\nu^2 + \frac{40L^2d\mu^2}{\sigma^A_{\min}\rho^2},
\end{align}
where $\vartheta_1 = L + \frac{20L}{\sigma^A_{\min}\rho}$, $\vartheta_2 = 1 + \frac{2q}{b\rho^2\sigma^A_{\min}}$,
and the second inequality holds by the above Lemma \ref{lem:M3}.

Let $b=dq$, $\eta = \frac{\alpha\sigma_{\min}(G)}{4L} \ (0<\alpha \leq 1)$ and $\rho = \frac{2\sqrt{237}\kappa_G L}{\sigma^A_{\min}\alpha}$. Since $m$ denotes the number of nonsmooth regularization functions, it is relatively small.
It is easily verified that
$\vartheta_1=O(1)$, $\vartheta_2=O(1)$, $\beta_{\max} = O(1)$ and $\gamma=O(1)$, which are independent on $n$ and $K$. Thus, we obtain
\begin{align}
\mathbb{E}\big[ \mbox{dist}(0,\partial L(x_\zeta,y_{[m]}^\zeta,\lambda_\zeta))^2\big] = \frac{1}{K}\sum_{k=1}^K \mathbb{E}\big[ \mbox{dist}(0,\partial L(x_k,y_{[m]}^k,\lambda_k))^2\big] \leq  O(\frac{1}{K})+ O(d^2\nu^2) + O(d\mu^2),
\end{align}
where $\{x_\zeta,y_{[m]}^{\zeta},\lambda_\zeta\}$ is chosen uniformly randomly from $\{x_{k},y_{[m]}^{k},\lambda_k\}_{k=0}^{K-1}$.

\end{proof}

\subsection{ Convergence Analysis of ZOO-ADMM+(CooGE) Algorithm }
\label{Appendix:A3}
In this subsection, we study the convergence properties of the ZOO-ADMM+(CooGE) Algorithm.
We begin with giving some useful lemmas as follows:

\begin{lemma} \label{lem:H1}
Suppose the zeroth-order gradient $v_k$ be generated from the Algorithm \ref{alg:3}, we have
\begin{align}
 \mathbb{E}\| v_k - \nabla f(x_k) \|^2\leq \frac{6L^2}{b_2} \sum_{i=c_kq}^{k-1} \mathbb{E}\|x_{i+1}-x_i\|^2 + \frac{12qL^2d\mu^2}{b_2} + 2 L^2d\mu^2 + \frac{6}{b_1} (2L^2d\mu^2 + \sigma^2).
\end{align}
\end{lemma}
\begin{proof}
Since the online problem \eqref{eq:1} can be viewed as
having infinite samples, \emph{i.e.,} $n=+\infty$, we have $\mathcal{I}(|\mathcal{S}_1|<n)=1$.
In Algorithm \ref{alg:1}, we have $|S_1|=b_1$ and $|S_2|=b_2$.  According to Lemma \ref{lem:A3}, and
Assumption \ref{ass:5}, i.e., $\mathbb{E} \|\nabla f(x;\xi)-\nabla f(x)\|^2 \leq \sigma^2$ for all $x \in \mathbb{R}^d$, we have
\begin{align} \label{eq:H1}
\mathbb{E}\| v_k - \hat{\nabla}_{\texttt{coo}}f(x_k) \|^2 & \leq \frac{3L^2}{b_2} \sum_{i=c_kq}^{k-1} \mathbb{E}\|x_{i+1}-x_i\|^2 + (k-c_k q)\frac{6L^2d\mu^2}{b_2}
+ \frac{3\mathcal{I}(|\mathcal{S}_1|<n)}{b_1} (2L^2d\mu^2 + \sigma^2) \nonumber \\
& \leq \frac{3L^2}{b_2} \sum_{i=c_kq}^{k-1} \mathbb{E}\|x_{i+1}-x_i\|^2 + \frac{6 q L^2d\mu^2}{b_2} + \frac{3}{b_1} (2L^2d\mu^2 + \sigma^2),
\end{align}
where the second inequality follows by $k-c_kq \leq q$, and $\mathcal{I}(|\mathcal{S}_1|<n)=1$.

It follows that
\begin{align}
\mathbb{E}\| v_k - \nabla f(x_k) \|^2 & = \mathbb{E}\| v_k - \hat{\nabla}_{\texttt{coo}}f(x_k) + \hat{\nabla}_{\texttt{coo}}f(x_k) -  \nabla f(x_k) \|^2 \nonumber \\
 & \leq  2\mathbb{E}\| v_k - \hat{\nabla}_{\texttt{coo}}f(x_k)\|^2 + 2\mathbb{E}\| \hat{\nabla}_{\texttt{coo}}f(x_k) -  \nabla f(x_k) \|^2 \nonumber \\
 & \leq \frac{6L^2}{b_2} \sum_{i=c_kq}^{k-1} \mathbb{E}\|x_{i+1}-x_i\|^2 + \frac{12qL^2d\mu^2}{b_2} + 2 L^2d\mu^2 + \frac{6}{b_1} (2L^2d\mu^2 + \sigma^2),
\end{align}
where the second inequality holds by the above inequality \eqref{eq:H1} and the Lemma \ref{lem:A1}.
\end{proof}

\begin{lemma} \label{lem:H2}
Suppose the sequence $\{x_k,y_{[m]}^k,\lambda_k\}_{k=1}^K$ be generated from the Algorithm \ref{alg:3}, it holds that
\begin{align}
\mathbb{E}\|\lambda_{k+1} - \lambda_k\|^2 \leq & \frac{60L^2}{b_2\sigma^A_{\min}} \sum_{i=c_kq}^{k-1} \mathbb{E}\|x_{i+1}-x_i\|^2
+ (\frac{5L^2}{\sigma^A_{\min}}+\frac{5\sigma^2_{\max}(G)}{\sigma^A_{\min}\eta^2})\|x_k-x_{k-1}\|^2 \nonumber \\
& + \frac{5\sigma^2_{\max}(G)}{\sigma^A_{\min}\eta^2}\|x_{k+1}-x_{k}\|^2  + \frac{120qL^2d\mu^2}{b_2\sigma^A_{\min}} + \frac{20 L^2d\mu^2}{\sigma^A_{\min}}
+ \frac{60}{\sigma^A_{\min}b_1} (2L^2d\mu^2 + \sigma^2).
\end{align}
\end{lemma}
\begin{proof}
The proof of this lemma is the same as the proof of Lemma \ref{lem:B2}.
\end{proof}

\begin{lemma} \label{lem:H3}
Suppose the sequence $\{x_k,y_{[m]}^k,\lambda_k\}_{k=1}^K$ be generated from the Algorithm \ref{alg:3},
and define a \emph{Lyapunov} function $\Gamma_k$ as follows:
\begin{align}
 \Gamma_k = \mathbb{E} \big[\mathcal{L}_{\rho} (x_k,y_{[m]}^k,\lambda_k) + (\frac{5L^2}{\sigma^A_{\min}\rho}+\frac{5\sigma^2_{\max}(G)}{\sigma^A_{\min}\eta^2\rho})\|x_k-x_{k-1}\|^2 + \frac{12L^2}{\sigma^A_{\min}\rho b_2} \sum_{i=c_kq}^{k-1}\|x_{i+1}-x_i\|^2 \big]. \nonumber
\end{align}
Let $b_2=q$, $\eta=\frac{\alpha\sigma_{\min}(G)}{4L} \ (0<\alpha \leq 1)$ and $\rho = \frac{2\sqrt{237}\kappa_GL}{\sigma^A_{\min}\alpha}$,
then we have
\begin{align}
\frac{1}{K}\sum_{k=0}^{K-1} (\|x_{k+1} - x_{k}\|^2 + \sum_{j=1}^m \|y_j^k-y_j^{k+1}\|^2) \leq & \frac{\Gamma_{0} - \Gamma^*}{K\gamma}
 + \frac{7\vartheta_1L^2d\mu^2}{\gamma} + \frac{3\vartheta_1(2L^2d\mu^2 + \sigma^2)}{b_1\gamma},
\end{align}
where $\vartheta_1 = \frac{1}{L} + \frac{20}{\sigma^A_{\min}\rho}$,
$\gamma = \min(\chi,\sigma_{\min}^H)$ with $\chi\geq \frac{\sqrt{237}\kappa_GL}{2\alpha}$, and $\Gamma^*$ is a lower bound of the function $\Gamma_k$.
\end{lemma}

\begin{proof}
The proof of this lemma is the similar to Lemma \ref{lem:B3}.
By the optimal condition of step 9 in Algorithm \ref{alg:3},
we obtain
\begin{align} \label{eq:J1}
\mathcal{L}_{\rho} (x_k,y^{k+1}_{[m]},\lambda_k) \leq \mathcal{L}_{\rho} (x_k,y^k_{[m]},\lambda_k)
- \sigma_{\min}^H\sum_{j=1}^m \|y_j^k-y_j^{k+1}\|^2,
\end{align}
where $\sigma_{\min}^H=\min_{j\in[m]}\sigma_{\min}(H_j)$.

By using the optimal condition of step 10 in Algorithm \ref{alg:3},
we have
\begin{align} \label{eq:J2}
\mathcal{L}_{\rho} (x_{k+1},y_{[m]}^{k+1},\lambda_k) & \leq \mathcal{L}_{\rho} (x_k,y_{[m]}^{k+1},\lambda_k)
- (\frac{\sigma_{\min}(G)}{\eta}+\frac{\rho \sigma^A_{\min}}{2}-L ) \|x_{k+1} - x_{k}\|^2 \nonumber \\
& \quad + \frac{3L}{b_2} \sum_{i=c_kq}^{k-1} \mathbb{E}\|x_{i+1}-x_i\|^2 + \frac{6qLd\mu^2}{b_2} + Ld\mu^2 + \frac{3}{b_1L} (2L^2d\mu^2 + \sigma^2).
\end{align}

By using the step 10 in Algorithm \ref{alg:3}, we have
\begin{align} \label{eq:J3}
&\mathcal{L}_{\rho} (x_{k+1},y_{[m]}^{k+1},\lambda_{k+1}) -\mathcal{L}_{\rho} (x_{k+1},y_{[m]}^{k+1},\lambda_k)
= \frac{1}{\rho}\|\lambda_{k+1}-\lambda_k\|^2 \nonumber \\
& \leq \frac{60L^2}{b_2\sigma^A_{\min}\rho} \sum_{i=c_kq}^{k-1} \mathbb{E}\|x_{i+1}-x_i\|^2
+ (\frac{5L^2}{\sigma^A_{\min}\rho}+\frac{5\sigma^2_{\max}(G)}{\sigma^A_{\min}\eta^2\rho})\|x_k-x_{k-1}\|^2 \nonumber \\
& \quad + \frac{5\sigma^2_{\max}(G)}{\sigma^A_{\min}\eta^2\rho}\|x_{k+1}-x_{k}\|^2  + \frac{120qL^2d\mu^2}{b_2\sigma^A_{\min}\rho} + \frac{20 L^2d\mu^2}{\sigma^A_{\min}\rho}
+ \frac{60}{\sigma^A_{\min}\rho b_1} (2L^2d\mu^2 + \sigma^2),
\end{align}
where the above inequality holds by Lemma \ref{lem:H1}.

Combining \eqref{eq:J1}, \eqref{eq:J2} and \eqref{eq:J3}, we have
\begin{align} \label{eq:J4}
\mathcal{L}_{\rho} (x_{k+1},y_{[m]}^{k+1},\lambda_{k+1}) & \leq \mathcal{L}_{\rho} (x_k,y_{[m]}^k,\lambda_k)
 - \sigma_{\min}^H\sum_{j=1}^m \|y_j^k-y_j^{k+1}\|^2 + (\frac{3L}{b_2} + \frac{60L^2}{\sigma^A_{\min}b_2\rho})\sum_{i=c_kq}^{k-1}\mathbb{E}\|x_{i+1}-x_i\|^2 \nonumber \\
& \quad + (\frac{5L^2}{\sigma^A_{\min}\rho}+\frac{5\sigma^2_{\max}(G)}{\sigma^A_{\min}\eta^2\rho})\|x_k-x_{k-1}\|^2 + \big( \frac{6qL}{b_2} + L + \frac{120qL^2}{b_2\sigma^A_{\min}\rho}
+ \frac{20 L^2}{\sigma^A_{\min}\rho} \big) d\mu^2\nonumber \\
& \quad - (\frac{\sigma_{\min}(G)}{\eta}+\frac{\rho \sigma^A_{\min}}{2}-L - \frac{5\sigma^2_{\max}(G)}{\sigma^A_{\min}\eta^2\rho}) \|x_{k+1} - x_{k}\|^2 \nonumber \\
& \quad + (\frac{3}{b_1 L } + \frac{60}{\sigma^A_{\min}\rho b_1}) (2L^2d\mu^2 + \sigma^2).
\end{align}

Next, we define a useful \emph{Lyapunov} function $\Gamma_k$ as follows:
\begin{align}
\Gamma_k = \mathbb{E} \big[\mathcal{L}_{\rho} (x_k,y_{[m]}^k,\lambda_k) + (\frac{5L^2}{\sigma^A_{\min}\rho}+\frac{5\sigma^2_{\max}(G)}{\sigma^A_{\min}\eta^2\rho})\|x_k-x_{k-1}\|^2 + \frac{12L^2}{\sigma^A_{\min}\rho b_2} \sum_{i=c_kq}^{k-1}\|x_{i+1}-x_i\|^2 \big].
\end{align}
Following the above lemma \ref{lem:B3}, it is easily verified that the function $\Gamma_k$ is bounded from below.
Let $\Gamma^*$ denotes a lower bound of function $\Gamma_k$.
By the inequality \eqref{eq:J4}, we have
\begin{align} \label{eq:J5}
\Gamma_{k+1} & = \mathbb{E} \big[ \mathcal{L}_{\rho} (x_{k+1},y_{[m]}^{k+1},\lambda_{k+1}) + (\frac{5L^2}{\sigma^A_{\min}\rho} + \frac{5\sigma^2_{\max}(G)}{\sigma^A_{\min}\eta^2\rho})\|x_{k+1}-x_{k}\|^2
+ \frac{12L^2}{\sigma^A_{\min}\rho b_2} \sum_{i=c_kq}^{k}\|x_{i+1}-x_i\|^2 \big] \nonumber \\
& \leq \mathbb{E} \big[ \mathcal{L}_{\rho} (x_k,y_{[m]}^k,\lambda_k) + (\frac{5L^2}{\sigma^A_{\min}\rho}+\frac{5\sigma^2_{\max}(G)}{\sigma^A_{\min}\eta^2\rho})\mathbb{E}\|x_{k}-x_{k-1}\|^2
+ \frac{12L^2}{\sigma^A_{\min}\rho b_2} \sum_{i=c_kq}^{k-1}\|x_{i+1}-x_i\|^2  \big]\nonumber \\
& \quad - \sigma_{\min}^H\sum_{j=1}^m \|y_j^k-y_j^{k+1}\|^2 - (\frac{\sigma_{\min}(G)}{\eta}+\frac{\rho \sigma^A_{\min}}{2}-L - \frac{10\sigma^2_{\max}(G)}{\sigma^A_{\min}\eta^2\rho}-\frac{5L^2}{\sigma^A_{\min}\rho}
 - \frac{12L^2}{\sigma^A_{\min}\rho b_2})\mathbb{E} \|x_{k+1} - x_{k}\|^2 \nonumber \\
& \quad + (\frac{3L}{b_2} + \frac{60L^2}{\sigma^A_{\min}b_2\rho})\sum_{i=c_kq}^{k-1}\mathbb{E}\|x_{i+1}-x_i\|^2 + \big( \frac{6qL}{b_2} + L + \frac{120qL^2}{b_2\sigma^A_{\min}\rho}
+ \frac{20 L^2}{\sigma^A_{\min}\rho} \big) d\mu^2 + (\frac{3}{b_1 L } + \frac{60}{\sigma^A_{\min}\rho b_1}) (2L^2d\mu^2 + \sigma^2) \nonumber \\
& \leq \Gamma_k - \sigma_{\min}^H\sum_{j=1}^m \|y_j^k-y_j^{k+1}\|^2 - (\frac{\sigma_{\min}(G)}{\eta}+\frac{\rho \sigma^A_{\min}}{2}- L
- \frac{10\sigma^2_{\max}(G)}{\sigma^A_{\min}\eta^2\rho}-\frac{5L^2}{\sigma^A_{\min}\rho}
 - \frac{12L^2}{\sigma^A_{\min}\rho b_2})\mathbb{E} \|x_{k+1} - x_{k}\|^2 \nonumber \\
& \quad+ (\frac{3L}{b_2} + \frac{60L^2}{\sigma^A_{\min}b_2\rho})\sum_{i=c_kq}^{k-1}\mathbb{E}\|x_{i+1}-x_i\|^2 + \big( \frac{6qL}{b_2} + L + \frac{120qL^2}{b_2\sigma^A_{\min}\rho}
+ \frac{20L^2}{\sigma^A_{\min}\rho} \big) d\mu^2 \nonumber \\
& \quad + (\frac{3}{b_1 L } + \frac{60}{\sigma^A_{\min}\rho b_1}) (2L^2d\mu^2 + \sigma^2).
\end{align}

Telescoping inequality \eqref{eq:J5} over $k$ from $c_kq$ to $k$, we have
\begin{align} \label{eq:J6}
\Gamma_{k+1} & \leq \Gamma_{c_kq} - (\frac{\sigma_{\min}(G)}{\eta}+\frac{\rho \sigma^A_{\min}}{2}-L - \frac{10\sigma^2_{\max}(G)}{\sigma^A_{\min}\eta^2\rho}-\frac{5L^2}{\sigma^A_{\min}\rho}
  -\frac{12L^2}{\sigma^A_{\min}\rho b_2} ) \sum_{t=c_kq}^k \mathbb{E}\|x_{t+1} - x_{t}\|^2 \nonumber \\
& \quad - \sigma_{\min}^H\sum_{t=c_kq}^{k}\sum_{j=1}^m \|y_j^t-y_j^{t+1}\|^2 + (\frac{3L}{b_2} + \frac{60L^2}{\sigma^A_{\min}b_2\rho})\sum_{t=c_kq}^{k}\sum_{i=c_kq}^{k-1}\mathbb{E}\|x_{i+1}-x_i\|^2 \nonumber \\
& \quad + \sum_{t=c_kq}^{k}\big( \frac{6qL}{b_2} + L + \frac{120qL^2}{b_2\sigma^A_{\min}\rho} + \frac{20 L^2}{\sigma^A_{\min}\rho} \big) d\mu^2
+ \sum_{t=c_kq}^{k}(\frac{3}{b_1 L } + \frac{60}{\sigma^A_{\min}\rho b_1}) (2L^2d\mu^2 + \sigma^2) \nonumber \\
& \leq \Gamma_{c_kq} - (\frac{\sigma_{\min}(G)}{\eta}+\frac{\rho \sigma^A_{\min}}{2}-L - \frac{10\sigma^2_{\max}(G)}{\sigma^A_{\min}\eta^2\rho}-\frac{5L^2}{\sigma^A_{\min}\rho}
 -\frac{12L^2}{\sigma^A_{\min}\rho b_2}) \sum_{i = c_kq}^k \mathbb{E}\|x_{i+1} - x_{i}\|^2 \nonumber \\
&  \quad - \sigma_{\min}^H\sum_{i=c_kq}^{k-1}\sum_{j=1}^m \|y_j^i-y_j^{i+1}\|^2 + (\frac{3Lq}{b_2} + \frac{60L^2q}{\sigma^A_{\min}b_2\rho})\sum_{i=c_kq}^{k}\mathbb{E}\|x_{i+1}-x_i\|^2 \nonumber \\
& \quad + \big( \frac{6qL}{b_2} + L + \frac{120qL^2}{b_2\sigma^A_{\min}\rho} + \frac{20 L^2}{\sigma^A_{\min}\rho} \big) qd\mu^2 + q(\frac{3}{b_1 L } + \frac{60}{\sigma^A_{\min}\rho b_1}) (2L^2d\mu^2 + \sigma^2)\nonumber \\
& = \Gamma_{c_kq} \!-\! \underbrace{ (\frac{\sigma_{\min}(G)}{\eta} \!+\! \frac{\rho\sigma^A_{\min}}{2}\!-\!L \!-\! \frac{10\sigma^2_{\max}(G)}{\sigma^A_{\min}\eta^2\rho} \!-\! \frac{5L^2}{\sigma^A_{\min}\rho}
\!-\! \frac{12L^2}{\sigma^A_{\min}\rho b_2} \!-\! \frac{3Lq}{b_2} \!-\! \frac{60L^2q}{\sigma^A_{\min}b_2\rho} )}_{\chi}  \sum_{i=c_kq}^k \mathbb{E}\|x_{i+1} \!-\! x_{i}\|^2 \nonumber \\
& \quad - \sigma_{\min}^H\sum_{i=c_kq}^{k-1}\sum_{j=1}^m \|y_j^i-y_j^{i+1}\|^2 + \big( \frac{6qL}{b_2} + L + \frac{120qL^2}{b_2\sigma^A_{\min}\rho}
+ \frac{20L^2}{\sigma^A_{\min}\rho} \big) qd\mu^2 \nonumber \\
& \quad + q(\frac{3}{b_1 L } + \frac{60}{\sigma^A_{\min}\rho b_1}) (2L^2d\mu^2 + \sigma^2).
\end{align}

Let $b_2=q$, $\eta = \frac{\alpha\sigma_{\min}(G)}{4L} \ (0<\alpha \leq 1)$ and
$\rho = \frac{2\sqrt{237}\kappa_GL}{\sigma^A_{\min}\alpha}$,
we obtain $\chi \geq \frac{\sqrt{237}\kappa_GL}{2\alpha}$.
Telescoping inequality \eqref{eq:J6} over $k$ from $0$ to $K$, we have
\begin{align}
\Gamma_{K} - \Gamma_{0} & = \Gamma_{q} - \Gamma_{0} + \Gamma_{2q} - \Gamma_{q} + \cdots + \Gamma_{K} - \Gamma_{c_Kq} \nonumber \\
& \leq - \sum_{i=0}^{q-1} (\chi\|x_{i+1} - x_{i}\|^2 + \sigma_{\min}^H\sum_{j=1}^m \|y_j^i-y_j^{i+1}\|^2) - \sum_{i=q}^{2q-1} ( \chi\|x_{i+1} - x_{i}\|^2
+ \sigma_{\min}^H\sum_{j=1}^m \|y_j^i-y_j^{i+1}\|^2) \nonumber \\
& \quad - \cdots - \sum_{i=c_Kq}^{K-1} ( \chi \|x_{i+1} - x_{i}\|^2  + \sigma_{\min}^H\sum_{j=1}^m \|y_j^i-y_j^{i+1}\|^2) + K d\mu^2\big( \frac{6qL}{b} + L + \frac{120qL^2}{b\sigma^A_{\min}\rho}
+ \frac{20 L^2}{\sigma^A_{\min}\rho} \big) \nonumber \\
& \quad + K(\frac{3}{b_1 L } + \frac{60}{\sigma^A_{\min}\rho b_1}) (2L^2d\mu^2 + \sigma^2) \nonumber \\
& = - \sum_{i=0}^{K-1} (\chi \|x_{i+1} - x_{i}\|^2+ \sigma_{\min}^H\sum_{j=1}^m \|y_j^i-y_j^{i+1}\|^2) + K d\mu^2\big( \frac{6qL}{b_2} + L + \frac{120qL^2}{b_2\sigma^A_{\min}\rho}
+ \frac{20 L^2}{\sigma^A_{\min}\rho} \big) \nonumber \\
& \quad + K(\frac{3}{b_1 L } + \frac{60}{\sigma^A_{\min}\rho b_1}) (2L^2d\mu^2 + \sigma^2).
\end{align}
Since $b_2=q$, we obtain
\begin{align}
\frac{1}{K}\sum_{k=0}^{K-1} (\|x_{k+1} - x_{k}\|^2 + \sum_{j=1}^m \|y_j^k-y_j^{k+1}\|^2) \leq & \frac{\Gamma_{0} - \Gamma^*}{K\gamma}
 + \frac{7\vartheta_1L^2d\mu^2}{\gamma} + \frac{3\vartheta_1(2L^2d\mu^2 + \sigma^2)}{b_1\gamma},
\end{align}
where $\vartheta_1 = \frac{1}{L} + \frac{20}{\sigma^A_{\min}\rho}$ and $\gamma = \min(\chi,\sigma_{\min}^H)$ with $\chi\geq \frac{\sqrt{237}\kappa_GL}{2\alpha}$.

\end{proof}

Next, based on the above lemmas, we study the convergence properties of ZOO-ADMM+(CooGE) algorithm.
First, we define a useful variable $\theta_k = \mathbb{E}\big[ \|x_{k+1}-x_{k}\|^2+\|x_{k}-x_{k-1}\|^2+\frac{1}{q}\sum_{i=c_kq}^k \|x_{i+1}-x_i\|^2 + \sum_{j=1}^m \|y_j^k-y_j^{k+1}\|^2 \big]$.

\begin{theorem} \label{th:A3}
(Restatement of Theorem \ref{th:3})
Suppose the sequence $\{x_k,y_{[m]}^k,\lambda_k)_{k=1}^K$ be generated from the Algorithm \ref{alg:3}. Under the above assumptions,
the function $\Gamma_k$ is a lower bound.
Further, let $b_2=q$, $\eta = \frac{\alpha\sigma_{\min}(G)}{4L} \ (0<\alpha \leq 1)$ and $\rho = \frac{2\sqrt{237} \kappa_G L}{\sigma^A_{\min}\alpha}$,
we have
\begin{align}
\mathbb{E}\big[ \mbox{dist}(0,\partial L(x_\zeta,y_{[m]}^\zeta,\lambda_\zeta))^2\big] =\frac{1}{K}\sum_{k=1}^K \mathbb{E}\big[ \mbox{dist}(0,\partial L(x_k,y_{[m]}^k,\lambda_k))^2\big] \leq  O(\frac{1}{K}) + O(d\mu^2)+ O(\frac{1}{b_1}),
\end{align}
where $\{x_\zeta,y_{[m]}^{\zeta},\lambda_\zeta\}$ is chosen uniformly randomly from $\{x_{k},y_{[m]}^{k},\lambda_k\}_{k=0}^{K-1}$.
It implies that given
$K = O(\frac{1}{\epsilon})$, $b_1 = O(\frac{1}{\epsilon})$ and $\mu = O(\sqrt{\frac{\epsilon}{d}})$,
then $(x_{k^*},y^{k^*}_{[m]},\lambda_{k^*})$ is an $\epsilon$-approximate stationary point of the problem \eqref{eq:1},
where $k^* = \mathop{\arg\min}_{k}\theta_{k}$.
\end{theorem}

\begin{proof}
By the optimal condition of the step 9 in Algorithm \ref{alg:3}, we have, for all $j\in [m]$
\begin{align} \label{eq:K1}
\mathbb{E}\big[\mbox{dist}(0,\partial_{y_j} L(x,y_{[m]},\lambda))^2\big]_{k+1} & = \mathbb{E}\big[\mbox{dist} (0, \partial \psi_j(y_j^{k+1})-B_j^T\lambda_{k+1})^2\big] \nonumber \\
& = \|B_j^T\lambda_k -\rho B_j^T(Ax_k + \sum_{i=1}^jB_iy_i^{k+1} + \sum_{i=j+1}^mB_iy_i^{k} -c) - H_j(y_j^{k+1}-y_j^k) -B_j^T\lambda_{k+1}\|^2 \nonumber \\
& = \|\rho B_j^TA(x_{k+1}-x_{k}) + \rho B_j^T \sum_{i=j+1}^m B_i (y_i^{k+1}-y_i^{k})- H_j(y_j^{k+1}-y_j^k) \|^2 \nonumber \\
& \leq m\rho^2\sigma^{B_j}_{\max}\sigma^A_{\max}\|x_{k+1}\!-\!x_k\|^2 \!+\! m\rho^2\sigma^{B_j}_{\max}\sum_{i=j+1}^m \sigma^{B_i}_{\max}\|y_i^{k+1}\!-\!y_i^{k}\|^2 \!+\! m\sigma^2_{\max}(H_j)\|y_j^{k+1}\!-\!y_j^k\|^2\nonumber \\
& \leq m\big(\rho^2\sigma^B_{\max}\sigma^A_{\max} + \rho^2(\sigma^B_{\max})^2 + \sigma^2_{\max}(H)\big) \theta_{k},
\end{align}

By the step 10 of Algorithm \ref{alg:3}, we have
\begin{align} \label{eq:K2}
\mathbb{E}\big[\mbox{dist}(0,\nabla_x L(x,y_{[m]},\lambda))^2\big]_{k+1} & = \mathbb{E}\|A^T\lambda_{k+1}-\nabla f(x_{k+1})\|^2  \nonumber \\
& = \mathbb{E}\|v_k - \nabla f(x_{k+1}) - \frac{G}{\eta} (x_k - x_{k+1})\|^2 \nonumber \\
& = \mathbb{E}\|v_k - \nabla f(x_{k}) + \nabla f(x_{k})- \nabla f(x_{k+1}) - \frac{G}{\eta}(x_k-x_{k+1})\|^2  \nonumber \\
& \leq 3\mathbb{E}\|v_k \!-\! \nabla f(x_{k})\|^3 \!+\!3\mathbb{E}\|\nabla f(x_{k})\!-\! \nabla f(x_{k+1})\|^2 \!+\! \frac{3\sigma^2_{\max}(G)}{\eta^2}\|x_k\!-\!x_{k+1}\|^2  \nonumber \\
& \leq \frac{18L^2}{b_2} \sum_{i=c_kq}^{k-1} \mathbb{E}\|x_{i+1}-x_i\|^2 + \frac{36qL^2d\mu^2}{b_2} + 6L^2d\mu^2 + \frac{18}{b_1} (2L^2d\mu^2 + \sigma^2) \nonumber \\
& \quad + 3(L^2+ \frac{\sigma^2_{\max}(G)}{\eta^2})\|x_k-x_{k+1}\|^2 \nonumber \\
& \leq 18(L^2+ \frac{\sigma^2_{\max}(G)}{\eta^2})\theta_{k} + 42L^2d\mu^2 + \frac{18}{b_1} (2L^2d\mu^2 + \sigma^2),
\end{align}
where the second inequality holds by Lemma \ref{lem:H1}, and the last inequality is due to the definition of $\theta_k$ and $b_2=q$.

By the step 11 of Algorithm \ref{alg:3}, we have
\begin{align} \label{eq:K3}
\mathbb{E}\big[\mbox{dist}(0,\nabla_\lambda L(x,y_{[m]},\lambda))^2\big]_{k+1}  & = \mathbb{E}\|Ax_{k+1}+ \sum_{j=1}^m B_jy_j^{k+1}-c\|^2 \nonumber \\
&= \frac{1}{\rho^2} \mathbb{E} \|\lambda_{k+1}-\lambda_k\|^2  \nonumber \\
& \leq \frac{60L^2}{b_2\sigma^A_{\min}\rho^2} \sum_{i=c_kq}^{k-1} \mathbb{E}\|x_{i+1}\!-\!x_i\|^2
+ (\frac{5L^2}{\sigma^A_{\min}\rho^2} \!+\!\frac{5\sigma^2_{\max}(G)}{\sigma^A_{\min}\eta^2\rho^2})\|x_k\!-\!x_{k-1}\|^2  \nonumber \\
& \quad +\! \frac{5\sigma^2_{\max}(G)}{\sigma^A_{\min}\eta^2\rho^2}\|x_{k+1}\!-\!x_{k}\|^2+ \frac{120qL^2d\mu^2}{b_2\sigma^A_{\min}\rho^2} \!+\! \frac{20 L^2d\mu^2}{\sigma^A_{\min}\rho^2}
\!+\! \frac{60}{\sigma^A_{\min}\rho b_1} (2L^2d\mu^2 \!+\! \sigma^2) \nonumber \\
& \leq \big( \frac{60 L^2}{\sigma^A_{\min} \rho^2} + \frac{5\sigma^2_{\max}(G) }{\sigma^A_{\min}\eta^2\rho^2} \big) \theta_{k}
+ \frac{140L^2d\mu^2}{\sigma^A_{\min}\rho^2} + \frac{60}{\sigma^A_{\min}\rho^2 b_1} (2L^2d\mu^2 + \sigma^2),
\end{align}
where the first inequality follows by Lemma \ref{lem:H2}, and the last inequality is due to the definition of $\theta_k$ and $b_2=q$.

Let $\beta_{\max}=\max\{\beta_1,\beta_2,\beta_3\}$ with
\begin{align}
  \beta_1 \!=\! m\big(\rho^2\sigma^B_{\max}\sigma^A_{\max} \!+\! \rho^2(\sigma^B_{\max})^2 \!+\! \sigma^2_{\max}(H)\big), \ \beta_2 \!=\! 18(L^2 \!+\! \frac{\sigma^2_{\max}(G)}{\eta^2}), \
  \beta_3 \!=\! \frac{60L^2 }{\sigma^A_{\min} \rho^2} \!+\! \frac{5\sigma^2_{\max}(G) }{\sigma^A_{\min}\eta^2\rho^2}. \nonumber
\end{align}
Combining the above inequalities \eqref{eq:K1}, \eqref{eq:K2} and \eqref{eq:K3}, we have
\begin{align}
\frac{1}{K}\sum_{k=1}^K \mathbb{E}\big[ \mbox{dist}(0,\partial L(x_k,y_{[m]}^k,\lambda_k))^2\big] & \leq \frac{\beta_{\max}}{K}\sum_{k=1}^{K-1} \theta_k
 + 14\vartheta_2 L^2d\mu^2 + \frac{6\vartheta_2(2L^2d\mu^2 + \sigma^2)}{b_1} \nonumber \\
& \leq  \frac{3\beta_{\max}(\Gamma_{0} - \Gamma^*)}{K \gamma} + \frac{21\vartheta_1 L^2 d\mu^2}{\gamma}
+ \frac{9\vartheta_1(2L^2d\mu^2 + \sigma^2)}{b_1\gamma} \nonumber \\
& \quad + 14\vartheta_2 L^2d\mu^2 + \frac{6\vartheta_2(2L^2d\mu^2 + \sigma^2)}{b_1},
\end{align}
where the second inequality holds by the above Lemma \ref{lem:H3} and $\sum_{k=0}^{K-1}\sum_{i=c_kq}^k \|x_{i+1}-x_i\|^2 \leq q\sum_{k=0}^{K-1} \|x_{k+1}-x_k\|^2$; $\vartheta_1 = \frac{1}{L} + \frac{20}{\sigma^A_{\min}\rho}$, $\vartheta_2=\max\{3,\frac{10}{\sigma^A_{\min}\rho^2}\}$ and
$\gamma \geq \frac{\sqrt{237}\kappa_GL}{2\alpha}$.

Let $\eta = \frac{\alpha\sigma_{\min}(G)}{4L} \ (0<\alpha \leq 1)$ and $\rho = \frac{2\sqrt{237} \kappa_G L}{\sigma^A_{\min}\alpha}$. Since $m$ denotes the number of nonsmooth regularization functions, it is relatively small.
It is easily verified that
$\vartheta_1=O(1)$, $\vartheta_2=O(1)$, $\beta_{\max} = O(1)$ and $\gamma=O(1)$, which are independent on $n$ and $K$. Thus, we obtain
\begin{align}
\mathbb{E}\big[ \mbox{dist}(0,\partial L(x_\zeta,y_{[m]}^\zeta,\lambda_\zeta))^2\big] = \frac{1}{K}\sum_{k=1}^K \mathbb{E}\big[ \mbox{dist}(0,\partial L(x_k,y_{[m]}^k,\lambda_k))^2\big] \leq  O(\frac{1}{K})+ O(d\mu^2)+ O(\frac{1}{b_1}),
\end{align}
where $\{x_\zeta,y_{[m]}^{\zeta},\lambda_\zeta\}$ is chosen uniformly randomly from $\{x_{k},y_{[m]}^{k},\lambda_k\}_{k=0}^{K-1}$.
\end{proof}

\subsection{ Convergence Analysis of ZOO-ADMM+(CooGE+UniGE) Algorithm }
\label{Appendix:A4}
In this subsection, we analyze convergence of the ZOO-ADMM+(CooGE+UniGE) algorithm.
We first provide an upper bound of variance of stochastic zeroth-order gradient $v_k$.

\begin{lemma} \label{lem:P1}
 Suppose the zeroth-order stochastic gradient $v_t$ be generated from the Algorithm \ref{alg:4}, we have
 \begin{align}
  \mathbb{E} \|\nabla f(x_k)-v_k\|^2 \leq \frac{6dL^2}{b_2}\sum_{i=c_kq}^{k-1}\mathbb{E}\|x_{i+1}-x_i\|^2 + \frac{3\nu^2 L^2 d^2}{2} + \frac{3qL^2d^2\nu^2}{b_2}  + \frac{8\sigma^2}{b_1} + 8L^2d\mu^2.
 \end{align}
\end{lemma}

\begin{proof}
 We first define $f_{\nu}(x)=\mathbb{E}_{u\sim U_B}[f(x+\nu u)]$ be a smooth approximation of $f(x)$,
 where $U_B$ is the uniform distribution over the $d$-dimensional unit Euclidean ball $B$. By the Lemma \ref{lem:A3},
 we have $\mathbb{E}_{(u,\xi)}[\hat{\nabla}_{\texttt{uni}} f(x;\xi)]=f_{\nu}(x)$.

 Next, we give an upper bound of $\mathbb{E}\|v_k - \nabla f_{\nu}(x_{k})\|^2$. By the definition of $v_k$,
 we have
 \begin{align} \label{eq:P1}
  &\mathbb{E}\| \nabla f_{\nu}(x_{k})-v_k\|^2 \nonumber \\
  & = \mathbb{E}\|\underbrace{\nabla f_{\nu}(x_{k}) -\nabla f_{\nu}(x_{k-1}) -\frac{1}{b_2} \sum_{j \in \mathcal{S}_2} [\hat{\nabla}_{\texttt{uni}} f(x_k;\xi_j) - \hat{\nabla}_{\texttt{uni}} f(x_{k-1};\xi_j)] }_{=T_3} + \underbrace{\nabla f_{\nu}(x_{k-1})- v_{k-1}}_{=T_4}\|^2 \nonumber \\
  & = \mathbb{E}\|\nabla f_{\nu}(x_{k}) -\nabla f_{\nu}(x_{k-1}) -\frac{1}{b_2} \sum_{j \in \mathcal{S}_2} [\hat{\nabla}_{\texttt{uni}} f(x_k;\xi_j) - \hat{\nabla}_{\texttt{uni}} f(x_{k-1};\xi_j)] \|^2
  + \mathbb{E}\| \nabla f_{\nu}(x_{k-1})- v_{k-1}\|^2 \nonumber \\
  & = \frac{1}{b_2^2} \sum_{j \in \mathcal{S}_2}\mathbb{E}\| \nabla f_{\nu}(x_{k}) -\nabla f_{\nu}(x_{k-1}) - \hat{\nabla}_{\texttt{uni}} f(x_k;\xi_j) + \hat{\nabla}_{\texttt{uni}} f(x_{k-1};\xi_j) \|^2
  + \mathbb{E}\| \nabla f_{\nu}(x_{k-1}) - v_{k-1}\|^2 \nonumber \\
  & \leq\frac{1}{b_2^2} \sum_{j \in \mathcal{S}_2}\mathbb{E}\| \hat{\nabla}_{\texttt{uni}} f(x_k;\xi_j) - \hat{\nabla}_{\texttt{uni}} f(x_{k-1};\xi_j) \|^2
  + \mathbb{E}\| \nabla f_{\nu}(x_{k-1}) - v_{k-1}\|^2 \nonumber \\
  & \leq \frac{1}{b_2} \big( 3dL^2\|x_{k}-x_{k-1}\|^2 + \frac{3L^2d^2\nu^2}{2} \big) +  \mathbb{E}\| \nabla f_{\nu}(x_{k-1})- x_{k-1}\|^2,
 \end{align}
 where the second equality follows by $\mathbb{E}[T_3]=0$ and $T_4$ is independent to $\mathcal{S}_2$, and
 the third equality holds by the above Lemma \ref{lem:A0}, and the second inequality holds by the Lemma \ref{lem:A3}.

 When $k=c_kq$, we have $v_k = \frac{1}{b_1}\sum_{j\in\mathcal{S}_1}\hat{\nabla}_{\texttt{coo}} f(x_{k};\xi_j)$.
 Then we have
 \begin{align}
  \mathbb{E}\|\nabla f_{\nu}(x_{c_kq})-v_{c_kq}\|^2 & = \mathbb{E}\|\nabla f_{\nu}(x_{c_kq}) - \nabla f(x_{c_kq}) + \nabla f(x_{c_kq})-\frac{1}{b_1}\sum_{j\in\mathcal{S}_1}\hat{\nabla}_{\texttt{coo}} f(x_{c_kq};\xi_j)\|^2 \nonumber \\
  & \leq 2\mathbb{E}\|\nabla f_{\nu}(x_{c_kq}) - \nabla f(x_{c_kq})\|^2 + 2\mathbb{E}\|\nabla f(x_{c_kq})-\frac{1}{b_1}\sum_{j\in\mathcal{S}_1}\hat{\nabla}_{\texttt{coo}} f(x_{c_kq};\xi_j)\|^2 \nonumber \\
  & \leq \frac{L^2d^2\nu^2}{2} + 2\mathbb{E}\|\nabla f(x_{c_kq})-\frac{1}{b_1}\sum_{j\in\mathcal{S}_1}\nabla f(x_{k};\xi_j)+\frac{1}{b_1}\sum_{j\in\mathcal{S}_1}\nabla f(x_{c_kq};\xi_j)-\frac{1}{b_1}\sum_{j\in\mathcal{S}_1}\hat{\nabla}_{\texttt{coo}} f(x_{c_kq};\xi_j)\|^2 \nonumber \\
  & \leq \frac{L^2d^2\nu^2}{2} + 4\frac{1}{b_1^2}\sum_{j\in\mathcal{S}_1}\mathbb{E}\|\nabla f(x_{c_kq})-\nabla f(x_{c_kq};\xi_j)\|^2
  + \frac{4}{b_1}\sum_{j\in\mathcal{S}_1}\mathbb{E}\|\nabla f(x_{k};\xi_j)-\hat{\nabla}_{\texttt{coo}} f(x_{c_kq};\xi_j)\|^2 \nonumber \\
  & \leq \frac{L^2d^2\nu^2}{2} + \frac{4\sigma^2}{b_1} + 4L^2d\mu^2,
 \end{align}
 where the third inequality holds by Lemma \ref{lem:A0}.

 By recursion to \eqref{eq:P1}, we have
 \begin{align}
 \mathbb{E}\|\nabla f_{\nu}(x_{k})-v_k\|^2 & \leq (k-c_kq)\frac{3L^2d^2\nu^2}{2b_2}
 + \frac{3dL^2}{b_2}\sum_{i=c_kq}^{k-1}\mathbb{E}\|x_{i+1}-x_{i}\|^2 \big)
 + \mathbb{E}\| \nabla f_{\nu}(x_{c_kq})- v_{c_kq}\|^2 \nonumber \\
 & \leq \frac{3qL^2d^2\nu^2}{2b_2} + \frac{3dL^2}{b_2}\sum_{i=c_kq}^{k-1}\mathbb{E}\|x_{i+1}-x_{i}\|^2 + \frac{L^2d^2\nu^2}{2} + \frac{4\sigma^2}{b_1} + 4L^2d\mu^2,
 \end{align}
 where the last inequality holds by the Assumption 5, i.e., $\mathbb{E}\| v_{c_kq}-\nabla f_{\nu}(x_{c_kq})\|^2 \leq \frac{\sigma^2}{b_1}$ and $k-c_kq \leq q$.

 Finally, we have
 \begin{align}
  \mathbb{E} \|\nabla f(x_k)-v_k\|^2 &=\mathbb{E} \|\nabla f(x_k)- \nabla f_{\nu}(x_{k}) + \nabla f_{\nu}(x_{k})- v_k\|^2 \nonumber\\
  & \leq 2\mathbb{E} \|\nabla f(x_k)- \nabla f_{\nu}(x_{k})\|^2 + 2\mathbb{E} \| \nabla f_{\nu}(x_{k})- v_k\|^2 \nonumber\\
  & \leq \frac{3\nu^2 L^2 d^2}{2} + \frac{3qL^2d^2\nu^2}{b_2} + \frac{6dL^2}{b_2}\sum_{i=c_kq}^{k-1}\mathbb{E}\|x_{i+1}-x_{i}\|^2
  + \frac{8\sigma^2}{b_1} + 8L^2d\mu^2,
 \end{align}
 where the last inequality holds by the above Lemma \ref{lem:A3}.
\end{proof}

\begin{lemma} \label{lem:P2}
Suppose the sequence $\{x_k,y_{[m]}^k,\lambda_k\}_{k=1}^K$ be generated from the Algorithm \ref{alg:4}, it holds that
\begin{align}
\mathbb{E}\|\lambda_{k+1} - \lambda_k\|^2 \leq & \frac{60dL^2}{\sigma^A_{\min}b_2}\sum_{i=c_kq}^{k-1}\mathbb{E}\|x_{i+1}-x_i\|^2 + \big(\frac{5\sigma^2_{\max}(G)}{\sigma^A_{\min}\eta^2}+\frac{5L^2}{\sigma^A_{\min}}\big)\mathbb{E}\|x_k-x_{k-1}\|^2\nonumber \\
& + \frac{5\sigma^2_{\max}(G)}{\sigma^A_{\min}\eta^2}\|x_{k+1}-x_{k}\|^2
+ \frac{15\nu^2 L^2 d^2}{\sigma^A_{\min}} + \frac{30qL^2d^2\nu^2}{b_2\sigma^A_{\min}}  + \frac{80\sigma^2}{b_1\sigma^A_{\min}} + \frac{80L^2d\mu^2}{\sigma^A_{\min}}.
\end{align}
\end{lemma}
\begin{proof}
The proof of this lemma is the same as the proof of Lemma \ref{lem:B2}.
\end{proof}

\begin{lemma} \label{lem:P3}
Suppose the sequence $\{x_k,y_{[m]}^k,\lambda_k\}_{k=1}^K$ be generated from the Algorithm \ref{alg:4},
and define a \emph{Lyapunov} function $\Psi_k$ as follows:
\begin{align}
\Psi_k = \mathbb{E} \big[ \mathcal{L}_{\rho} (x_k,y_{[m]}^k,\lambda_k) + (\frac{5L^2}{\sigma^A_{\min}\rho}+\frac{5\sigma^2_{\max}(G)}{\sigma^A_{\min}\eta^2\rho})\|x_k-x_{k-1}\|^2 + \frac{12dL^2}{b_2\sigma^A_{\min}\rho}\sum_{i=c_kq}^{k-1}\mathbb{E}\|x_{i+1}-x_i\|^2 \big]. \nonumber
\end{align}
Let $b_2=qd$, $\eta=\frac{\alpha\sigma_{\min}(G)}{4L} \ (0<\alpha \leq 1)$ and $\rho = \frac{2\sqrt{237}\kappa_GL}{\sigma^A_{\min}\alpha}$,
we have
\begin{align} \label{eq:P9}
\frac{1}{K}\sum_{k=0}^{K-1} \big(\mathbb{E}\|x_{k+1} - x_{k}\|^2 + \sum_{j=1}^m \|y_j^k-y_j^{k+1}\|^2\big) & \leq \frac{\Psi_{0} - \Psi^*}{K\gamma}
+ \frac{\vartheta_1\vartheta_2d^2\nu^2}{\gamma} + \frac{\vartheta_1}{\gamma}\big(\frac{4\sigma^2}{b_1} + 4L^2d\mu^2\big),
\end{align}
where $\vartheta_1 = \frac{1}{L} + \frac{20}{\sigma^A_{\min}\rho}$, $\vartheta_2=\frac{3L^2}{4}+\frac{3qL^2}{2b_2}$, $\gamma = \min(\chi,\sigma_{\min}^H)$ with $\chi\geq \frac{\sqrt{237}\kappa_GL}{2\alpha}$,
and $\Psi^*$ is a lower bound of the function $\Psi_k$.
\end{lemma}

\begin{proof}
The proof of this lemma is the similar to the proof of Lemma \ref{lem:B3}.
Similarly, we can obtain
\begin{align} \label{eq:P3}
\mathcal{L}_{\rho} (x_k,y^{k+1}_{[m]},\lambda_k) \leq \mathcal{L}_{\rho} (x_k,y^k_{[m]},\lambda_k)
- \sigma_{\min}^H\sum_{j=1}^m \|y_j^k-y_j^{k+1}\|^2,
\end{align}
where $\sigma_{\min}^H=\min_{j\in[m]}\sigma_{\min}(H_j)$.

Similarly, we have
\begin{align}
0 & \leq \mathcal{L}_{\rho} (x_k,y_{[m]}^{k+1},\lambda_k) -  \mathcal{L}_{\rho} (x_{k+1},y_{[m]}^{k+1},\lambda_k)
- (\frac{\sigma_{\min}(G)}{\eta}+\frac{\rho \sigma^A_{\min}}{2}-L ) \|x_{k+1} - x_{k}\|^2 + \frac{1}{2L}\|v_k - \nabla f(x_k)\|^2 \nonumber \\
& \leq \mathcal{L}_{\rho} (x_k,y_{[m]}^{k+1},\lambda_k) -  \mathcal{L}_{\rho} (x_{k+1},y_{[m]}^{k+1},\lambda_k)
- (\frac{\sigma_{\min}(G)}{\eta}+\frac{\rho \sigma^A_{\min}}{2}-L ) \|x_{k+1} - x_{k}\|^2
+ \frac{3dL}{b_2}\sum_{i=c_kq}^{k-1}\mathbb{E}\|x_{i+1}-x_i\|^2 \nonumber \\
& \quad + \frac{3\nu^2 L d^2}{4} + \frac{3qLd^2\nu^2}{2b_2}  + \frac{4\sigma^2}{b_1L} + 4Ld\mu^2, \nonumber
\end{align}
where the final inequality holds by Lemma \ref{lem:P1}.
It follows that
\begin{align} \label{eq:P4}
\mathcal{L}_{\rho} (x_{k+1},y_{[m]}^{k+1},\lambda_k) & \leq \mathcal{L}_{\rho} (x_k,y_{[m]}^{k+1},\lambda_k)
- (\frac{\sigma_{\min}(G)}{\eta}+\frac{\rho \sigma^A_{\min}}{2}-L )\mathbb{E}\|x_{k+1} - x_{k}\|^2 \nonumber \\
& \quad + \frac{3dL}{b_2}\sum_{i=c_kq}^{k-1}\mathbb{E}\|x_{i+1}-x_i\|^2 + \frac{3\nu^2 L d^2}{4} + \frac{3qLd^2\nu^2}{2b_2}
+ \frac{4\sigma^2}{b_1L} + 4Ld\mu^2.
\end{align}

By using the step 10 in Algorithm \ref{alg:4}, we have
\begin{align} \label{eq:P5}
\mathcal{L}_{\rho} (x_{k+1},y_{[m]}^{k+1},\lambda_{k+1}) -\mathcal{L}_{\rho} (x_{k+1},y_{[m]}^{k+1},\lambda_k)
& = \frac{1}{\rho}\mathbb{E}\|\lambda_{k+1}-\lambda_k\|^2  \\
& \leq \frac{60dL^2}{\sigma^A_{\min}b_2\rho}\sum_{i=c_kq}^{k-1}\mathbb{E}\|x_{i+1}-x_i\|^2 + \big(\frac{5\sigma^2_{\max}(G)}{\sigma^A_{\min}\eta^2\rho}+\frac{5L^2}{\sigma^A_{\min}\rho}\big)\mathbb{E}\|x_k-x_{k-1}\|^2\nonumber \\
& \quad + \frac{5\sigma^2_{\max}(G)}{\sigma^A_{\min}\eta^2\rho}\mathbb{E}\|x_{k+1}-x_{k}\|^2 + \frac{15\nu^2 L^2 d^2}{\sigma^A_{\min}\rho} + \frac{30qL^2d^2\nu^2}{b_2\sigma^A_{\min}\rho}  \nonumber \\
&\quad + \frac{80\sigma^2}{b_1\sigma^A_{\min}\rho} + \frac{80L^2d\mu^2}{\sigma^A_{\min}\rho}, \nonumber
\end{align}
where the above inequality holds by Lemma \ref{lem:P1}.

Combining \eqref{eq:P3}, \eqref{eq:P4} and \eqref{eq:P5}, we have
\begin{align} \label{eq:P6}
\mathcal{L}_{\rho} (x_{k+1},y_{[m]}^{k+1},\lambda_{k+1}) & \leq \mathcal{L}_{\rho} (x_k,y_{[m]}^k,\lambda_k)
 \!-\! \sigma_{\min}^H\sum_{j=1}^m \|y_j^k-y_j^{k+1}\|^2 \!-\! \big(\frac{\sigma_{\min}(G)}{\eta}+\frac{\rho \sigma^A_{\min}}{2}-L - \frac{5\sigma^2_{\max}(G)}{\sigma^A_{\min}\eta^2\rho}\big)
 \mathbb{E}\|x_{k+1} - x_{k}\|^2 \nonumber \\
& \quad + \big(\frac{5L^2}{\sigma^A_{\min}\rho}+\frac{5\sigma^2_{\max}(G)}{\sigma^A_{\min}\eta^2\rho}\big)\mathbb{E}\|x_k-x_{k-1}\|^2 + \big(\frac{60dL^2}{\sigma^A_{\min}b_2\rho} + \frac{3dL}{b_2}  \big)\sum_{i=c_kq}^{k-1}\mathbb{E}\|x_{i+1}-x_i\|^2 \nonumber \\
& \quad + \big( \frac{3L}{4} + \frac{3qL}{2b_2} + \frac{15 L^2}{\sigma^A_{\min}\rho}+ \frac{30qL^2}{b_2\sigma^A_{\min}\rho}\big) d^2\nu^2+ (\frac{1}{L } + \frac{20}{\sigma^A_{\min}\rho})\big(\frac{4\sigma^2}{b_1} + 4L^2d\mu^2\big).
\end{align}

We define a useful \emph{Lyapunov} function $\Psi_k$ as follows:
\begin{align}
\Psi_k = \mathbb{E} \big[ \mathcal{L}_{\rho} (x_k,y_{[m]}^k,\lambda_k) + (\frac{5L^2}{\sigma^A_{\min}\rho}+\frac{5\sigma^2_{\max}(G)}{\sigma^A_{\min}\eta^2\rho})\|x_k-x_{k-1}\|^2 + \frac{12dL^2}{b_2\sigma^A_{\min}\rho}\sum_{i=c_kq}^{k-1}\mathbb{E}\|x_{i+1}-x_i\|^2 \big].
\end{align}
Following the above lemma \ref{lem:B3}, it is easily verified that the function $\Psi_k$ is bounded from below.
Let $\Psi^*$ denotes a lower bound of function $\Psi_k$.
It follows that
\begin{align} \label{eq:P7}
\Psi_{k+1} & = \mathbb{E} \big[ \mathcal{L}_{\rho} (x_{k+1},y_{[m]}^{k+1},\lambda_{k+1}) + (\frac{5L^2}{\sigma^A_{\min}\rho} + \frac{5\sigma^2_{\max}(G)}{\sigma^A_{\min}\eta^2\rho})\|x_{k+1}-x_{k}\|^2
+ \frac{12dL^2}{b_2\sigma^A_{\min}\rho}\sum_{i=c_kq}^{k}\mathbb{E}\|x_{i+1}-x_i\|^2 \big]  \\
& \leq \mathbb{E} \big[ \mathcal{L}_{\rho} (x_k,y_{[m]}^k,\lambda_k) + (\frac{5L^2}{\sigma^A_{\min}\rho}+\frac{5\sigma^2_{\max}(G)}{\sigma^A_{\min}\eta^2\rho})\|x_{k}-x_{k-1}\|^2
+ \frac{12dL^2}{b_2\sigma^A_{\min}\rho}\sum_{i=c_kq}^{k-1}\mathbb{E}\|x_{i+1}-x_i\|^2 \big]   \nonumber \\
& \quad - \sigma_{\min}^H\sum_{j=1}^m \|y_j^k-y_j^{k+1}\|^2 - \big(\frac{\sigma_{\min}(G)}{\eta}+\frac{\rho \sigma^A_{\min}}{2}-L - \frac{10\sigma^2_{\max}(G)}{\sigma^A_{\min}\eta^2\rho}-\frac{5L^2}{\sigma^A_{\min}\rho}
 - \frac{12dL^2}{b_2\sigma^A_{\min}\rho} \big) \mathbb{E}\|x_{k+1} - x_{k}\|^2 \nonumber \\
& \quad + \big(\frac{60dL^2}{\sigma^A_{\min}b_2\rho} + \frac{3dL}{b_2}  \big)\sum_{i=c_kq}^{k-1}\mathbb{E}\|x_{i+1}-x_i\|^2 + \big( \frac{3L}{4} + \frac{3qL}{2b_2} + \frac{15 L^2}{\sigma^A_{\min}\rho}+ \frac{30qL^2}{b_2\sigma^A_{\min}\rho}\big) d^2\nu^2+ (\frac{1}{L } + \frac{20}{\sigma^A_{\min}\rho})\big(\frac{4\sigma^2}{b_1} + 4L^2d\mu^2\big) \nonumber \\
& = \Psi_k - \sigma_{\min}^H\sum_{j=1}^m \|y_j^k-y_j^{k+1}\|^2 - \big(\frac{\sigma_{\min}(G)}{\eta}+\frac{\rho \sigma^A_{\min}}{2}-L - \frac{10\sigma^2_{\max}(G)}{\sigma^A_{\min}\eta^2\rho}-\frac{5L^2}{\sigma^A_{\min}\rho}
 - \frac{12dL^2}{b_2\sigma^A_{\min}\rho} \big) \mathbb{E}\|x_{k+1} - x_{k}\|^2 \nonumber \\
& \quad + \big(\frac{60dL^2}{\sigma^A_{\min}b_2\rho} + \frac{3dL}{b_2} \big)\sum_{i=c_kq}^{k-1}\mathbb{E}\|x_{i+1}-x_i\|^2 + \big( \frac{3L}{4} + \frac{3qL}{2b_2} + \frac{15 L^2}{\sigma^A_{\min}\rho}+ \frac{30qL^2}{b_2\sigma^A_{\min}\rho}\big) d^2\nu^2 + (\frac{1}{L} + \frac{20}{\sigma^A_{\min}\rho})\big(\frac{4\sigma^2}{b_1} + 4L^2d\mu^2\big),  \nonumber
\end{align}
where the first inequality holds by the inequality \eqref{eq:P6} and the equality
$$\sum_{i=c_kq}^{k}\mathbb{E}\|x_{i+1}-x_i\|^2 = \sum_{i=c_kq}^{k-1}\mathbb{E}\|x_{i+1}-x_i\|^2 + \mathbb{E}\|x_{k+1}-x_k\|^2.$$

Since $c_kq \leq k \leq (c_k+1)q -1$, and let $c_kq \leq t \leq (c_k+1)q-1$,
then telescoping inequality \eqref{eq:P7} over $k$ from $c_kq$ to $k$, we have
\begin{align} \label{eq:P8}
\Psi_{k+1} &
\leq \Psi_{c_kq} - \big(\frac{\sigma_{\min}(G)}{\eta}+\frac{\rho \sigma^A_{\min}}{2}-L - \frac{10\sigma^2_{\max}(G)}{\sigma^A_{\min}\eta^2\rho}-\frac{5L^2}{\sigma^A_{\min}\rho}
 - \frac{12dL^2}{b_2\sigma^A_{\min}\rho} \big) \sum_{t=c_kq}^k \mathbb{E}\|x_{t+1} - x_{t}\|^2 \nonumber \\
& \quad - \sigma_{\min}^H\sum_{t=c_kq}^{k}\sum_{j=1}^m \|y_j^t\!-\!y_j^{t+1}\|^2 + \big(\frac{60dL^2}{\sigma^A_{\min}b_2\rho} + \frac{3dL}{b_2} \big)
\sum_{t=c_kq}^{k}\sum_{i=c_kq}^{k-1}\mathbb{E}\|x_{i+1}-x_i\|^2 \nonumber \\
& \quad + \sum_{t=c_kq}^{k}\big( \frac{3L}{4} + \frac{3qL}{2b_2} + \frac{15 L^2}{\sigma^A_{\min}\rho}+ \frac{30qL^2}{b_2\sigma^A_{\min}\rho}\big) d^2\nu^2 + \sum_{t=c_kq}^{k}(\frac{1}{L} + \frac{20}{\sigma^A_{\min}\rho})\big(\frac{4\sigma^2}{b_1} + 4L^2d\mu^2\big) \nonumber \\
& \leq \Psi_{c_kq} - (\frac{\sigma_{\min}(G)}{\eta}+\frac{\rho \sigma^A_{\min}}{2}-L - \frac{10\sigma^2_{\max}(G)}{\sigma^A_{\min}\eta^2\rho}-\frac{5L^2}{\sigma^A_{\min}\rho}
 -\frac{12dL^2}{\sigma^A_{\min}\rho b_2}) \sum_{i = c_kq}^k \mathbb{E}\|x_{i+1} - x_{i}\|^2 \nonumber \\
& \quad - \sigma_{\min}^H\sum_{i=c_kq}^{k-1}\sum_{j=1}^m \|y_j^i-y_j^{i+1}\|^2 + \big(\frac{60qdL^2}{\sigma^A_{\min}b_2\rho} + \frac{3qdL}{b_2} \big)\sum_{i=c_kq}^{k}\mathbb{E}\|x_{i+1}-x_i\|^2 \nonumber \\
& \quad + \big( \frac{3L}{4} + \frac{3qL}{2b_2} + \frac{15 L^2}{\sigma^A_{\min}\rho}+ \frac{30qL^2}{b_2\sigma^A_{\min}\rho}\big) qd^2\nu^2
+ q(\frac{1}{L} + \frac{20}{\sigma^A_{\min}\rho})\big(\frac{4\sigma^2}{b_1} + 4L^2d\mu^2\big) \nonumber \\
& = \Psi_{c_kq} - \underbrace{ \big(\frac{\sigma_{\min}(G)}{\eta}+\frac{\rho \sigma^A_{\min}}{2}-L - \frac{10\sigma^2_{\max}(G)}{\sigma^A_{\min}\eta^2\rho}-\frac{5L^2}{\sigma^A_{\min}\rho}
- \frac{12dL^2}{b_2\sigma^A_{\min}\rho}-\frac{60qdL^2}{\sigma^A_{\min}b_2\rho} - \frac{3qdL}{b_2} \big)}_{\chi}  \sum_{i=c_kq}^k \|x_{i+1} - x_{i}\|^2 \nonumber \\
& \quad - \sigma_{\min}^H\sum_{i=c_kq}^{k-1}\sum_{j=1}^m \|y_j^i-y_j^{i+1}\|^2 + \big( \frac{3L}{4} + \frac{3qL}{2b_2} + \frac{15 L^2}{\sigma^A_{\min}\rho}+ \frac{30qL^2}{b_2\sigma^A_{\min}\rho}\big) qd^2\nu^2
+ q(\frac{1}{L} + \frac{20}{\sigma^A_{\min}\rho})\big(\frac{4\sigma^2}{b_1} + 4L^2d\mu^2\big),
\end{align}
where the second inequality holds by the fact that
\begin{align}
\sum_{j=c_kq}^{k}\sum_{i=c_kq}^{k-1}\mathbb{E}\|x_{i+1}-x_i\|^2 \leq \sum_{j=c_kq}^{k}\sum_{i=c_kq}^{k}\mathbb{E}\|x_{i+1}-x_i\|^2 \leq q\sum_{i=c_kq}^{k}\mathbb{E}\|x_{i+1}-x_i\|^2. \nonumber
\end{align}

Given $b_2=dq$, $\eta=\frac{\alpha\sigma_{\min}(G)}{4L} \ (0<\alpha \leq 1)$ and $\rho = \frac{2\sqrt{237}\kappa_GL}{\sigma^A_{\min}\alpha}$, similarly,
we have $\chi \geq \frac{\sqrt{237}\kappa_GL}{2\alpha}>0$.

Telescoping inequality \eqref{eq:P8} over $k$ from $0$ to $K$, we have
\begin{align}
\Psi_{K} - \Psi_{0}  &= \Psi_{q} - \Psi_{0} + \Psi_{2q} - \Psi_{q} + \cdots + \Psi_{K} - \Psi_{c_Kq}\nonumber \\
& \leq - \sum_{i=0}^{q-1} (\chi\|x_{i+1} - x_{i}\|^2 + \sigma_{\min}^H\sum_{j=1}^m \|y_j^i-y_j^{i+1}\|^2) - \sum_{i=q}^{2q-1} ( \chi\|x_{i+1} - x_{i}\|^2
+ \sigma_{\min}^H\sum_{j=1}^m \|y_j^i-y_j^{i+1}\|^2) - \cdots \nonumber \\
& \quad - \sum_{i=c_Kq}^{K-1} ( \chi \|x_{i+1} - x_{i}\|^2  + \sigma_{\min}^H\sum_{j=1}^m \|y_j^i-y_j^{i+1}\|^2) + K\big( \frac{3L}{4} + \frac{3qL}{2b_2} + \frac{15 L^2}{\sigma^A_{\min}\rho}+ \frac{30qL^2}{b_2\sigma^A_{\min}\rho}\big) d^2\nu^2 \nonumber \\
& \quad + K(\frac{1}{L} + \frac{20}{\sigma^A_{\min}\rho})\big(\frac{4\sigma^2}{b_1} + 4L^2d\mu^2\big) \nonumber \\
& = - \sum_{i=0}^{K-1} (\chi \|x_{i+1} - x_{i}\|^2+ \sigma_{\min}^H\sum_{j=1}^m \|y_j^i-y_j^{i+1}\|^2) + K\big( \frac{3L}{4} + \frac{3qL}{2b_2} + \frac{15 L^2}{\sigma^A_{\min}\rho}+ \frac{30qL^2}{b_2\sigma^A_{\min}\rho}\big) d^2\nu^2 \nonumber \\
& \quad + K(\frac{1}{L} + \frac{20}{\sigma^A_{\min}\rho})\big(\frac{4\sigma^2}{b_1} + 4L^2d\mu^2\big).
\end{align}
Thus, we have
\begin{align} \label{eq:P9}
\frac{1}{K}\sum_{k=0}^{K-1} \big(\mathbb{E}\|x_{k+1} - x_{k}\|^2 + \sum_{j=1}^m \|y_j^k-y_j^{k+1}\|^2\big) & \leq \frac{\Psi_{0} - \Psi^*}{K\gamma}
+ \frac{\vartheta_1\vartheta_2d^2\nu^2}{\gamma} + \frac{\vartheta_1}{\gamma}\big(\frac{4\sigma^2}{b_1} + 4L^2d\mu^2\big),
\end{align}
where $\vartheta_1 = \frac{1}{L} + \frac{20}{\sigma^A_{\min}\rho}$, $\vartheta_2=\frac{3L^2}{4}+\frac{3qL^2}{2b_2}$, $\gamma = \min(\chi,\sigma_{\min}^H)$ with $\chi\geq \frac{\sqrt{237}\kappa_GL}{2\alpha}$.

\end{proof}

Next, based on the above lemmas, we study the convergence properties of the ZOO-ADMM+(CooGE+UniGE) Algorithm.
First, we define a useful variable $\theta_k = \mathbb{E}\big[ \|x_{k+1}-x_{k}\|^2 + \|x_{k}-x_{k-1}\|^2 + \frac{1}{q}\sum_{i={c_kq}^k}\|x_{i+1}-x_{i}\|^2 + \sum_{j=1}^m \|y_j^k-y_j^{k+1}\|^2 \big]$.

\begin{theorem} \label{th:A4}
(Restatement of Theorem \ref{th:4})
Suppose the sequence $\{x_k,y_{[m]}^k,\lambda_k)_{k=1}^K$ be generated from the Algorithm \ref{alg:4}. The function $\Psi_k$ has a lower bound.
Further, let $b_2=qd$, $\eta = \frac{\alpha\sigma_{\min}(G)}{4L} \ (0<\alpha \leq 1)$ and $\rho = \frac{2\sqrt{237} \kappa_G L}{\sigma^A_{\min}\alpha}$,
we have
\begin{align}
\mathbb{E}\big[ \mbox{dist}(0,\partial L(x_\zeta,y_{[m]}^\zeta,\lambda_\zeta))^2\big] =\frac{1}{K}\sum_{k=1}^K \mathbb{E}\big[ \mbox{dist}(0,\partial L(x_k,y_{[m]}^k,\lambda_k))^2\big] \leq  O(\frac{1}{K})+ O(d^2\nu^2) + O(d\mu^2) + O(\frac{1}{b_1}),
\end{align}
where $\{x_\zeta,y_{[m]}^{\zeta},\lambda_\zeta\}$ is chosen uniformly randomly from $\{x_{k},y_{[m]}^{k},\lambda_k\}_{k=0}^{K-1}$.
It implies that given
$K = O(\frac{1}{\epsilon}), \ \nu = O(\frac{\sqrt{\epsilon}}{d}), \ \mu=O(\sqrt{\frac{\epsilon}{d}}), \ b_1 = O(\frac{1}{\epsilon})$
then $(x_{k^*},y^{k^*}_{[m]},\lambda_{k^*})$ is an $\epsilon$-approximate stationary point of the problem \eqref{eq:1}, where $k^* = \mathop{\arg\min}_{k}\theta_{k}$.
\end{theorem}
\begin{proof}
By the optimal condition of the step 9 in Algorithm \ref{alg:4}, we
have, for all $j\in [m]$
\begin{align} \label{eq:Q1}
\mathbb{E}\big[\mbox{dist}(0,\partial_{y_j} L(x,y_{[m]},\lambda))^2\big]_{k+1} & = \mathbb{E}\big[\mbox{dist} (0, \partial \psi_j(y_j^{k+1})-B_j^T\lambda_{k+1})^2\big] \nonumber \\
& = \|B_j^T\lambda_k -\rho B_j^T(Ax_k + \sum_{i=1}^jB_iy_i^{k+1} + \sum_{i=j+1}^mB_iy_i^{k} -c) - H_j(y_j^{k+1}-y_j^k)
-B_j^T\lambda_{k+1}\|^2 \nonumber \\
& = \|\rho B_j^TA(x_{k+1}-x_{k}) + \rho B_j^T \sum_{i=j+1}^m B_i (y_i^{k+1}-y_i^{k})- H_j(y_j^{k+1}-y_j^k) \|^2 \nonumber \\
& \leq m\rho^2\sigma^{B_j}_{\max}\sigma^A_{\max}\|x_{k+1}\!-\!x_k\|^2 + m\rho^2\sigma^{B_j}_{\max}\!\sum_{i=j+1}^m \!\sigma^{B_i}_{\max}\|y_i^{k+1}\!-\!y_i^{k}\|^2
+ m\sigma^2_{\max}(H_j)\|y_j^{k+1}\!-\!y_j^k\|^2\nonumber \\
& \leq m\big(\rho^2\sigma^B_{\max}\sigma^A_{\max} + \rho^2(\sigma^B_{\max})^2 + \sigma^2_{\max}(H)\big) \theta_{k},
\end{align}
where the first inequality follows by the inequality $\|\sum_{i=1}^r \alpha_i\|^2 \leq r\sum_{i=1}^r \|\alpha_i\|^2$.

By the step 10 of Algorithm \ref{alg:4}, we have
\begin{align} \label{eq:Q2}
\mathbb{E}\big[\mbox{dist}(0,\nabla_x L(x,y_{[m]},\lambda))^2\big]_{k+1} & = \mathbb{E}\|A^T\lambda_{k+1}-\nabla f(x_{k+1})\|^2  \nonumber \\
& = \mathbb{E}\|v_k - \nabla f(x_{k+1}) - \frac{G}{\eta} (x_k - x_{k+1})\|^2 \nonumber \\
& = \mathbb{E}\|v_k - \nabla f(x_{k}) + \nabla f(x_{k})- \nabla f(x_{k+1}) - \frac{G}{\eta}(x_k-x_{k+1})\|^2  \nonumber \\
& \leq 3\mathbb{E}\|v_k - \nabla f(x_{k})\|^2 +  3\mathbb{E}\|\nabla f(x_{k})- \nabla f(x_{k+1})\|^2 + \frac{3\sigma^2_{\max}(G)}{\eta^2}\|x_k-x_{k+1}\|^2  \nonumber \\
& \leq \frac{18dL^2}{b_2}\sum_{i=c_kq}^{k-1}\mathbb{E}\|x_{i+1}-x_i\|^2 + \frac{9\nu^2 L^2 d^2}{2} + \frac{9qL^2d^2\nu^2}{b_2}
+ \frac{24\sigma^2}{b_1} + 24L^2d\mu^2 \nonumber \\
& \quad + 3(L^2+ \frac{\sigma^2_{\max}(G)}{\eta^2})\|x_k-x_{k+1}\|^2 \nonumber \\
& \leq 18(L^2+ \frac{\sigma^2_{\max}(G)}{\eta^2})\theta_{k} + \frac{9\nu^2 L^2 d^2}{2} + \frac{9qL^2d^2\nu^2}{b_2}
+ \frac{24\sigma^2}{b_1} + 24L^2d\mu^2,
\end{align}
where the second inequality holds by Lemma \ref{lem:P1}, and the last inequality is due to the definition of $\theta_k$ and $b_2=qd$.

By the step 11 of Algorithm \ref{alg:4}, we have
\begin{align} \label{eq:Q3}
\mathbb{E}\big[\mbox{dist}(0,\nabla_\lambda L(x,y_{[m]},\lambda))^2\big]_{k+1} & = \mathbb{E}\|Ax_{k+1}+ \sum_{j=1}^m B_jy_j^{k+1}-c\|^2 \nonumber \\
&= \frac{1}{\rho^2} \mathbb{E} \|\lambda_{k+1}-\lambda_k\|^2  \nonumber \\
& \leq \frac{60dL^2}{\rho^2\sigma^A_{\min}b_2}\sum_{i=c_kq}^{k-1}\mathbb{E}\|x_{i+1}-x_i\|^2 + \big(\frac{5\sigma^2_{\max}(G)}{\rho^2\sigma^A_{\min}\eta^2}+\frac{5L^2}{\rho^2\sigma^A_{\min}}\big)\mathbb{E}\|x_k-x_{k-1}\|^2\nonumber \\
& \quad + \frac{5\sigma^2_{\max}(G)}{\rho^2\sigma^A_{\min}\eta^2}\|x_{k+1}-x_{k}\|^2 + \frac{15\nu^2 L^2 d^2}{\sigma^A_{\min}} + \frac{30qL^2d^2\nu^2}{b_2\sigma^A_{\min}}  + \frac{80\sigma^2}{b_1\sigma^A_{\min}} + \frac{80L^2d\mu^2}{\sigma^A_{\min}} \nonumber \\
& \leq \big( \frac{60L^2}{\sigma^A_{\min} \rho^2} + \frac{5\sigma^2_{\max}(G) }{\sigma^A_{\min}\eta^2\rho^2} \big) \theta_{k}
+ \frac{15\nu^2 L^2 d^2}{\sigma^A_{\min}\rho^2} + \frac{30qL^2d^2\nu^2}{b_2\sigma^A_{\min}\rho^2}  + \frac{80\sigma^2}{b_1\sigma^A_{\min}\rho^2} + \frac{80L^2d\mu^2}{\sigma^A_{\min}\rho^2},
\end{align}
where the first inequality follows by Lemma \ref{lem:P2}, and the last inequality is due to the definition of $\theta_k$ and $b_2=qd$.

Let $\beta_{\max}=\max\{\beta_1,\beta_2,\beta_3\}$ with
\begin{align}
  \beta_1 = m\big(\rho^2\sigma^B_{\max}\sigma^A_{\max} + \rho^2(\sigma^B_{\max})^2 + \sigma^2_{\max}(H)\big), \ \beta_2 = 18(L^2 + \frac{\sigma^2_{\max}(G)}{\eta^2}), \
  \beta_3 = \frac{60L^2 }{\sigma^A_{\min} \rho^2} + \frac{5\sigma^2_{\max}(G) }{\sigma^A_{\min}\eta^2\rho^2}. \nonumber
\end{align}
Combining the above inequalities \eqref{eq:Q1}, \eqref{eq:Q2} and \eqref{eq:Q3}, we have
\begin{align}
\frac{1}{K}\sum_{k=1}^K \mathbb{E}\big[ \mbox{dist}(0,\partial L(x_k,y_{[m]}^k,\lambda_k))^2\big] & \leq \frac{\beta_{\max}}{K}\sum_{k=1}^{K-1} \theta_k
 + \big(\frac{9}{2} + \frac{15}{\sigma^A_{\min}\rho^2}\big)(1+\frac{2q}{b_2})L^2d^2\nu^2 + \big(\frac{\sigma^2}{b_1} + L^2d\mu^2\big)(24 + \frac{80}{\sigma^A_{\min}\rho^2}) \nonumber \\
& \leq  \frac{3\beta_{\max}(\Psi_{0} - \Psi^*)}{K \gamma} + \frac{3\beta_{\max}\vartheta_1\vartheta_2d^2\nu^2}{\gamma} + \frac{3\beta_{\max}\vartheta_1}{\gamma}\big(\frac{4\sigma^2}{b_1} + 4L^2d\mu^2\big) \nonumber \\
& \quad + \big(\frac{9}{2} + \frac{15}{\sigma^A_{\min}\rho^2}\big)(1+\frac{2q}{b_2})L^2d^2\nu^2 + \big(\frac{\sigma^2}{b_1} + L^2d\mu^2\big)(24 + \frac{80}{\sigma^A_{\min}\rho^2}),
\end{align}
where $\vartheta_1 = \frac{1}{L} + \frac{20}{\sigma^A_{\min}\rho}$, $\vartheta_2=\frac{3L^2}{4}+\frac{3qL^2}{2b_2}$,
and the second inequality holds by the above Lemma \ref{lem:P3}.

Let $b_2=dq$, $\eta = \frac{\alpha\sigma_{\min}(G)}{4L} \ (0<\alpha \leq 1)$ and $\rho = \frac{2\sqrt{237}\kappa_G L}{\sigma^A_{\min}\alpha}$. Since $m$ denotes the number of nonsmooth regularization functions, it is relatively small.
It is easily verified that
$\vartheta_1=O(1)$, $\vartheta_2=O(1)$, $\beta_{\max} = O(1)$ and $\gamma=O(1)$, which are independent on $n$ and $K$. Thus, we obtain
\begin{align}
\mathbb{E}\big[ \mbox{dist}(0,\partial L(x_\zeta,y_{[m]}^\zeta,\lambda_\zeta))^2\big] = \frac{1}{K}\sum_{k=1}^K \mathbb{E}\big[ \mbox{dist}(0,\partial L(x_k,y_{[m]}^k,\lambda_k))^2\big] \leq  O(\frac{1}{K})+ O(d^2\nu^2) + O(d\mu^2) + O(\frac{1}{b_1}),
\end{align}
where $\{x_\zeta,y_{[m]}^{\zeta},\lambda_\zeta\}$ is chosen uniformly randomly from $\{x_{k},y_{[m]}^{k},\lambda_k\}_{k=0}^{K-1}$.

\end{proof}

\end{appendices}

\end{onecolumn}

\end{document}